%%%%%%%%%%%%%%%%%%%%%%%%%%%%%%%%%%%%%%%%%%%%%%%%%%%%%%%%%%%%%%%%%%%%%%%%%%%
%% Abhyankar, Shreeram S.
%% 
%% Galois theory on the line in nonzero characteristic
%% 
%% publ:  Bull. Amer. Math. Soc. (N.S.) 27(1992) no. 1
%% pp:    68-133
%% type:  Research-Expository Paper    markup: amstex    file size: 233K
%% 
%% copyright: American Math. Society copyright; see end of article
%% 
%% Include files necessary for this article: bull-ppt.tex
%% 
%%%%%%%%%%%%%%%%%%%%%%%%%%%%%%%%%%%%%%%%%%%%%%%%%%%%%%%%%%%%%%%%%%%%%%%%%%%
% Date: 22-AUG-1991
%   : 0   \newpage: 0   \displaybreak: 0
%   \eject: 0   \bye: 0   \break: 0   \allowbreak: 0
%   \allowdisplaybreak: 0   \allowdisplaybreaks: 0
%   \allowlinebreak: 0   \allowmathbreak: 0
%   \smallpagebreak: 0   \medpagebreak: 0   \bigpagebreak: 0
%   \smallbreak: 0   \medbreak: 0   \bigbreak: 0   
%\goodbreak: 0
%   : 0   : 0   \newline: 29
%   \magnification: 1   \mag: 0
%   \baselineskip: 0   \normalbaselineskip: 0
%   \hsize: 0   \vsize: 0   \pagewidth: 0   \pageheight: 0
%   \hoffset: 0   \voffset: 0   \hcorrection: 0   
%\vcorrection: 0
%   \parindent: 0   \parskip: 0
%   \vfil: 0   \vfill: 0   \vskip: 0
%   \smallskip: 0   \medskip: 0   \bigskip: 0
%   \sl: 2   \def: 86   \let: 0   \redefine: 0   
%\predefine: 0
%   \tolerance: 0   \pretolerance: 0
%   ASCII 13 (Control-M Carriage return): 0
%   ASCII 10 (Control-J Linefeed): 0
%   ASCII 12 (Control-L Formfeed): 0
%   ASCII 0 (Control-@): 0
%
%\input usefraktur.tex
\input amstex
\documentstyle{amsppt}
\input bull-ppt
\NoBlackBoxes
\keyedby{bull270e/paz}

\define\ds{\displaystyle}
\define\Aut{\text{\rm Aut}\,}
\define\CHAR{\text{\rm char}\,}
\redefine\deg{\text{\rm deg}\,}
%\define\det{\text{\rm det}\,}

\define\im{\text{\rm im}\,}

%\define\ker{\text{\rm ker}\,}

\define\cf{{\Cal F}}

\define\cp{{\Cal P}}
%\define\cr{{\Cal R}}

\define\a{\alpha}
\redefine\b{\beta}
%\redefine\c{\psi}
\define\uc{\Psi}
\redefine\d{\delta}
\define\ud{\Delta}
\define\e{\epsilon}
\define\f{\phi}
\define\uf{\Phi}
\define\g{\gamma}
\define\ug{\Gamma}
\define\h{\theta}
\define\uh{\Theta}
\define\k{\kappa}
\redefine\l{\lambda}
\redefine\ul{\Lambda}
\define\m{\mu}
\define\n{\nu}
\define\p{\pi}
\define\r{\rho}
\define\s{\sigma}
\redefine\t{\tau}
\redefine\v{\chi}
\define\w{\omega}
\define\uw{\Omega}
\define\x{\xi}
\define\y{\eta}
\define\z{\zeta}

\define\bc{{\Bbb C}}

\define\bq{{\Bbb Q}}

\define\od{\bar d}
\define\of{\bar f}
\define\ok{\bar k}
\define\OR{\bar r}
\define\oF{\overline{F}}
\define\oG{\overline{G}}
\define\oK{\overline{K}}
\define\tF{\widetilde{F}}
\define\tG{\widetilde{G}}

\define\hv{\hat v}
\define\hA{\widehat{A}}
\define\hF{\widehat{F}}
\define\hG{\widehat{G}}
\define\hK{\widehat{K}}
\define\hM{\widehat{M}}

\topmatter
\cvol{27}
\cvolyear{1992}
\cmonth{July}
\cyear{1992}
\cvolno{1}
\cpgs{68-133}
%\centerline{{\bf 
\title Galois theory on the line in nonzero\\
characteristic\endtitle
%\centerline{
\author Shreeram S. Abhyankar\endauthor
\shortauthor{S. S. Abhyankar}
%\footnote""{
%1980 Mathematical Subject Classification (1985 Revision). 
\subjclass Primary 12F10, 14H30, 20D06, 20E22\endsubjclass
%\newline\indent
\thanks This work was partly supported by NSF grant DMS 
88-16286\endthanks
\date December 19, 1990, and, in revised form, on July 31, 
1991\enddate
%\centerline{
\address {\rm{(Shreeram S. Abhyankar)}} Mathematics 
Department, Purdue 
University, West Lafayette, Indiana 47907\endaddress
%\centerline{}
%\centerline{\it Dedicated to}
\dedicto Walter Feit, J-P. Serre, and e-mail\ededicto
%\centerline{}
%\centerline{}
\endtopmatter

\document
%\centerline{{\bf }}
\heading 1. What is Galois theory?\endheading

Originally, the equation $Y^2 + 1 = 0$ had no solution. 
Then the two solutions 
$i$ and $-i$ were created.  But there is absolutely no way 
to tell 
who is $i$ and who is $-i$.  
\footnote{Only the physicist can tell the difference by 
declaring that
$i$ is up and $-i$ is down. But then, what is up and what 
is down?
%}
%\footnote{
%\newline\indent
What I am saying is that intrinsically there is no way to 
distinguish,
in an abstractly given copy of the complex numbers, 
between the two 
square-roots of minus one. However, practically everyone 
has their own
favorite concrete model of $\Bbb C$, perhaps $\Bbb 
R\times\Bbb R$, perhaps
$\Bbb R[Y]/(Y^2+1)$, and then there is no problem in 
pointing to say
$(0,1)$ in the first case or $[Y]$ in the second and 
calling that $i$.}
That is Galois Theory.

Thus, Galois Theory tells you how far we cannot 
distinguish between
the roots of an equation. This is codified in the Galois 
Group. 

%\centerline{}
%\centerline{{\bf }}
\heading 2. Galois groups\endheading

More precisely, consider an equation 
$$
Y^n + a_1Y^{n-1} + \dots + a_n = 0
$$
and let $\a_1,\dots,\a_n$ be its roots, which are assumed 
to be
distinct.  By definition, the Galois Group $G$ of this 
equation
consists of those permutations of the roots which preserve 
all
relations between them. Equivalently, $G$ is the set of 
all those permutations
$\s$ of the symbols $\{1,2,\dots,n\}$ such that 
$\f(\a_{\s(1)},\dots,\a_{\s(n)})=0$ for every $n$-variable 
polynomial $\f$
for which $\f(\a_1,\dots,\a_n)=0$. The coefficients of 
$\f$ are supposed to be 
in a field $K$ which contains the coefficients 
$a_1,\dots,a_n$ of the 
given polynomial 
$$
f=f(Y)=Y^n+a_1Y^{n-1}+\dots+a_n.
$$ 
We call $G$ the Galois Group of $f$ over $K$ and denote it 
by 
$\operatorname{Gal}_Y(f,K)$ or 
$\operatorname{Gal}(f,K)$. This is Galois' original 
concrete definition. 

According to the modern abstract definition, the Galois 
Group of a 
normal extension $L$ of a field $K$ is defined to be the 
group of
all $K$-automorphisms of $L$ and is denoted by 
$\operatorname{Gal}(L,K)$. Note that a normal
extension $L$ of a field $K$ is a field obtained by 
adjoining to $K$ all
\pagebreak
the roots of a bunch of univariate polynomials with 
coefficients in $K$.
To relate the two definitions, let $L=K(\a_1,\dots,\a_n)$ 
and note that we get
an isomorphism of $\operatorname{Gal}
(L,K)$ onto $\operatorname{Gal}(f,K)$ by sending any 
$\t\in\operatorname{Gal}(L,K)$ to that 
$\s\in\operatorname{Gal}(f,K)$ for which 
$\t(\a_i)=\a_{\s(i)}$ for $1\leq i\leq n$.

%\centerline{{\bf }}
\heading 3. Permutation groups\endheading

The above concrete definition brings out the close 
connection
between group theory and theory of equations. To wit, the 
Galois Group 
$\operatorname{Gal}
(f,K)$ is now a subgroup of $S_n$ where, as usual, $S_n$ 
denotes the
symmetric group of degree $n$, i.e., the group of all 
permutations of
$n$ symbols; note that the order of $S_n$ is $n!$. Quite 
generally,
a subgroup of $S_n$ is called a (permutation) group of 
degree $n$.
Here the use of the word degree is meant to remind us that 
potentially it
comes from an equation of degree $n$. To convert this 
potentiality into
actuality, in various situations, constitutes Inverse 
Galois Theory.	
To further bring out the parallelism between group theory 
and the theory of
equations, we note that.

(1) $f$ is irreducible iff $\operatorname{Gal}(f,K)$ is 
transitive.

\noindent
Here, a permutation group $G$ ``acting'' 
\footnote{Self-advice: Don't be so scared of the term
``acting.'' It is simply the modern substitute for 
%%\newline 
``permuting.''}
on the set $\uw=\{1,2,\dots,n\}$ is {\it transitive} if 
for all $i,j$ in $\uw$,
there exists $\s\in G$ such that $\s(i)=j$. Likewise, $G$ 
is 2-transitive 
(or doubly transitive) if for all $i\neq i'$ and $j\neq 
j'$ in $\uw$, there
exists $\s\in G$ such that $\s(i)=j$ and $\s(i')=j'$. Quite
generally, $G$ is $l$-{\it transitive} for a positive 
integer $l\leq n$, if
for all pairwise distinct elements $i_1,i_2,\dots,i_l$ in 
$\uw$
and pairwise distinct elements $j_1,j_2,\dots,j_l$ in 
$\uw$, there
exists $\s\in G$ such that $\s(i_e)=j_e$ for $1\leq e\leq 
l$.
This brings us to MTR, i.e., the method of

%\centerline{}
%\centerline{{\bf }}
\heading 4. Throwing away roots\endheading 

Assuming $f$ to be irreducible in $K[Y]$, let us ``throw 
away'' a root of
$f$, say $\a_1$, and get
$$
f_1=f_1(Y)=\frac{f(Y)}{(Y-\a_1)}
= Y^{n-1}+b_1Y^{n-2}+\dots+b_{n-1}\in K(\a_1)[Y].
$$
In continuation of (1), we see that.

(2) $f$ and $f_1$ are irreducible in $K[Y]$ and 
$K(\a_1)[Y]$ respectively 
iff$\operatorname{Gal}(f,K)$ is 
%\newline 
2-transitive.

It may be noted that, assuming $f$ to be irreducible, it 
does not
matter which root of $f$ we throw away; for instance, the 
irreducibility
of $f_1$ in $K(\a_1)[Y]$ and, up to isomorphism, the 
Galois group
$\operatorname{Gal}(f_1,K(\a_1))$ are independent of which 
root we call $\a_1$.

Likewise, by throwing away $s$ roots of $f_0=f$ we get
$$
\split
f_s=f_s(Y)&=\frac{f(Y)}{(Y-\a_1)\cdots(Y-\a_s)}\\
&=Y^{n-s}+d_1Y^{n-s-1}+\dots+d_{n-s}\in 
K(\a_1,\dots,\a_s)[Y]
\endsplit
$$
and then:

(3) $f_s$ is irreducible in $K(\a_1,\dots,\a_s)[Y]$ for 
$0\leq s<l$ iff 
$\operatorname{Gal}(f,K)$ is $l$-transitive.

%\centerline{}
%\centerline{{\bf }
\heading 5. Classification theorems\endheading 

Now you would have thought that you could (easily or 
possibly)
construct a 
polynomial $f=f_0$, say of degree 20, such that 
$f_0,f_1,\dots,f_9$ are
irreducible whereas $f_{10}$ is reducible. But No! And 
that is the
surprise of the century. You cannot!! So says the CT, 
i.e., the
recently established Classification Theorem of Finite 
Simple Groups,
which was a magnificent piece of ``team work''.  According 
to the staggering 
statistics as reported by Coach Gorenstein \cite{G3},
the CT took 30 years (1950--1980), 100 authors, 500 
papers, and 15000 
pages!
Yet several hundred more pages are required to prove the 
implications 
$$ 
\text{CT}\Rightarrow\text{CDT}\Rightarrow\text{CTT}%
\Rightarrow 
\text{CQT}\Rightarrow\text{CFT}\Rightarrow\text{CST} 
$$ 
where CDT (resp: CTT, CQT, CFT, and CST) stands for the 
Classification 
Theorem of Doubly (resp: Triply, Quadruply, Fivefold, and 
Sixfold)
transitive permutation groups. 
\footnote{This is an expanded version of a lecture given 
at Walter Feit's 
60$^{\text{th}}$
birthday conference in Oxford, England. In addition to 
Walter Feit, the audience
included the group theorists Peter Cameron, Michael 
Collins, Sandy Green, 
Graham Higman, Peter Neumann, Ron Solomon, and John 
Thompson; 
my talking group theory in this meet of topnotch group 
theorists was 
like carrying coal to Newcastle, or bringing holy water to 
the Ganges!!} 
Promising to come back to CDT to CFT
in a moment, let us state CST. It simply says that the 
symmetric group $S_n$ 
for $n\geq 6$ and the alternating group
$A_n$ for $n\geq 8$ are the only sixfold transitive 
groups!! 
Vis-a-vis equations, what we are saying is that if $f$ is 
an irreducible 
polynomial of degree $n>7$ such that $f_1,f_2,f_3,f_4$ and 
$f_5$ are irreducible
then so are $f_6,f_7,\dots,f_{n-3}$.
\footnote{Self-Challenge = challenge to the extollers of 
high school
algebra: prove that by high school algebra if you can! Of 
course we
can simply decree CT be high school algebra!!} 

To take a first shot at CT, in addition to permutation 
groups, we should also
consider groups of matrices over the (Galois) Field $
\operatorname{GF}(q)$ of $q$ elements 
where $q$ is a power of a prime. So let 
$$
\align
\operatorname{GL}(m,q)=&\text{ the {\it general linear 
group} of degree }m
\text{ over }\operatorname{GF}(q) \\ 
=&\text{the group of all nonsingular }m\text{ by }m
\text{ matrices}\\
&\text{with entries in }\operatorname{GF}(q).
\endalign
$$
Here the multiplicative group $\operatorname{GF}
(q)^*$ comes in two ways. Firstly,
thinking of scalar matrices, $\operatorname{GF}
(q)^*$ becomes a normal subgroup
(and, in fact, the center) of $\operatorname{GL}
(m,q)$. Secondly, taking determinants 
we get a surjective homomorphism $\operatorname{GL}
(m,q)\rightarrow \operatorname{GF}(q)^*$. This 
motivates the definitions 
$$
\align
\operatorname{PGL}
(m,q)=&\text{the {\it projective general linear group} of 
degree }m\\
&\text{over } \operatorname{GF}(q) \\
=&\operatorname{GL}(m,q)/\operatorname{GF}(q)^*
\endalign
$$
and 
$$
\align
\operatorname{SL}
(m,q)=&\text{ the {\it special linear group} of degree 
}m\text{ over } 
\operatorname{GF}(q) \\
=&\text{ ker }\operatorname{GL}(m,q)\rightarrow 
\operatorname{GF}(q)^*.
\endalign
$$
Combining these two roles of $\operatorname{GF}(q)^*$ we get
$$
\align
\operatorname{PSL}(m,q)=&\text{ the {\it projective 
special linear group} of
degree }m\\
&\text{over } \operatorname{GF}(q)\\
=&\operatorname{SL}(m,q)/(\operatorname{SL}(m,q)\cap 
\operatorname{GF}(q)^*).
\endalign
$$
In group theory parlance
\footnote{In learning group theory, I am following the 
traditional
Indian method: memorize things by heart and the meaning will
eventually be revealed to you. Moreover, every subject has 
its lingo. Thus 
$\operatorname{GL}$, $\operatorname{SL}$ and 
$\operatorname{PSL}$ 
are the Tom, Dick, and Harry of group theory.}
$$
\align
L_m(q)=&\text{ the {\it linear group} of degree }m\text{ 
over }\operatorname{GF}(q) \\
=&\operatorname{PSL}(m,q).
\endalign
$$
In yet another notation
$$
A_n(q) = L_{n+1}(q)=\operatorname{PSL}(n+1,q).
$$

With these preliminaries,
 
{\it CT essentially says that $A_n$ and $A_n(q)$
are the only \RM(finite\/\RM) simple groups.\/}
 
Here we have to exclude ``small'' cases; namely, from the 
alternating group 
$A_n$ exclude
$n\leq 4$, and from $A_n(q)$ exclude $n=1$ and $q\leq 3$ 
(we define $A_n(q)$
only for $n\geq 1$). Moreover, ``essentially'' means that 
with $A_n(q)$
we have to include its relatives and incarnations, to be 
discussed later.
Finally, in addition to these infinite families, there are 
26 ``sporadics'',
again to be discussed later.
\footnote{Does the number 26 vindicate the spread of the 
English language
which has exactly that many letters?}

%\centerline{}
%\centerline{{\bf }}
\heading 6. Brief thirty year history\endheading

First, in the fifties, there was the fundamental work of 
Brauer \cite{B} and 
Chevalley \cite{Ch}. In 1962 this was followed by the 
path-breaking odd order
paper of Feit and Thompson \cite{FT}. Then came the large 
team coached by
Gorenstein \cite{G4}. At any rate, we are meeting here to 
felicitate our
friend Walter Feit on the occasion of his forthcoming 
sixtieth
birthday.
\footnote{My fondest memory of Walter is that in 1957, 
when we were both 
at Cornell, we decided to go on a diet together and the 
one who lost more
weight was to get a quarter. At the end of one month, 
Walter gained
one pound and I gained two. Who lost and who won?}

%\centerline{}
%\centerline{{\bf }}
\heading 7. Primitive groups\endheading

As another example of the parallelism between group theory 
and
theory of equations, let us note that a permutation group 
$G$ is
said to be primitive if it is transitive, and the  {\it 
one-point
stabilizer} $G_1$ of $G$ is a maximal subgroup of $G$. 
Here we are assuming $G$ 
to be a subgroup of the symmetric group $S_n$ acting on 
$\{1,2,\dots,n\}$, and
then by definition, $G_1= G\cap S_{n-1}$ with 
$S_{n-1}=\{\s\in S_n:\s(1)=1\}$.
Just as it did not matter which root of the irreducible 
equation we threw
away, so in the present situation, if $G$ is transitive 
then we may replace
$G_1$ by any $G_i=\{\s\in G:\s(i)=i\}$, which is called 
the {\it stabilizer} of
$i$ in $G$. Clearly.

(4) If $f$ is irreducible and $G=\operatorname{Gal}
(f,K)$ then $G_1=\operatorname{Gal}(f_1,K(\a_1))$.

Moreover;

(5) 
$\operatorname{Gal}(f,K)$ is primitive iff $f$ is 
irreducible and there is no field 
between $K$ and $K(\a_1)$. 

Not to get completely lost in group theory, let us revert 
to algebraic 
geometry by talking about 

%\centerline{}
%\centerline{{\bf }}
\heading 8. Fundamental groups\endheading

In my 1957 paper on ``Coverings of Algebraic Curves'' in 
the American
Journal \cite{A3}, I considered the {\it algebraic 
fundamental group}
$\p_A(C)$ of a nonsingular curve $C$. Here $C$ is allowed 
to be ``open'', 
i.e., it may consist of a projective algebraic curve minus 
a finite number of 
points. We assume that $C$ is irreducible and defined over 
an algebraically 
closed ground field $k$. By $\p_A(C)$ we mean the family 
of Galois groups 
$\operatorname{Gal}
(L/k(C))$ as $L$ varies over all finite normal extensions 
of the function 
field $k(C)$ of $C$ such that no point of $C$ 
(equivalently, no valuation of 
$k(C)/k$ having center on $C$) is {\it ramified} in $L$.

In case $C$ is the (affine) line $L_k$ over $k$, or more 
generally if 
$C=L_{k,r}=$ the line $L_k$ minus $r$ points 
$\l_1,\dots,\l_r$ 
(with $\l_i\neq\l_j$ in $k$ for $1\leq i < j\leq r$) then 
this amounts to 
considering 
$\operatorname{Gal}(F,k(X))$ where
$$
F=F(X,Y)=Y^n+\f_1(X)Y^{n-1}+\dots+\f_n(X)
$$
is a bivariate polynomial with coefficients 
$\f_1(X),\dots,\f_n(X)$ in $k[X]$
such that $F$ is unramified at all (finite) values of $X$ 
other than
$\l_1,\dots,\l_r$, i.e., such that for every $\l$ in $k$ 
different from
$\l_1,\dots,\l_r$ we have
$$
F(X+\l,Y)=\prod_{i=1}^n(Y-\y_{\l}^{(i)}(X))
$$
with $\y_{\l}^{(i)}(X)$ in the (formal) power series ring 
$k[[X]]$.
We may call $F$ an {\it unramified covering} of $L_{k,r}$.
  
For any group $G$\<, let $G_h$ denote the family of finite 
homomorphic images of
$G$. Let $\cf_r$ be the free group on $r$ generators, and 
let $J_r$ be the
family of all finite groups generated by $r$ generators, 
and note that then
$J_r=\cf_{rh}$. By the Riemann Existence Theorem etc., we 
see that
if $k$ = the field of complex numbers $\bc$, then for any 
(irreducible) nonsingular curve $C$ over $k$ we have
$\p_A(C)=\p_1(C)_h$ where, as usual, $\p_1(C)$ denotes the 
(topological) fundamental group of $C$; see Serre 
\cite{S1}. Hence in particular
$\p_A(L_{\bc,r})=\p_1(L_{\bc,r})_h$, and clearly 
$\p_1(L_{\bc,r})=\cf_r$.
Therefore $\p_A(L_{\bc,r})=J_r$.

By taking $r=0$ or $1$ in the last equation, we get 
$\p_A(L_{\bc})=J_0=\{1\}$ and 
$\p_A(L_{\bc,1})=J_1=$ the family of all finite cyclic 
groups. 
These two facts can also be proved purely algebraically, 
by the genus formula 
due to Hurwitz--Riemann--Zeuthen; see Serre \cite{S2}. 
The said genus formula actually shows that
$$
\p_A^*(L_k)=J_0=\{1\}
$$ 
and 
$$
\align
\p_A^*(L_{k,1})=&J_1[\text{char }k]^* \\
             =&\text{ the family of all finite cyclic 
groups with order}\\
              &\text{ nondivisible by char }k,
\endalign
$$
where char $k$ is the characteristic of $k$ 
and where for any (irreducible) nonsingular curve $C$ over 
any 
algebraically closed field $k$ we define
$$
\align
\p_A^*(C) =&\text{ the {\it globally tame fundamental 
group} of }C \\
          =&\text{ the family of all the members of }
            \p_A(C)\text{ whose order is}\\
           &\text{ nondivisible by char }k
\endalign
$$
and for any nonnegative integer $m$ and any family of 
finite groups
$J$ we put
$$
J[m]^*=\text{ the family of all the members of } J
\text{ whose order is nondivisible by } m.
$$

%\centerline{}
%\centerline{{\bf }}
\heading 9. Quasi $p$-groups\endheading

Given any prime number $p$ we put
$$
Q(p)=\text{ the family of all quasi }p\text{-groups},
$$
where by a {\it quasi} $p$-{\it group} 
we mean a finite group which is generated by
all of its $p$-Sylow subgroups, and for any finite group 
$G$ we put
$$
p(G) =\text{ the (normal) subgroup of } G\text{ generated 
by all of its }
p\text{-Sylow subgroups}
$$
and for any family of finite groups $J$ we put
$$
J(p)=\text{ the family of all finite groups }G\text{ such 
that }G/p(G)\in J
$$
and we note that then
$J_0(p)=Q(p)$
and for every nonnegative integer $r$ we have
$$
\align
J_r(p) =&J_r[p]^*(p) \\ 
=&\text{ the family of all finite groups }G \text{ such 
that }
                    G/p(G)\text{ is generated}\\
 &\text{ by }r\text{ generators}
\endalign
$$
and
$$
\{p(G):G\in J_r(p)\}=Q(p).
$$

%\centerline{}
%\centerline{{\bf }}
\heading 10. Conjectures\endheading
 
For any algebraically closed field $k$ of characteristic 
$p\neq 0$, in the 
above cited 1957 paper, I conjectured that 
$\p_A(L_{k,r})=\p_A(L_{\bc,r})(p)$, i.e., equivalently, 
\dfn{General Conjecture} 
For every nonnegative integer $r$ we have 
$\p_A(L_{k,r})=J_r(p)$.
%\endproclaim
Hence in particular
\enddfn
\subheading{Quasigroup Conjecture} 
$\p_A(L_k)=Q(p)$.
%\endproclaim

Now a (finite) simple group whose order is divisible by 
$p$ is obviously a
quasi $p$-group, and therefore QC (= the Quasigroup 
Conjecture) subsumes the
\dfn{Simple Group Conjecture} $\p_A(L_k)$ contains every 
simple group 
whose order is divisible by $p$.
\enddfn
%\endproclaim
In particular
\subheading{Alternating Group Conjecture} For every 
integer $n\geq p$ we have
$A_n\in \p_A(L_k)$ except when $p=2<n<5$.
%\endproclaim

Here, as usual, $A_n$ denotes the alternating group of 
degree $n$, i.e.,
the group of all even permutations on $n$ symbols, and we 
note that the
order of $A_n$ is $n!/2$ or $1$ according as $n>1$ or 
$n=1$. 
As proved by Galois, $A_n$ is a simple group for every 
$n>4$,
and hence AGC (= the Alternating Group Conjecture) is a 
special case of SGC 
(= the Simple Group Conjecture) except when $(p,n)=$ 
(3,4); in this exceptional
case, $A_n$ is obviously a quasi $p$-group and so we can 
directly fall back 
upon QC. The symmetric
group $S_n$ of degree $n$ is obviously a quasi $2$-group 
for every $n$, and
hence QC subsumes the
\dfn{Even Prime Symmetric Group Conjecture} 
If $p=2$ then for every 
integer $n\geq 2$ we have $S_n\in\p_A(L_k)$.
\enddfn
%\endproclaim 

To match up with EPSGC (= the Even Prime Symmetric Group 
Conjecture), let
us divide AGC into EPAGC and OPAGC, i.e., into the 
following two
conjectures respectively.
\subheading{Even Prime Alternating Group Conjecture} If 
$p=2$ then for every 
integer $n\geq 2$ except when $2<n<5$ we have 
$A_n\in\p_A(L_k)$.
%\endproclaim
\subheading{Odd Prime Alternating Group Conjecture} If 
$p>2$ then for every 
integer $n\geq p$ we have $A_n\in\p_A(L_k)$.
%\endproclaim

%\centerline{}
%\centerline{{\bf }}
\heading 11. Again some history\endheading

Let $0\neq a\in k=$ an algebraically closed field of 
characteristic 
$p\neq 0$, and let $n,s,t$ be positive integers such that 
$t\not\equiv 0(p)$.
Now in support of the above conjectures, in the above 
cited 1957 paper, I 
had written down the following two examples of unramified 
coverings
of $L_k$:
$$
\hF_n=Y^n-aX^sY^t+1\quad\text{with }n=p+t
$$ 
and
$$
\tF_n=Y^n-aY^t+X^s\text{ with }t<n\equiv 0(p)\text{ and 
}\operatorname{GCD}(n,t)=1\text{ and }
s\equiv 0(t)
$$
and  had suggested that their Galois groups 
$\hG_n=\operatorname{Gal}(\hF_n,k(X))$
and $\tG_n=\operatorname{Gal}(\tF_n,k(X))$ should be 
calculated.

Now, after a gap of thirty years, with Serre's 
encouragement and with the help
of CT, I can calculate these Galois groups, and the 
answers are as follows.

%\centerline{{\bf (I)}}
\heading (I)\endheading

(I.1) If $t=1$, then $\hG_n=\operatorname{PSL}(2,p)=%
\operatorname{PSL}(2,n-1)$.

(I.2) If $t=2$ and $p=7$, then $\hG_n=%
\operatorname{PSL}(2,8)=\operatorname{PSL}(2,n-1)$.

(I.3) If $t=2$ and $p\neq 7$, then $\hG_n=A_n$.

(I.4) If $t>2$ and $p\neq 2$, then $\hG_n=A_n$.

(I.5) If $p=2$, then $\hG_n=S_n$.

%\centerline{{\bf (II)}}
\heading (II)\endheading

(II.1) If $1<t<4$ and $p\neq 2$, then $\tG_n=A_n$. 

(II.2) If $1<t<n-3$ and $p\neq 2$, then $\tG_n=A_n$. 

(II.3) If $1<t=n-3$ and $p\neq 2$ and $11\neq p \neq 23$, 
then $\tG_n=A_n$. 

(II.4) If $1<t<4<n$ and $p=2$, then $\tG_n=A_n$ or $S_n$.

(II.5) If $1<t<n-3$ and $p=2$, then $\tG_n=A_n$ or $S_n$.

Actually, $\hF_n$ is a slight generalization of the 
original equation
$$
\oF_n=Y^n-XY^t+1\quad\text{with }n=p+t
$$
written down in the 1957 paper. This equation $\oF_n$ was 
discovered by taking 
a section of a surface extracted from my 1955 paper 
``Ramification of 
Algebraic Functions'' in the American Journal \cite{A1}, 
which was the second 
part of my Ph.D. Thesis written under the able guidance of 
Oscar Zariski. In 
the 1956 paper in the Annals of Mathematics \cite{A2}, 
which  was the first 
part of my Ph.D Thesis, I proved resolution of 
singularities of algebraic 
surfaces in nonzero characteristic. In the 1955 American 
Journal paper, 
I was showing why 
Jung's method of surface resolution in the complex case 
does not generalize
to nonzero characteristic, because the local fundamental 
group, above a normal 
crossing of the branch locus, in the former case is abelian 
whereas in the latter case it can even
be unsolvable. It was a surface constructed for this 
purpose whose section
I took in the 1957 paper.

Although the second equation $\tF_n$ is also a slight 
generalization of an 
equation occuring in the 1957 paper but, amusingly,
it got rediscovered in 1989 as a variation of the first 
equation $\hF_n$. 

Now (I.1) was originally proved by Serre and when he told 
me about it in
September 1988, that is what started off my calculations 
after a thirty
year freeze! As a slight generalization of (I.1), in the 
case of $t=1$, I can 
also calculate the Galois group 
$\hG_{n,q}=\operatorname{Gal}(\hF_{n,q},k(X))$ of 
the unramified covering of $L_k$ given by 
$$
\hF_{n,q}=Y^n-aX^{-s}Y^t+1\quad\text{with }n=q+t,
$$
where $q$ is any positive power of $p$, and it turns out 
that

%\centerline{{\bf (III)}}
\heading (III)\endheading

(III.1) If $t=1$, then $\hG_{n,q}=\operatorname{PSL}(2,q)=%
\operatorname{PSL}(2,n-1)$.

Again, $\hF_{n,q}$ is a slight generalization of the 
equation
$$
\oF_{n,q}=Y^n-XY^t+1\quad\text{with }n=q+t,
$$
which also occurs in the 1957 paper. Note that then
$$
\oF_n=\oF_{n,p}.
$$

%\centerline{}
%\centerline{{\bf }}
\heading 12. Using MRT\endheading

Applying MRT (= method of removing tame ramification 
through cyclic
compositums = so called Abhyankar's Lemma) to the 
one-point stabilizer of 
$\hF_{n+1}$ we get the monic polynomial of degree $n$ in 
$Y$ with
coefficients in $k(X)$ given by
$$
\hF'_n=h(Y)(Y+b)^p-aX^{-s}Y^t\quad\text{with } 0\neq b\in k,
$$
where $h(Y)$ is the monic polynomial of degree $n-p$ in 
$Y$ with coefficients 
in $k$ given by
$$
h(Y)=\frac{(Y+n+1)^{n+1-p}-Y^{n+1-p}}{(n+1)^2}
$$
and we let $\hG'_n=\operatorname{Gal}(\hF'_n,k(X))$. 
As an immediate consequence of (I) we now get

%\centerline{{\bf (IV)}}
\heading (IV)\endheading

Assuming that $n+1\not\equiv 0(p)$, in the following cases
$\hF'_n$ gives an unramified covering of $L_k$ with the 
indicated Galois group.

(IV.1) If $n+1-p=t>2\neq p$ and $b=t$ and $s\equiv 0(p-1)$ 
and $s\equiv 0(t)$, 
then $\hG'_n=A_n$.

(IV.2) If $n+1-p=t=2$ and $p\neq 7$ and $b=t$ and $s\equiv 
0(p-1)$, 
then $\hG'_n=A_n$.

(IV.3) If $n=p+1$ and $p>5$, then $t$ can be chosen so 
that $1<t<\frac{p+1}{2}$
and $\operatorname{GCD}(p+1,t)=1$, and for any such $t$, 
upon assuming
$b=\frac{t}{t-1}$ and $s\equiv 0(t(p+1-t))$, we have 
$\hG'_n=A_n$. 

(IV.4) If $n+1-p=t$ and $p=2$ and $b=t$ and $s\equiv 0(t)$, 
then $\hG'_n=S_n$.

%\centerline{}
%\centerline{{\bf }}
\heading 13. Unramified coverings\endheading

We have the following four corollaries of calculations (I) 
through (IV).

In Calculations (I) through (IV), we use a lot of 
Ramification Theory, or,
equivalently, CS (= Cycle Structure). In addition to CS + 
MTR + MRT,
in the proofs of Calculations (I) to (IV) we also use CT.
In our original version of these proofs, the use of CT was 
heavy. 
Gradually the use of CT decreased, but could not be 
removed completely.
However, by traversing a delicate path through 
calculations (I) through (IV),
we have arranged a proof of the First and the Second 
Corollaries independent
of CT.

\proclaim{First Corollary} OPAGC is true. Equivalently, 
for any $n\geq p >2$,
there exists an unramified covering of the affine line in 
characteristic
$p$ whose Galois group is the alternating group $A_n$ of 
degree $n$.
\endproclaim
\proclaim{Second Corollary} EPSGC is true. Equivalently, 
for any $n\geq p =2$,
there exists an unramified covering of the affine line in 
characteristic
$p$ whose Galois group is the symmetric group $S_n$ of 
degree $n$.
\endproclaim
\proclaim{Third Corollary} Unramified coverings of the 
affine line in 
characteristic $p$ with a few more Galois groups have been 
constructed.
\endproclaim 
\dfn{Definition-Remark} By the {\it minimal index\/} 
of a finite group $G$ we mean the 
smallest number $d$ such that $G$ has a subgroup of index 
$d$ which does not
contain any nonidentity normal subgroup of $G$. 
\footnote{In other words, we are minimizing the index over 
subgroups of $G$ 
which do not contain any minimal normal subgroups of 
$G$\<, where we recall 
that a {\it minimal normal subgroup\/} of a group $G$ is a 
nonidentity 
normal subgroup $N$ of $G$ such that $N$ does not contain 
any nonidentity
normal subgroup of $G$ other than $N$ itself.} 
Now QC is obviously equivalent to saying that for every 
integer 
$d\geq p$ the following is true.
\enddfn

\proclaim{QC($d$)} Every quasi $p$-group of minimal index 
$d$ belongs to
$\p_A(L_k)$, i.e., occurs as the Galois group of an 
unramified covering of
the affine line in characteristic $p$.
\endproclaim
 
Therefore it is interesting to point out that as a 
corollary of the above 
results we have the following:

\proclaim{Fourth Corollary} QC$(p+1)$ is true for every 
$p$ which is not a
Mersenne prime 
\footnote{A {\it Mersenne prime} is a prime $p$ 
of the form $p=2^{\m}-1$ for 
some positive integer $\m$. Note that then $\m$ is 
necessarily prime,
because otherwise by factoring
$\m=\m'\m''$ with $\m'>1$ and $\m''>1$, we would get a 
factorization
$2^{\m}-1=l'l''$ with $l'=2^{\m'}-1>1$ and 
$l''=1+2^{\m'}+2^{2\m'}+\dots+2^{\m'(\m''-1)}>1$.}
and which is different from $11$ and $23$.
More generally, if $G$ is a quasi $p$-group containing a 
subgroup $H$ of
index $p+1$ such that $H$ does not contain any nonidentity 
normal subgroup
of $G$, and if $p$ is not a Mesenne prime and $p$ is 
different from $11$ and
$23$, then there exists an unramified covering of the 
affine line in
characteristic $p$ having $G$ as the Galois group.
\endproclaim

The proofs of the above four Corollaries and the four 
claims (I) to (IV)
will be completely given in this paper with the exception 
that the
proof of claim (I.2) will be completed in my forthcoming 
paper \cite{A7}.
In connection with (II.4), it may be noted that the case 
of $p=11$ or $23$
is still open.
\footnote{As will become apparent later, the reason for 
this, as well as
for the exclusion of these values of $p$ from the Fourth 
Corollary, 
is the existence of the ``Mathieu Groups''.}

%\centerline{}
%\centerline{{\bf }}
\heading 14. History of a pilgrimage\endheading

When I said that  ``Now...I can calculate these Galois 
groups'', what I really
meant was that, from September 1988 to August 1989, I 
undertook a
pilgrimage (physical as well as mental)
\footnote{mental = e-mail + s-mail. s-mail = snail mail = 
usual mail.}
to seek the help of lots of mathematicians, and then I 
simply collated the
help so obtained. In chronological order, these 
mathematicians were.
Serre (oh yes, very much Serre), Kantor, Feit, Cameron, 
Sathaye, Eakin, 
Stennerson, Gorenstein, O'Nan, Mulay, and Neumann.

The pilgrimage started when in September and October of 
1988, Serre
sent me one after another four long letters briefly saying 
that

{\it {\rm``}In your {\rm1957} paper you suggested that the 
Galois group of
$\oF_n$ should be calculated. I can now prove that for 
$t=1$ it is
$\operatorname{PSL}(2,p)$. Can you calculate it for other 
values of $t${\rm?} 
Also, the conjectures in your paper include AGC. Can you 
now prove
AGC\/}{\rm?''}

Fortunately, in his last letter, Serre added a sentence 
saying that
``my e-mail is...\,.'' 

%\centerline{}
%\centerline{{\bf }}
\heading 15. Multiply transitive groups\endheading

Having already commented on the significance of 
transitivity for Galois theory,
before proceeding further with calculations of Galois 
groups, let us give a 
brief review of multiply transitive groups.

So let $G$ be a permutation group, say of degree $n$, 
i.e., let $G$ act on
$\uw = \{1,2,\dots,n \}$.  In analogy with the concept of 
transitivity
introduced in \S3, we say that $G$ is $\l$-{\it 
antitransitive}
(or $\l$-fold antitransitive) for a positive integer $\l 
\leq n$, if for
all pairwise distinct elements $i_1, i_2, \dots , i_{\l}$ 
in $\uw$ we 
have that the identity is the only member of $G$ which 
keeps them fixed.
\footnote{ That is, if $\s \in G$ is such that 
$\s(i_e)=i_e$ for $1 \leq
e \leq \l$ then we must have $\s = 1$.}
Moreover, for positive integers $l \leq \l \leq n$, we say 
that $G$ is
$(l,\l)$-{\it transitive} if $G$ is $l$-transitive and 
$\l$-antitransitive;
we may express this by simply saying that $G$ is $(l,\l)$.
\footnote{Note that if $G$ is $l$-transitive and 
$m$-antitransitive for
positive integers $l \leq n$ and $m \leq n$, then 
automatically
$l \leq m$.  Also note that if $G$ is $l$-transitive for 
some positive 
integer $l \leq n$ then $G$ is $l'$-transitive for every 
positive
integer $l' \leq l$.  Likewise, if $G$ is 
$\l$-antitransitive for a positive 
integer $\l$ then $G$ is $\l'$-antitransitive for every 
positive integer $\l'
\geq \l$ with $\l' \leq n$.}
Finally, $G$ is {\it sharply $l$-transitive} means $G$ is 
$(l,l)$.  Now
if $G$ is $l$-transitive, with $l > 1$, then the one-point 
stabilizer
of $G$ is obviously $(l-1)$-transitive
as a permutation group of degree $n-1$, acting on the 
``remaining'' $n-1$
elements; conversely, if $G$ is transitive and its 
one-point stabilizer is
$(l-1)$-transitive then $G$ is $l$-transitive.  
Similarly, if $G$ is $\l$-antitransitive, with $\l > 1$, 
then the one-point stabilizer of $G$ is 
$(\l-1)$-antitransitive  
as a permutation group of degree $n-1$, acting on the 
``remaining'' $n-1$
elements; conversely, if $G$ is transitive and its 
one-point stabilizer is
$(l-1)$-antitransitive then $G$ is $l$-antitransitive.  
Thus, to classify all $(l,\l)$ groups, we can make 
induction and each time 
increase $l$, $\l$,  and $n$ by one.  So we start with 
$(1,1)$.

By definition, $G$ is {\it regular} means $G$ is $(1,1)$.  
Now the 
classification of $(1,1)$ groups is either obvious or 
impossible.  
Obvious because it is so easy to define what is $(1,1)$.  
Indeed,
the usual proof of the usual theorem which says that every 
finite group 
is a permutation group, amounts to representing the given 
group
as a regular permutation group.
\footnote{ By making $G$ act on itself by right or left 
(but not
both) multiplication.}
Hence impossible because it would amount to classifying 
all finite groups.

$G$ is {\it Frobenius\/} means $G$ is $(1,2)$ but not 
$(1,1)$.  A 
prototype of a Frobenius group is the group of all affine 
linear 
transformations $ax + b$ with $a,b$ in a finite field, or 
more generally
in a finite near-field.  Zassenhaus \cite{Z2}, in his 1936 
Thesis written
under Artin, proved that the converse is true for $(2,2)$ 
groups.

\proclaim{Zassenhaus' Theorem} 
%{\eightrm $G$ is $(2,2) \Leftrightarrow
%G = \operatorname{AGLNF}(1,\uc)$ for some finite 
%near-field $\uc$.}
$G$ is $(2,2) \Leftrightarrow
G = \operatorname{AGLNF}(1,\uc)$ for some finite 
near-field $\uc$.
\endproclaim 

Here, by $\operatorname{AGLNF}(1,\uc)$ we are denoting the 
group of all affine linear
transformations of degree $1$ over the near-field $\uc$, 
where 
``near-field'' is a generalization of ``field'' obtained 
by weakening the
distributive law.  Namely, a {\it near-field} is an additive
abelian group $\uc$ in which the nonzero elements form a 
multiplicative 
group such that for all $a,b,c$ in $\uc$ we have $a(b+
c)=ab+ac$.  Thus
we are not assuming the other distributive law $(b+c)a = 
ba + bc$
which is not a consequence of the first distributive law 
because the 
multiplication is not required to be commutative.  It can 
easily be 
seen that, as in the case of a field, the number of 
elements in any
finite near-field is a power of a prime number.  By an 
affine linear
transformation of degree $1$ over the near-field $\uc$ we 
mean a map
$\uc \to \uc$ given by $x \mapsto ax+b$ with $a,b$ in 
$\uc$ and $a \neq 0$.
It can easily be seen that distinct $(a,b)$ give distinct 
maps
$\uc \to \uc$.  A proof of Zassenhaus' Theorem is given in 
9.10 on page
424 of volume III of Huppert-Blackburn \cite{HB}.
\footnote{ Although Dickson found all finite near-fields, 
it was left to
Zassenhaus to prove that there were no more.  E. H. Moore 
\cite{Mo} of
the newly opened University of Chicago classified finite 
fields, around
1895, and then, around 1905, his students Wedderburn 
\cite{W} and Dickson
\cite{D2} studied skew fields and near-fields 
respectively.}  
In the proof of this Theorem, as well as in the proofs of 
various other
theorems on multitransitive groups, an important role is 
played by the
following Theorem of Frobenius (1901) \cite{Fr} for a 
proof of which we
refer to 8.2 on page 496 of volume I of Huppert--Blackburn 
\cite{HB}.

\proclaim{Frobenius' Theorem}  A Frobenius group $G$ 
always has a $(1,1)$
normal subgroup.  More precisely, the subset of $G$ 
consisting of the identity
together with those elements $\s$ which fix no letter 
\RM(i.e., $\s(i) \neq i$
for $i=1,2, \dots,n)$ forms a regular normal subgroup of 
$G$.
\endproclaim

In the notation $\operatorname{AGLNF}(1,\uc)$, 
the letters $\operatorname{NF}$ are meant to remind us
of a near-field.  In case $\uc = $ a field $\uf$, we may 
write $\operatorname{AGL}(1,\uf)$
instead of $\operatorname{AGLNF}(1,\uc)$.

To put the notation $\operatorname{AGL}(1,\uf)$ in proper 
perspective, first we remark that
for a field $\uf$ and a positive integer $m$ the groups 
$\operatorname{GL}(m,\uf), 
\operatorname{PGL}(m,\uf), \operatorname{SL}(m,\uf)$,
and $\operatorname{PSL}(m,\uf)$ are defined by replacing 
$\operatorname{GF}(q)$ by $\uf$ in 
\S5.
\footnote{Writing $q$ for $\operatorname{GF}(q)$ is 
justified because for any prime power
$q$, up to isomorphism, there is exactly one field with 
$q$ elements.}
Note that $Z(\operatorname{GL}(m,\uf))=$ the set of all 
scalar matrices, and $Z(\operatorname{SL}(m,\uf))
= \operatorname{SL}(m,\uf) \cap 
Z(\operatorname{GL}(m,\uf))$, where the {\it center} of 
any group
$\ug$ is denoted by $Z(\ug)$, i.e., $Z(\ug)$ is the normal 
subgroup of
$\ug$ given by putting $Z(\ug) =\{a \in \ug : ab=ba \text{ 
for all } 
b \in \ug \}$.  Now a nonsingular $m$ by $m$ matrix $\a 
\in \operatorname{GL}(m,\uf)$
corresponds to the bijection $\uf^m \to \uf^m$ which sends 
any $1$ by
$m$ matrix $\x \in \uf^m$ to the matrix product $\x \a \in 
\uf^m$.  In
this manner, the group $\operatorname{GL}(m,\uf)$, and 
hence also the subgroup $\operatorname{SL}(m,\uf)$,
may be regarded as a permutation group on $\uf^m$.  Let 
$\cp(\uf^m)$
be the $(m-1)$-dimensional {\it projective space} over 
$\uf$, where
we think of $\cp(\uf^m)$ as the set of all one-dimensional 
subspaces
of $\uf^m$.  Now the bijection $\uf^m \to \uf^m$ 
corresponding to any
$\a \in \operatorname{GL}(m,\uf)$ clearly induces a 
bijection $\cp(\uf^m) \to
\cp(\uf^m)$;  moreover, if $\a^* \in 
\operatorname{GL}(m,\uf)$ differs from $\a$ by a 
scalar matrix, then $\a$ and $\a^*$ induce the same 
bijection
$\cp(\uf^m) \to \cp(\uf^m)$.  Thus the group 
$\operatorname{PGL}(m,\uf)$, and hence also
the subgroup $\operatorname{PSL}(m,\uf)$, becomes a 
permutation group on $\cp(\uf^m)$.  
Members of $\operatorname{PGL}(m,\uf)$ are called {\it 
projective transformations} of
$\cp(\uf^m)$.  Note that for any $\g \in 
\operatorname{PGL}(m,\uf)$ and $\z \in \cp(\uf^m)$
we have $\g(\z) \in \cp(\uf^m)$.

Now the {\it affine general linear group\/} 
$\operatorname{AGL}(m,\uf)$ of degree $m$
over $\uf$ may be introduced as the {\it semidirect 
product\/}
$\uf^m \rtimes \operatorname{GL}(m,\uf)$ of $\uf^m$ by 
$\operatorname{GL}(m,\uf)$ 
with the obvious action
of $\operatorname{GL}(m,\uf)$ on $\uf^m$, where we recall 
that a group $\ug$ is said
to be the (internal) semidirect product of a normal 
subgroup $\uh$ by
a subgroup $\ud$ provided $\ug = \uh \ud$ and $\uh \cap 
\ud = 1$, and
we note that in this case $\ud$ acts on $\uh$ by 
conjugation.  Concretely,
$\operatorname{AGL}(m,\uf)$ may be regarded as the set of 
all $m$ by $m+1$ matrices whose
entries are in the $\uf$ and whose $m$ by $m$ piece is 
nonsingular;  in
other words, we think of $\operatorname{AGL}(m,\uf)$ as 
the set of all pairs $(\a,\b)$
with $\a \in \operatorname{GL}(m,\uf)$ and $\b \in \uf^m$, 
where multiplication is
defined by $(\a,\b)(\a',\b') = (\a\a',\b\a'+\b')$.
\footnote{This is an example of an ``external'' semidirect 
product.  For
further elucidation see Suzuki \cite{Su2}, 
Huppert-Blackburn \cite{HB},
and Wielandt \cite{Wi}.  A lot of the group theory 
background, required
in this paper, I learned in the last two years from these 
nice books.  I
highly recommend them.  It may be noted that a Frobenius 
group $G$ is the
semidirect product of the Frobenius kernel of $G$ by a 
Frobenius
complement of $G$, where by the {\it Frobenius kernel} of 
$G$ we mean
the $(1,1)$ normal subgroup of $G$, and by a {\it 
Frobenius complement}
of $G$ we mean a $1$-point stabilizer of $G$.}
To $(\a,\b) \in \operatorname{AGL}(m,\uf)$ there 
corresponds the bijection
$\uf^m \to \uf^m$ which sends every $\x \in \uf^m$ to 
$\x\a+\b \in \uf^m$.
\footnote{ In the above case of $m=1$, this reduces to $x 
\mapsto ax+b$
by taking $x = \x \in \uf$, $0 \neq a = \a \in \uf$, and 
$b = \b \in \uf$.}
Thus $\operatorname{AGL}(m,\uf)$ also becomes a 
permutation group on $\uf^m$.  Members
of $\operatorname{AGL}(m,\uf)$ are called {\it affine 
transformations} of $\uf^m$.
\footnote{Transformations, or substitutions, of the type 
$x' = ax + by + c$ 
and $y' = a'x+b'y+c'$ are familiar to us from high school. 
 By
adding ``points at infinity'' to the ordinary plane we get 
the projective
plane.  To distinguish between the ordinary plane and the 
projective plane,
the ordinary plane is called the affine plane and the 
above transformations 
are called affine transformations of the plane $\uf^2$.  
To know that they
can be redefined in terms of semidirect products should 
help to make
this notion friendly.}

To make another example of a semidirect product, let $\Aut 
\uf$ be the
{\it group of all automorphisms} of the field $\uf$.  For 
every matrix
$\g$ with entries in $\uf$ and for every $g \in \Aut \uf$, 
let $\g^g$ be 
the matrix obtained by applying $g$ to each entry of $\g$. 
 This gives an
action of $\Aut \uf$ on $\operatorname{GL}(m,\uf)$.  The 
semidirect product $\operatorname{GL}(m,\uf) \rtimes
\Aut \uf$ is denoted by $\ug\operatorname{L}
(m,\uf)$ and members of $\ug\operatorname{L}(m,\uf)$ are
called {\it semilinear transformations} of $\uf^m$.  A 
member of
$\ug \operatorname{L}(m,\uf)$ may be thought of as a pair 
$(g,\a)$ with $g \in \Aut \uf$ and
$\a \in \operatorname{GL}(m,\uf)$, and the corresponding 
bijection $\uf^m \to \uf^m$
sends every $\x \in \uf^m$ to $\x^g  \a \in \uf^m$.  The 
multiplication
in $\ug\operatorname{L}(m,\uf)$ is given by 
$(g,\a)(g',\a')=(gg',\a^{g'}\a')$.  
Thus
$\ug\operatorname{L}(m,\uf)$, which may be called the {\it 
semilinear group} of degree $m$
over $\uf$, acts on $\uf^m$.  With this action we can form 
the semidirect
product $\operatorname{A}\ug\operatorname{L}
(m,\uf) = \uf^m \rtimes \ug \operatorname{L}(m, \uf)$ and 
call it
the {\it affine semilinear group} of degree $m$ over 
$\uf$.  A member of
$\operatorname{A}\ug \operatorname{L}(m,\uf)$ is called an 
{\it affine semilinear transformation} of
$\uf^m$ and it may be represented as a triple $(g,\a,\b)$ 
with $g \in \Aut \uf$,
$\a \in \operatorname{GL}(m,\uf)$, $\b \in \uf^m$;  the 
corresponding bijection
$\uf^m \to \uf^m$ sends every $\x \in \uf^m$ to $\x^g\a + 
\b \in \uf^m$, and
the multiplication in $\operatorname{A}\ug 
\operatorname{L}(m,\uf)$ is given by 
$(g,\a,\b)(g',\a',\b')=(gg',\a^{g'}\a',\b^{g'}\a'+\b')$. 
Thus $\operatorname{A}
\ug \operatorname{L}(m,\uf)$
is a permutation group on $\uf^m$, and in a natural manner 
we have
$$
\operatorname{SL}(m,\uf) \triangleleft 
\operatorname{GL}(m,\uf) \triangleleft \ug 
\operatorname{L}(m,\uf) < 
\operatorname{A} \ug \operatorname{L}(m,\uf)
$$
and
$$
\operatorname{GL}(m,\uf) < \operatorname{AGL}(m,\uf) < 
\operatorname{A}
 \ug \operatorname{L}(m,\uf),
$$
where $<$ and $\triangleleft$ denote {\it subgroup} and 
{\it normal subgroup}
respectively.

To construct one more example of a semidirect product, we 
note that the action
of $\Aut \uf$ on $\operatorname{GL}(m,\uf)$ obviously 
induces an action of $\Aut \uf$ on
the factor group $\operatorname{PGL}(m,\uf)$.  With this 
induced action we form the
semidirect product $\operatorname{P}\ug\operatorname{L}
(m,\uf) = \operatorname{PGL}(m,\uf) \rtimes \Aut \uf$ and
call it the {\it projective semilinear group} of degree 
$m$ over $\uf$.
A member of $\operatorname{P}\ug\operatorname{L}
(m,\uf)$ is called a {\it projective semilinear 
transformation} of $\cp (\uf^m)$ and it may be represented 
by a pair
$(\g,g)$ with $\g \in \operatorname{PGL}(m,\uf)$ and $g 
\in \Aut \uf$; the corresponding
bijection $\cp(\uf^m) \to \cp(\uf^m)$ sends every $\z \in 
\cp(\uf^m)$ to
$\g(\z^g) \in \cp(\uf^m)$ where $\z^g\in\cp(\uf^m)$ is 
given by 
$\z^g=\{\x^g:\x \in \z \}$.  Thus
$\operatorname{P}\ug\operatorname{L}
(m,\uf)$ becomes a permutation group on $\cp (\uf^m)$, and 
in a
natural manner we have
$$
\operatorname{PSL}(m,\uf)\triangleleft 
\operatorname{PGL}(m,\uf) \triangleleft 
\operatorname{P}\ug\operatorname{L}
(m,\uf).
$$

Although we have spoken of
$\operatorname{SL}(m\!,\uf)\!,\operatorname{GL}(m\!,\uf),%
\ug\operatorname{L}(m
\!,\uf)\!,\operatorname{AGL}(m\!,\uf)\!,\operatorname{A}\ug
\operatorname{L}(m,\uf)$ as permutation groups on $\uf^m$, 
and $\operatorname
{PSL}(m,\uf), \operatorname{PGL}(m,\uf),\operatorname{P}\ug
\operatorname{L}(m,\uf)$ as permutation groups on $\cp 
(\uf^m)$, this is 
relevant mainly when $\uf = \operatorname{GF}(q)$ for some 
prime power $q$ and 
in that case we may write 
$\operatorname{SL}(m,q)$,$\operatorname{GL}
(m,q)$, $\ug\operatorname{L}(m,q)$, 
$\operatorname{AGL}(m,q)$,
$\ug\operatorname{L}(m,q)$, $\operatorname{PSL}(m,q)$, 
$\operatorname{PGL}(m,q)$, 
and $\operatorname{P}\ug\operatorname{L}(m,q)$ for 
$\operatorname{SL}(m,\uf)$, $\operatorname{GL}
(m,\uf)$, $\ug\operatorname{L}(m,\uf)$, 
$\operatorname{AGL}(m,\uf)$,
$\operatorname{A}\ug \operatorname{L}(m,\uf)$,
%\newline 
$\operatorname{PSL}(m,\uf)$, $\operatorname{PGL}(m,\uf)$, 
and $\operatorname{P}
 \ug \operatorname{L}(m,\uf)$ 
respectively.

Henceforth by a {\it permutation group} we shall again 
mean a permutation group
on a finite set.

%\centerline{}
%\centerline{{\bf }}
\heading 16. Zassenhaus groups\endheading

Having talked about
$(1,1)$, $(1,2)$, and $(2,2)$ groups, let us now discuss 
$(2,3)$ groups
which are not $(2,2)$.  Basically they fall into the 
following three classes.

(i) A {\it Feit group} is defined to be a $(2,3)$ group 
which is 
not $(2,2)$ but has a $(1,1)$ normal subgroup.  As a 
prototype 
we have the group $\operatorname{A}\ug\operatorname{L}
(1,2^p)$ where $p$ is a prime number. 
This consists of all transformations $x \mapsto ax^g+b$ 
with $0 \neq
a \in \operatorname{GF}(2^p)$ and $b \in 
\operatorname{GF}(2^p)$ and $g \in \Aut 
\operatorname{GF}(2^p)$.  Now
$|\operatorname{GF}(2^p)|=2^p$ and $|\Aut 
\operatorname{GF}(2^p)|=p$ where $|\;|$ denotes 
cardinality.
Thus $\operatorname{A}\ug \operatorname{L}(1,2^p)$ is a 
permutation group of degree $2^p$
and order $2^p(2^p-1)p$.  

(ii) A {\it sharp Zassenhaus group} is defined to be a 
$(3,3)$ group; such a 
group is clearly a $(2,3)$ group; moreover, it is a 
$(2,2)$ group only when 
its degree is 3 and in that case it is simply $S_3$.
The Fundamental Theorem of Projective Geometry says that, 
on the projective
line over a field $\uf$, any three points can be sent to 
any other
three points by one and only one projective 
transformation.  In case,
$\uf = \operatorname{GF}(q)$, where $q$ is a prime power, 
\footnote{That is $q$ is a positive integral power of a 
prime number.
Clearly then $q$ is an odd prime power or an even prime 
power according as the
corresponding prime number is even or odd. When $q$ is an 
odd prime power,
it can be an even power of an odd prime or an odd power of 
an odd prime.} 
this amounts to saying that the group 
$\operatorname{PGL}(2,q)$ is a $(3,3)$ group, where we 
regard  $\operatorname{PGL}(2,q)$ as a permutation group 
of degree $q+1$; 
\footnote{$\operatorname{PGL}(2,q)$ acts on the projective 
line over $\operatorname{GF}(q)$ which
has $q+1$ points on it, out of which $q$ are at ``finite 
distance''
and one is the point at ``infinity''.} 
it is easily seen that $\operatorname{PGL}(2,q)$ has no 
$(1,1)$ 
normal subgroup; clearly the order of 
$\operatorname{PGL}(2,q)$ is $(q+1)q(q-1)$.
\footnote{Obviously the order of any $(l,l)$ group of 
degree $n$ is 
$n(n-1)\dotsb(n-l+1)$.}
It can be shown that, if $q$ is an even power of an odd 
prime, $\operatorname{P}\ug 
\operatorname{L}(2,q)$ 
has exactly one $3$-transitive subgroup, which we denote by 
$\operatorname{PML}(2,q)$, such 
that $\operatorname{PML}(2,q)\ne \operatorname{PGL}(2,q)$ 
and 
$\operatorname{PSL}(2,q)$ is a subgroup of 
$\operatorname{PML}
(2,q)$ of index $2$. Now $\operatorname{PML}(2,q)$ is 
also a $(3,3)$ group, where we again regard it as a 
permutation group of 
degree $q+1$; 
it is easily seen that $\operatorname{PML}(2,q)$ has no 
$(1,1)$ normal subgroup
provided $q\geq4$; clearly the order of 
$\operatorname{PML}(2,q)$ is $(q+1)q(q-1)$; we call 
$\operatorname{PML}(2,q)$ the 
{\it projective mock linear group} of degree $2$ over 
$\operatorname{GF}(q)$. 
For further discussion about $\operatorname{PML}(2,q)$ see 
page 163 of volume III of
Huppert-Blackburn \cite{HB} where it is denoted by $M(q)$.
 
(iii) A {\it strict Zassenhaus group} is defined to be a 
$(2,3)$ group 
which is neither $(2,2)$  nor $(3,3)$ and does not have 
any $(1,1)$ normal 
subgroup.  
\footnote{In Gorenstein's book \cite{G1} every $(2,3)$ 
group which
is not $(2,2)$ is called a Zassenhaus group.  In Huppert 
and Blackburn's book
\cite{HB} every $(2,3)$ group which is not $(2,2)$ and 
does not have
any $(1,1)$ normal subgroup is called a Zassenhaus group.  
We are calling
the groups mentioned in (i) Feit groups because they were 
completely
characterized by Feit \cite{F} in 1960.}
It can easily be seen that, for any odd prime power $q$, 
the group
$\operatorname{PSL}(2,q)$ is a strict Zassenhaus group of 
degree $q+1$.
To find the order of this group, we might as well start by 
calculating the 
order of $\operatorname{GL}(m,q)$ for any positive integer 
$m$ and any prime power $q$
which need not be odd. Now the number of ways of choosing 
the first column
of an element of $\operatorname{GL}(m,q)$ equals 
$|\operatorname{GF}(q)^m|-1=q^m-1$. The number of multiples 
of the first column is $q$ and hence, having chosen the 
first column, the 
number of ways of choosing the second column equals
$|\operatorname{GF}(q)^m|-|\operatorname{GF}(q)|=q^m-q$. 
More generally, the first $i$ columns generate an 
$i$-dimensional vector space over $\operatorname{GF}(q)$ 
and hence, having chosen
the first $i$ columns, the number of ways of choosing the 
$(i+1)\text{th}$
columns equals $|\operatorname{GF}(q)^m|-|%
\operatorname{GF}(q)^i|=q^m-q^i$. Therefore
$$
|\operatorname{GL}(m,q)|=(q^m-1)(q^m-q)\dotsb(q^m-q^{m-1}).
$$
Consequently
$$
|\operatorname{AGL}(m,q)|=|\operatorname{GF}(q)^m||
\operatorname{GL}(m,q)|=q^m(q^m-1)(q^m-q)%
\dotsb(q^m-q^{m-1}).
$$
Now $|\operatorname{GF}(q)^*|=q-1$, where 
$\operatorname{GF}(q)^*$ is the multiplicative group of
nonzero elements of $\operatorname{GF}(q)$, and we have 
obvious exact sequences of finite
groups $1\to \operatorname{SL}(m,q)\to 
\operatorname{GL}(m,q)\to \operatorname{GF}(q)^*\to 1$ and
$1\to \operatorname{GF}(q)^*\to \operatorname{GL}(m,q)\to 
\operatorname{PGL}(m,q)\to 1$, and hence
$$
\split
|\operatorname{PGL}(m,q)|&=|\operatorname{SL}(m,q)|=|%
\operatorname{GL}
 (m,q)|/(q-1)\\
&=(q^m-1)(q^m-q)\dotsb(q^m-q^{m-2})q^{m-1}.
\endsplit
$$
Let
$$
\split
Z_{\n}=&\text{the cyclic group of order }\n\text{ where }\n 
\text{ is either a positive }\\
&\text{integer or }\infty.
\endsplit
$$
Now clearly 
$$
\Aut \operatorname{GF}(q)=Z_{\m}\text{ where 
}q=p^{\m}\text{ with }p= \CHAR k
$$
and hence
$$
\split
|\ug \operatorname{L}
(m,q)|/|\operatorname{GL}(m,q)|&=|\operatorname{A}\ug 
\operatorname{L}
(m,q)|/|\operatorname{AGL}(m,q)|\\
&=
|\operatorname{P}\ug 
\operatorname{L}(m,q)|/\operatorname{PGL}(m,q)|=\m.
\endsplit
$$
Let $\operatorname{SL}(m,q)^*$ be the group of all $m$ by 
$m$ scalar matrices whose
entries are in $\operatorname{GF}(q)$ and whose 
determinant is 1. 
Then $\operatorname{SL}(m,q)^*$
is isomorphic to the group of all $m\text{th}$ roots of 1 
in $\operatorname{GF}(q)$
and hence 
$|\operatorname{SL}(m,q)^*|=\operatorname{GCD}(m,q-1)$. 
Also we have an obvious exact
sequence of finite groups $1\to \operatorname{SL}(m,q)^*\to 
\operatorname{SL}(m,q)\to \operatorname{PSL}(m,q)\to 1$
and hence
$$
\split
|\operatorname{PSL}(m,q)|&=|\operatorname{SL}(m,q)|/%
\operatorname{GCD}(m,q-1)\\
&=|\operatorname{PGL}(m,q)|/\operatorname{GCD}(m,q-1).
\endsplit
$$
So in particular
$$
|\operatorname{PSL}(2,q)|=(q+
1)q(q-1)/2\quad\text{or}\quad(q+1)q(q-1)
$$
according as $q\text{ is odd or even}$
i.e., according as $q-1$ is or is not divisible by 2.
If $q$ is even then we cannot divide $q-1$ by 2, and so we 
do the best we 
can; namely, assuming $q$ to be a square of a proper odd 
power 
\footnote{That is with an odd exponent $>1$.}
of 2, in the
expression $(q+1)q(q-1)/2$ we replace $(q-1)/2$ by 
$(q^{1/2}-1)$ to
get the expression $(q+1)q(q^{1/2}-1)$, and now to get rid 
of the
fractional power we write $q$ in place of $q^{1/2}$. For 
every $q$ which
is a proper odd power of 2, the resulting expression $(q^2+
1)q^2(q-1)$ is the 
order of a certain strict Zassenhaus group 
$\operatorname{Sz}(q)$ of degree $q^2+1$; 
this group $\operatorname{Sz}(q)$ is isomorphic 
\footnote{Not as a permutation group.}
to a certain subgroup of $\operatorname{GL}(4,q)$ and, 
since it was discovered by Suzuki \cite{Su1} in 1962, it 
is called the 
{\it Suzuki group} 
over $\operatorname{GF}(q)$.
We may think of $\operatorname{Sz}(q)$ as ersatz 
$\operatorname{PSL}(2,q^2)$; we have just
given a heuristic reason for its existence and a mnemonic 
device for 
remembering its order; to recapitulate
$$
|\operatorname{Sz}(q)|=(q^2+1)q^2(q-1)\quad\text{if $q$ is 
any odd power of 2}.
$$ 
As hinted in the above order formula, the definition of 
the Suzuki group 
$\operatorname{Sz}(q)$ can be extended so as
to include the case of $q=2$; in this case we still get a 
2-transitive 
permutation group of degree $q^2+1$ and of the above 
order, which is however
not a strict Zassenhaus group; indeed, as a permutation 
group, $\operatorname{Sz}(2)$
is isomorphic to the $(2,2)$ group 
$\operatorname{AGL}(1,5)$.

The following theorem of Zassenhaus \cite{Z1}, Feit 
\cite{F}, and 
Suzuki \cite{Su1} says that the above
examples of $(2,3)$ groups which are not $(2,2)$ are 
exhaustive; for a proof 
see 1.1 and 11.16 on pages 161 and  286 of volume III of 
Huppert-Blackburn \cite{HB}.
\proclaim{Zassenhaus-Feit-Suzuki Theorem}  
For a permutation group $G$ we have the following.
\roster
\item \<$G$ is Feit $\Leftrightarrow G=\operatorname{A}\ug 
\operatorname{L}
(1,2^p)$ for some prime $p$.
\item \<$G$ is sharp Zassenhaus $\Leftrightarrow 
G=\operatorname{PGL}(2,q)$ for some prime
power $q$, or $G=\operatorname{PML}(2,q)$ for some even 
power $q$ of an odd prime.
\item \<$G$ is strict Zassenhaus $\Leftrightarrow 
G=\operatorname{PSL}(2,q)$ for some odd 
prime power $q$, or $G=\operatorname{Sz}(q)$ for some 
proper odd power $q$ of $2$.
\endroster
Moreover, this gives an exhaustive and mutually exclusive 
listing of
$(2,3)$ groups which are not $(2,2)$, with the proviso 
that the group
$\operatorname{PSL}(2,2)=S_3$ is included in item $(2)$ 
even though it is a $(2,2)$
group in addition to being a $(3,3)$, and hence a $(2,3)$, 
group.
\endproclaim

%\centerline{}
%\centerline{{\bf }}
\heading 17. More about classification theorems\endheading

By analyzing the Suzuki group $\operatorname{Sz}(q)$,
a certain analogous 2-transitive permutation group 
$R_1(q)$ of degree 
$q^3+1$, for every odd power $q$ of 3,
was discovered by Ree \cite{R} in 1964; although $R_1(q)$ 
is 2-transitive, 
it is not 3-antitransitive; the group $R_1(q)$ is defined 
in terms of some 
7 by 7 matrices over
$\operatorname{GF}(q)$ and is called the Ree group over 
$\operatorname{GF}(q)$; the order of $R_1(q)$ is given
by
$$
|R_1(q)|=(q^3+1)q^3(q-1)\quad\text{if $q$ is any odd power 
of 3}.
$$
For every odd power
$q$ of 2, in terms of certain matrices over 
$\operatorname{GF}(q)$, Ree defined a group
$R_2(q)$ which is also called the Ree group over 
$\operatorname{GF}(q)$; the order is now
given by
$$
|R_2(q)|=q^{12}(q^6+1)(q^4-1)(q^3+1)(q-1)
$$
if $q$ is any odd power of 2.
All the Suzuki groups and all the Ree groups turn out to 
be a simple groups 
except: $\operatorname{Sz}(2)$ is solvable; $R_1(3)$ is 
isomorphic to 
$\operatorname{P}\ug \operatorname{L}(2,8)$ and hence
the simple group $\operatorname{PSL}(2,8)$ 
may be regarded as a normal subgroup
of $R_1(3)$ of index 3 and then $\operatorname{PSL}(2,8)$ 
is the unique minimal normal
subgroup
of $R_1(3)$; the commutator subgroup $R_2(2)'$ of $R_2(2)$ 
is a (normal)
subgroup of $R_2(2)$ of index 2, and the said commutator 
subgroup $R_2(2)'$
is actually a simple group called the Tits group.
Thus we have the following three infinite families of 
finite simple
groups deduced from matrices over $\operatorname{GF}(q)$. 
%\proclaim
\subheading{Suzuki and Ree groups} $\operatorname{Sz}(q)$ 
for every proper odd power $q$ of $2$.
$R_1(q)$ for every proper odd power $q$ of $3$.
$R_2(q)$ for every proper odd power $q$ of $2$, together 
with the commutator 
subgroup $R_2(2)'$.
%\endproclaim

Just before Suzuki and Ree found these groups, Steinberg 
\cite{St}, in
1959, discovered that some known groups together with some 
further ones
suggested by them could be organized into 
four other infinite families of finite simple
groups which are defined in terms of matrices over 
$\operatorname{GF}(q^2)$ and which 
are now labeled as follows. 
%\proclaim
\subheading{Steinberg groups}${}^2A_n(q), {}^2D_n(q), 
{}^3D_4(q), \text{ and }
{}^2E_6(q)$. Here $n$ is any positive integer and $q$ is 
any prime power
except: in case of $\,^2A_n(q)$ exclude 
$(n,q)=(1,2),(1,3),(2,2)$, and in case
of $\,^2D_n(q)$ exclude $n=1$.
%\endproclaim

These $3+4=7$ families are the ``twisted incarnations'' of 
the original
nine infinite families of ``relatives'' of $A_n(q)$.  
These nine infinite 
families of finite simple groups defined in terms of 
matrices over 
$\operatorname{GF}(q)$ 
were systematized by Chevalley \cite{Ch} in 1955, and they 
are labeled as 
follows. 
\footnote{When $q=$ a prime $p$, many of these were 
already studied by
Jordan \cite{J1} in the last century. For general $q$\<, 
some of them were
discussed by Dickson \cite{D1} at the turn of the century. 
Likewise the
first two families of Steinberg groups were already known 
to Jordan and
Dickson, while the last two were independently found by 
Tits.}
%\proclaim
\subheading{Chevalley groups} $A_n(q), B_n(q), C_n(q), 
D_n(q), E_6(q), E_7(q), 
E_8(q), F_4(q)$, and $G_2(q)$. Again here $n$ is any 
positive integer and $q$ is any prime power except: in 
case of $A_n(q)$ exclude 
$(n,q)=(1,2),(1,3)$; in case of $B_n(q)$ and $C_n(q)$  
exclude 
$(n,q)=(1,2),(1,3),(2,2)$; in case of $D_n(q)$ exclude 
$n<4$; and in case of
$G_2(q)$ exclude $q=2$.
%\endproclaim

Thus we have $7+9=16$ infinite families of finite simple 
groups defined
in terms of matrices over finite fields.  
It may be noted that, just as the group 
$A_n(q)=\operatorname{PSL}(n+1,q)$ is obtained
by projectivizing the special linear group 
$\operatorname{SL}(n+1,q)$, 
the group $C_n(q)=\operatorname{PSp}(2n,q)$ is obtained 
by projectivizing the ``symplectic group'' 
$\operatorname{Sp}(2n,q)$, and the groups
$B_n(q)=\operatorname{P}
\uw(2n+1,q)$ and $D_n(q)=\operatorname{P}
\uw^+(2n,q)$ are obtained by projectivizing 
the ``commutator groups'' $\uw(2n+1,q)$ and $\uw^+(2n,q)$ 
of the 
``orthogonal groups'' $\operatorname{O}
(2n+1,q)$ and $\operatorname{O}^+(2n,q)$ respectively.
The groups $A_n(q),\ B_n(q),\ C_n(q),$$D_n(q)$ are 
collectively called 
``classical groups'', and the remaining Chevalley groups 
$E_6(q),E_7
(q),E_8(q),F_4(q),G_2(q)$ are collectively called 
``exceptional 
groups.''

Likewise,
amongst the Steinberg groups, the group 
$^2A_n(q)=\operatorname{PSU}
(n+1,q)$ is obtained 
by projectivizing the ``special unitary group'' 
$\operatorname{SU}(n+1,q)$, and the group 
$^2D_n(q)=\operatorname{P}
\uw^{-}(2n,q)$ is obtained by projectivizing the 
``commutator 
group'' $\uw^-(2n,q)$ of the ``orthogonal group'' 
$\operatorname{O}
^-(2n,q)$.
The Suzuki and Ree groups have the alternative labels 
$\operatorname{Sz}(q)={}^2B_2(q)$,
$R_1(q)={}^2G_2(q)$, and $R_2(q)={}^2F_4(q)$. The groups 
${}^2A_n(q)$,
${}^2B_2(q)$, ${}^2D_n(q)$, ${}^3D_4(q)$ 
may collectively be called the ``twisted classical 
groups'', 
and the groups ${}^2E_6(q)$, ${}^2F_4((q)$, ${}^2G_2(q)$ 
may collectively be called the ``twisted exceptional 
groups.''

Now the projective mock linear group
$\operatorname{PML}(2,q)$ may also be called the {\it 
premathieu linear group} 
of degree $2$ over $\operatorname{GF}(q)$.  The reason for 
this nomenclature is that, 
around 1865, Mathieu \cite{Mat} found a transitive
extension of $\operatorname{PML}(2,9)$ which is denoted by 
$M_{11}$.  By a {\it transitive
extension} of a permutation group of degree $n$ we mean a 
transitive
permutation group of degree $n+1$ having the given group 
as a one-point
stabilizer.  Note that $M_{11}$ is a $(4,4)$ group of 
degree $11$ 
and hence its order is $11 \cdot 10 \cdot 9 \cdot 8$.  
Mathieu also 
found a transitive extension of $M_{11}$ which we denote 
by $M_{12}$.
Clearly $M_{12}$ is a $(5,5)$ group of degree $12$ and 
hence its order
is $12 \cdot 11 \cdot 10 \cdot 9 \cdot 8$. The permutation 
groups
$M_{11}$ and $M_{12}$ are called {\it Mathieu groups} of 
degree $11$ and $12$
respectively. It can be shown that $M_{12}$ has a noninner 
automorphism
$\a$ of order $2$ such that $\a(M_{11})$ is 3-transitive 
of degree 12.
The permutation group $\a(M_{11})$ may be regarded as an 
{\it incarnation} of 
the Mathieu group $M_{11}$ and we may denote it by 
$\hM_{11}$.
To recapitulate
$$
|M_{12}|=12 \cdot 11 \cdot 10 \cdot 9 \cdot 8 \qquad\text{ 
and }\qquad
|\hM_{11}|=|M_{11}| =11 \cdot 10 \cdot 9 \cdot 8.
$$

To introduce the remaining three groups discovered by 
Mathieu, let us first
note that for any positive integer $m$ and any prime power 
$q$ we have 
$|\operatorname{GF}(q)^m|
=q^m$ and hence for the corresponding $(m-1)$-dimensional
projective space we have
$$
|\cp(\operatorname{GF}(q)^m)|=(q^m-1)/(q-1)=q^{m-1}+
q^{m-2}+\cdots+q+1.
$$
So in particular $|\cp(\operatorname{GF}(4)^3)|=21$ and 
hence $\operatorname{PSL}(3,4)$, 
which acts on$\cp
(\operatorname{GF}(4)^3)$, is a permutation group of 
degree $21$; it can easily be seen
that it is $2$-transitive but not $3$-transitive; 
by the above order formula we also get
$$
|\operatorname{PSL}(3,4)|=(4^3-1)(4^3-4)4^2/%
\operatorname{GCD}(3,3)=21\cdot 20\cdot 48.
$$
Mathieu obtained a transitive extension $M_{22}$ of 
$\operatorname{PSL}(3,4)$, a transitive
extension $M_{23}$ of $M_{22}$, and a transitive extension 
$M_{24}$ of
$M_{23}$. Clearly $M_{22}$ is a 3-transitive but not 
4-transitive permutation group of degree 22, 
$M_{23}$ is a 4-transitive but not 5-transitive 
permutation of 
degree 23, and $M_{24}$ is
a 5-transitive but not 6-transitive permutation group of 
degree 24. 
These groups are called {\it Mathieu groups} of degree 22, 
23, 24
respectively, and obviously their orders are
\footnote{For a transitive permutation group $G$ of degree 
$n$ we clearly
have $|G|=n|G_1|$ where $G_1$ is a one-point stabilizer of 
$G$.}
$|M_{22}| = 22 \cdot 21 \cdot 20 \cdot 48$;
$|M_{23}| = 23 \cdot
22 \cdot 21 \cdot 20 \cdot 48$,
 and $|M_{24}| = 24 \cdot 23 \cdot 
22 \cdot 21 \cdot 20\cdot 48$. 
Note the striking similarity between the numbers $(2,9)$ 
occuring in the
``parent group'' $\operatorname{PML}(2,9)$ of $M_{11}$ and 
$M_{12}$, and the numbers $(3,4)$
occuring in the ``parent group'' $\operatorname{PSL}(3,4)$ 
of $M_{22}$,
$M_{23}$, and $M_{24}$.
All of the five Mathieu groups $M_{11}$, $M_{12}$, 
$M_{22}$, $M_{23}$, 
and $M_{24}$
turn out to be simple groups.  One hundred years after 
their discovery, 
during 1965 to 1975, twenty-one other finite simple {\it 
sporadic groups}, 
i.e., those which do not naturally fit in any infinite 
family, were
discovered by various people; the  
largest of the $21+5=26$ sporadic groups is called the 
{\it monster}
and its order is
$$
2^{46}\cdot 3^{20}\cdot 5^9\cdot 7^6\cdot 11^2\cdot 
13^3\cdot 17\cdot 
19\cdot 23\cdot 29\cdot 31\cdot 41\cdot 47\cdot 59\cdot 71;
$$
for details see Gorenstein \cite{G1}. 
We are now
ready to state, of course without proof, the 
Classification Theorem
of Finite Simple Groups.
\footnote{See Gorenstein \cite{G2} or Aschbacher 
\cite{As}. At least one
part of this extremely long proof, namely Mason's paper on 
quasi-thin
groups \cite{Mas}, is still to see the light of day!} 
\proclaim{CT} The following is a complete list of
finite simple groups.
\roster
\item The cyclic group $Z_p$ for every prime $p$.
\item The alternating group $A_n$ for every integer $n 
\geq 5$.
\item The sixteen infinite families of ``matrix'' groups 
mentioned above.
\footnote{I am using the more friendly term ``matrix 
groups'' instead of
the awe inspiring ``Lie type groups.''}
\item The twenty-six sporadics including the five Mathieus.
\endroster
\endproclaim

Just as $\operatorname{PSL}(m,q)$ is the typical example 
of a finite simple
group, so$\operatorname{PGL}(2,q)$ is the quintessential 
example of 
a $3$-transitive permutation group. As obvious variations 
of this, additional
3-transitive permutation groups are obtained by taking 
into account 
all the groups between $\operatorname{PGL}(2,q)$ and 
$\operatorname{P}
\ug \operatorname{L}(2,q)$ for 
every prime power $q$\!, and all the groups between 
$\operatorname{PML}(2\!,q)$ and
$\operatorname{P}\ug \operatorname{L}(2\!,q)$ for every 
even power $q$ of an odd prime; both these types 
have degree $q+1$. All these arise out
of the fact that, by a projective transformation, any 3 
points of a 
projective line can be sent to any other 3. Going to 
higher dimensions,
by a projective transformation, any 4 points of a 
projective plane can be
sent to any other 4, and any 5 points of a projective 
3-space can be sent to
any other 5, and so on. This would give rise to 
4-transitives,
5-transitives, and so on. But there is a flaw. 3 collinear 
points in a
projective plane, or in a projective 3-space, cannot be 
sent to 3
noncollinear points. 4 coplanar points in a projective 
3-space cannot be
sent to 4 noncoplanar points. And so on. Thus, because of 
questions of
linear independence, for every integer $m>1$ and every 
prime power $q$, 
the group $\operatorname{PGL}(m,q)$, instead of being 
$(m+1)$-transitive, is only 2-transitive, unless every 
line contains
only 2 distinct points, in which case it would be 
3-transitive.
Well, for $q=2$, a line should contain only 2 points. But 
even that is so
only in the affine case because then we don't have the 
point at infinity.
Thus, it is not $\operatorname{PGL}(m,2)$ which is 
3-transitive, but $\operatorname{AGL}(m,2)$.
We can see that $\operatorname{AGL}(m,2)$ is, however, not 
4-transitive unless $m=2$
in which case we actually have 
$\operatorname{AGL}(2,2)=S_4$.

Now
$$
\operatorname{AGL}(m,2)=\operatorname{GF}(2)^m \rtimes 
\operatorname{GL}(m,2)
=2^m \cdot L_m(2),
$$
where a  dot stands for the semidirect
product symbol $\rtimes$, and $2^m$ stands for $(Z_2)^m	$, 
i.e., for the direct
product $Z_2 \times Z_2 \times \dots \times Z_2$ of $m$ 
copies of
$Z_2$.  By the order formula we have 
$$
|L_4(2)|=(2^4-1)(2^4-2)(2^4-4)(2^4-8)=8!/2=|A_8|.
$$
Hence, by the philosophical principle that {\it two finite 
simple groups
of equal order are usually isomorphic}, we expect that 
$L_4(2) \approx
A_8$ where $\approx$ stands for isomorphism, and this can, 
in fact, be
easily proved.  Note that $L_4(2)=\operatorname{GL}(4,2)$ 
is a 1-point stabilizer of
$\operatorname{AGL}(4,2)$ and hence in this incarnation 
$A_8$ is only 2-transitive;
\footnote{By taking the stabilizer at the origin, 
$\operatorname{GL}(m,q)$ becomes the
1-point stabilizer of $\operatorname{AGL}(m,q)$ for every 
integer $m>1$ and every
prime power $q$. Likewise, by taking the stabilizer at the 
point at 
infinity, $\operatorname{AGL}(1,q)$ may be regarded as the 
1-point stabilizer of
$\operatorname{PGL}(2,q)$ for every prime power $q$.} 
let us denote this {\it incarnation} by $\hA_8$. 
\footnote{That is, as a permutation group, 
$\hA_8=\operatorname{GL}(4,2)=\operatorname{PGL}(4,2)%
\subset S_{15}$.} 
In a natural manner, $A_7$ may
be regarded as a subgroup of $A_8$ and then it turns out 
that the image of
$A_7$ under the said isomorphism is also only 
2-transitive; let us denote 
the said image by $\hA_7$.  
\footnote{In other words, $\hA_7$ is the image of $A_7$ 
under some injective
group homomorphism 
$A_7\to\hA_8=\operatorname{GL}(4,2)\subset S_{15},$ 
and $\hA_7$ is a 2-transitive
but not 3-transitive permutation group of degree 15.}
The corresponding subgroup of 
$2^4 \cdot \hA_8 = \operatorname{AGL}(4,2)$ may be denoted 
by $2^4 \cdot \hA_7$; 
this is a 3-transitive but not 4-transitive permutation 
group of degree
$2^4$ and order $2^3\cdot 7!$.
\footnote{As  permutation groups, $\hA_7$ and $2^4\cdot 
\hA_7$ 
are independent of the
injective group homomorphism $A_7\to\hA_8$ we choose for 
defining $\hA_7$.}

%As another way of defining $\hA_7$, first note that 
%$|\cp(\operatorname{GF}(2))^3|=7$ and
%hence we have several injective group hompomorphisms 
%$\operatorname{PGL}(3,2)\to S_7$; 
%it can be shown that the image of 
%$\operatorname{PGL}(3,2)$ is then a subgroup of $A_7$ of 
%index 15, and so by taking the (right) multiplication 
%action on the (right)
%cosets of the said subgroup, we get a permutation 
%representation of $A_7$ of 
%degree 15 which can be seen to be 2-transitive but not 
%3-transitive; 
%these representations of $A_7$ given by the various 
%injective group
%homomorphisms of $\operatorname{PGL}(3,2)$ give two 
%distinct 2-transitive but not
%3-transitive permutation 

As a consequence of CT, it can be shown that there are no 
more 3-transitive
permutation groups other than those we have already 
listed. In other
words we have the following detailed version of CTT, i.e., 
the 
Classification Theorem of Triply Transitive Permutation 
Groups;
this theorem was compiled from conversations with Cameron, 
Neumann,
and O'Nan.
\proclaim{CTT or Refined Fundamental Theorem of Projective 
Geometry}
The 
following is a complete list of \RM3-transitive 
permutation groups.
\roster
\item
For every prime power $q$, each group between 
$\operatorname{PGL}(2,q)$ and$\operatorname{P}
\ug\operatorname{L}(2,q)$ is a {\rm3}-transitive permutation
group of degree $q+1$, and we have 
$|\operatorname{PGL}(2,q)|=$ %\newline
$|\operatorname{P}\ug \operatorname{L}(2,q)|/\m=(q+1)q(q-1)$
where $q=p^{\m}$ with $p=\CHAR \operatorname{GF}(q)$.
Among these, $\operatorname{PGL}(2,3)$ and 
$\operatorname{P}\ug 
\operatorname{L}(2,4)$ are the only groups which are
{\rm4}-transitive, and for them we actually have 
$\operatorname{PGL}(2,3)=S_4$ and
$\operatorname{P}\ug \operatorname{L}(2,4)=S_5$.
\item
For every even power $q$ of an odd prime, each group 
between $\operatorname{PML}
(2,q)$ and
$\operatorname{P}\ug \operatorname{L}(2,q)$ is a  
{\rm3}-transitive but not
{\rm4}-transitive permutation
group of degree $q+1$, and we have 
$|\operatorname{PML}(2,q)|=|\operatorname{P}
\ug \operatorname{L}(2,q)|/\m= (q+1)q(q-1)$
where $q=p^{\m}$ with $p=\CHAR \operatorname{GF}(q)$.
\item
For every integer $m>1$, the group 
$\operatorname{AGL}(m,2)$ is a 
{\rm3}-transitive 
permutation group of degree $2^m$ and order 
$2^m(2^m-1)(2^m-2)\dotsb(2^m-2^{m-1})$.
This is {\rm4}-transitive only for $m=2$, and in that case 
we actually have
$\operatorname{AGL}(2,2)=S_4$.
\item
The group $2^4\cdot\hA_7$ is a {\rm3}-transitive but not 
{\rm4}-transitive permutation 
group of degree $2^4$ and order $2^4\cdot 7!$, and as a 
permutation group it is a subgroup of 
$2^4\cdot\hA_8=\operatorname{AGL}(4,2)$.
\item
The reincarnated Mathieu group $\hM_{11}$ is a 
{\rm3}-transitive but not
{\rm4}-transitive permutation group of degree $12$ and 
order 
$11\cdot 10\cdot 9\cdot 8=12\cdot 11\cdot 10\cdot 6$.
\item
The Mathieu group $M_{22}$ and its automorphism group 
$\Aut M_{22}$
are {\rm3}-transitive but not {\rm4}-transitive 
permutation groups of degree 
{\rm22} 
with $|M_{22}|=|\Aut M_{22}|/2 = 22\cdot 21\cdot 20\cdot 
48$.
\item
The Mathieu groups $M_{11}$ and $M_{23}$ are 
{\rm4}-transitive but not
{\rm5}-transitive permutation groups of degree {\rm11} and 
{\rm23} and order
$11\cdot 10\cdot 9\cdot 8$ and $23\cdot 22\cdot 21\cdot 
20\cdot 48$ 
respectively.
\item
The Mathieu groups $M_{12}$ and $M_{24}$ are 
{\rm5}-transitive but not
{\rm6}-transitive permutation groups of degree {\rm12} and 
{\rm24} and order
$12\cdot 11\cdot 10\cdot 9\cdot 8$ and 
$24\cdot 23\cdot 22\cdot 21\cdot 20\cdot 48$ respectively.
\item
For every integer $n\ge 5$, the alternating group $A_n$ is 
an
$(n-2)$-transitive but not $(n-1)$-transitive permutation 
group
of degree $n$ and order $n!/2$.
\item
For every integer $n\ge 3$, the symmetric group $S_n$ is an
$n$-transitive but not $(n+1)$-transitive permutation group
of degree $n$ and order $n!$.
\endroster
\endproclaim 

The above formulation of CTT obviously subsumes CQT, CFT, 
and CST.
In turn the CTT is of course subsumed under the CDT which 
is given by
Cameron \cite{C} and Kantor \cite{K2} and which, in 
addition to heavily using
the CT, is based on the previous work of 
Curtis-Kantor-Seitz \cite{CKS},
O'Nan \cite{O}, and others. The following weaker version 
of CDT, dealing mainly
with the degrees of 2-transitive permutation groups which 
are not
3-transitive, was communicated to me by Cameron.

\proclaim{Weak CDT} Concerning the degrees of 
\RM2-transitive permutation
groups we have the following.
\roster
\item
For every integer $m>1$ and every prime power $q$, each 
group between
$\operatorname{PSL}(m,q)$ and 
$\operatorname{P}\ug\operatorname{L}(m,q)$ is a
\RM2-transitive 
permutation group of degree 
$|\cp(\operatorname{GF}(q)^m)|=(q^m-1)/(q-1)$. Out of 
these,
only the groups listed in items {\rm(1)} and {\rm(2)} of 
CTT are \RM3-transitive.
In case of $m>2$, for each group between 
$\operatorname{PSL}(m,q)$ and 
$\operatorname{P}\ug\operatorname{L}(m,q)$,
by considering the action on ``hyperplanes'' in 
$\cp(\operatorname{GF}(q)^m)$, we get a
second representation as a \RM2-transitive but not 
\RM3-transitive permutation
group of degree $(q^m-1)/(q-1)$.
\item
For every integer $m>2$, the group 
$\operatorname{Sp}(2m,2)$ has  \RM2-transitive but
not \RM3-transitive permutation representations of degrees 
$2^{2m-1}+2^{m-1}$
and $2^{2m-1}-2^{m-1}$.
\item 
For every prime power $q=p^{\m}>2$ with prime $p$, each 
group between 
$\operatorname{PSU}(3,q)$ and its automorphism group $\Aut
\operatorname{PSU}(3,q)$ has a \RM2-transitive but
not \RM3-transitive permutation representation of degree 
$q^3+1$, and moreover
$\operatorname{PSU}(3,q)$ is a normal subgroup of index 
$[\operatorname{GCD}(3,q+1)]2\m$ in 
$\Aut \operatorname{PSU}(3,q)$.
\footnote{A {\it centerless} group $G$, i.e., a group $G$ 
whose center
is the identity, may be identified with its inner 
automorphism group and
hence may be regarded as a normal subgroup of its 
automorphism
group $\Aut G$.} 
\item
For every proper odd power $q=2^{\m}$ of $2$, each group 
between the Suzuki 
group $\operatorname{Sz}(q)$ and its automorphism group 
$\Aut \operatorname{Sz}(q)$ has 
a \RM2-transitive but not \RM3-transitive permutation 
group of degree $q^2+1$,
and moreover $\operatorname{Sz}(q)$ is a normal subgroup 
of index $\m$ in $\Aut \operatorname{Sz}(q)$.
\item
For every odd power $q=3^{\m}$ of $3$, each group between 
the Ree group 
$R_1(q)$ and its automorphism group $\Aut R_1(q)$ has a 
representation as a 
\RM2-transitive but not \RM3-transitive permutation group 
of degree $q^3+1$,
and moreover $R_1(q)$ is a normal subgroup of index $\m$ 
in $\Aut R_1(q)$.
\item
The group $\operatorname{PSL}(2,11)$ has two distinct 
\RM2-transitive but not 
\RM3-transitive 
permutation representations of degree {\rm11}.
\footnote{Two permutation representations $G\to S_n$ and 
$G\to S_{n'}$
of a (finite) group $G$ are {\it equivalent} if they 
differ by an isomorphism
$S_n\to S_{n'}$ induced by a bijection between the 
underlying sets of $S_n$
and $S_{n'}$; note that then automatically $n=n'$. Two 
representations are
{\it distinct} if they are not equivalent.}
\item
The alternating group $A_7$ has two distinct 
\RM2-transitive but not 
\RM3-transitive permutation representations of degree 
{\rm15;} 
both are equivalent to isomorphisms $A_7\to\hA_7\subset 
S_{15}$.
%one of which is $\hA_7$ and the other is its ``flip'' by 
%a noninner 
%automorphism; (but there is no canonical way to tell 
%which is $\hA_7$ and
%which is its flip). 
\item
The ``Higman-Sims'' sporadic group {\rm HS} has two 
different  
\RM2-transitive but 
not \RM3-transitive permutation representations of degree 
{\rm176}.
\item
The ``third Conway'' sporadic group $\operatorname{Co}_3$ 
has a 
\RM2-transitive but not
\RM3-transitive permutation representation of degree 
{\rm276}.
\endroster
The above items {\rm(1)} to {\rm(9)} contain a complete 
list of \RM2-transitive
but not \RM3-transitive permutation groups having a 
nonabelian 
minimal normal subgroup. The degree of a \RM2-transitive 
but not \RM3-transitive
permutation group $G$ having an abelian minimal normal 
subgroup is necessarily
some power $p^m$ of some prime $p$\RM; the said minimal 
normal subgroup is
isomorphic to $(Z_p)^m$ and a \RM1-point stabilizer of the 
group $G$ itself is
isomorphic to a subgroup of $\operatorname{GL}(m,p)$\RM;
moreover\RM: in the case $m=1$ we must have 
$G=\operatorname{AGL}(1,p)$, whereas, in 
the case
$m=$ a prime number and $p=2$  we must have 
$G=\operatorname{AGL}(1,2^m)$ or 
$\operatorname{A}\Gamma\operatorname{L}(1,2^m)$,
and finally, in the case $m=2$ and $p=3$ we must have 
$G=\operatorname{AGL}(1,9)$ or
$\operatorname{A}\Gamma\operatorname{L}(1,9)$ or 
$\operatorname{AML}(1,9)$
or $\operatorname{AGL}(2,3)$ or $\operatorname{ASL}(2,3)$, 
where we have put $\operatorname{AML}(1,9)$ = the \RM1-point
stabilizer of $\operatorname{PML}(2,9)$,
\footnote{In other words, $\operatorname{AML}(1,9)$ is the 
2-point stabilizer of the
Mathieu group $M_{11}$.}
and 
$\operatorname{ASL}(2,3)=(\operatorname{GF}(3))^2\rtimes 
\operatorname{SL}(2,3)$.
\endproclaim 
 
In connection with the above statement, we note the 
following Theorem of
Burnside which is really the starting point of the 
classification of
2-transitive permutation groups. Although most modern proofs
of this make use of Frobe- nius' Theorem (1901) \cite{Fr},
\footnote{For Frobenius' Theorem there is no ``character 
free'' proof.
As examples of modern proofs of Burnside's Theorem which 
seem to use
Frobenius' Theorem, see 12.4 on page 32 of Wielandt 
\cite{Wi} 
and 7.12 on page 233 of volume III of Huppert-Blackburn 
\cite{HB}.}
it is clear that Burnside's original proof did not, since 
it is already
given as Theorem IX on page 192 of the first edition of 
\cite{Bu}
published in 1897. In the second edition of \cite{Bu} 
published in
1911, it occurs as Theorem XIII on page 202, and there 
Burnside gives two 
proofs of it, one using Frobenius and the other without.
\footnote{In Burnside's classical style of writing, 
Theorem $x$ means
Theorem $x$ together with the discussion around it. In 
other words,
although everything is proved, only some of the 
conclusions are
called theorems. This ``classical'' style is quite 
different from the so 
called ``Landau Style'' of Satz-Beweis-Bemerkung. 
In the classical style, you first discuss things and then 
suddenly
say that you have proved such and such; in other words, 
the proof 
precedes the statement of a theorem.}
Thus Burnside gives an ``elementary'' proof of the 
following Theorem
without using ``character theory.''
\footnote{And certainly without using CT!}

\proclaim{Burnside's Theorem}
A \RM2-transitive permutation group has a unique minimal 
normal subgroup.
The said subgroup is either an elementary abelian group
\footnote{A group is {\it elementary abelian} if it is 
isomorphic to 
$(Z_p)^m$ for some positive integer $m$ and some prime $p$.}
or a nonabelian simple group.
\endproclaim

As a consequence of CTT and Weak CDT, and in view of a 
simple numerical
lemma, we have the following.
%\proclaim
\thm{Special CDT}
Given any prime $p$ and any positive integer $\m$, 
concerning
\RM2-transitive but not \RM3-transitive permutation groups 
of degree $q + 1$, 
where $q = p^{\m}$, we have the following. 
If $\m = 1$ and $p$ is a Mersenne prime, 
then $\operatorname{PSL}(2,p), \operatorname{AGL}(1,p+1)$ 
and 
$\operatorname{A}\Gamma\operatorname{L}(1,p+1)$ are the 
only such groups. If $\m = 1$ but $p$ is not a Mersenne 
prime, then $\operatorname{PSL}(2,p)$ is the only such 
group, except that for $p=2$ this 
group is ``accidentally'' \RM3-transitive because it 
coincides with $S_3$. If 
$\m > 1$ then, in addition to the relevant groups listed 
in items {\rm(1), 
(3), (4), (5)} of Weak CDT, the only other such groups are 
the groups
$\operatorname{AGL}(1,9)$, 
$\operatorname{A}\Gamma\operatorname{L}(1,9)$,
$\operatorname{AML}(1,9)$, 
$\operatorname{AGL}(2,3)$, and $\operatorname{ASL}(2,3)$, 
which occur when $(\m,p)=(3,2)$,
and the group $\operatorname{AGL}(1,q+1)$ which occurs 
when $q+1$ is a Fermat prime.
\footnote{A {\it Fermat prime} is a prime of the form 
$2^{\m}+1$ for some 
positive integer $\m$. It follows then that $\m$ must be a 
power of $2$,
because otherwise $\m=\m'\m''$ where $\m'$ is even and 
$\m''>1$ is odd and
this would give the nontrivial factorization
$2^{\m}+1=(2^{\m'}+1)(2^{\m'(\m''-1)}-2^{\m'(\m''-2)}+
\dots+1)$.}  
\ethm
%\endproclaim

Here is the said 
\proclaim{Simple Numerical Lemma}
Given any primes $p$ and $\p$ and any positive integers 
$\m$ and $u$,
such that $p^{\m}+1=\p^u$, we have the following. If 
$p>2$, then $\m=1$ and
$\p=2$ and $u=$ a prime number, and so $p$ must be a 
Mersenne prime.
If $p=2$ and $u>1$, then $\m=3$ and $\p=3$ and $u=2$. If 
$p=2$ and $u=1$, then
$\m=$ a power of \RM2, and so $\p$ must be a Fermat prime.
\endproclaim

To see this, first suppose that $p>2$; now $p$ is odd and 
hence $p^{\m}+1$
is even and hence $\p=2$; since $p$ is odd, we have 
$p\equiv 1(4)$ or
$p\equiv 3(4)$ and hence $p^2\equiv 1(4)$, and therefore 
if $\m$ is even then
$2^u=p^{\m}+1\equiv 2(4)$ and this would imply $u=1$ which 
would be 
absurd; on the other hand, if $\m$ is odd then
$2^u=p^{\m}+1=(p+1)(p^{\m-1}-p^{\m-2}+\dots+1)$ where the 
second parenthesis
contains an odd number of odd terms and hence its value is 
odd, but since
that value divides $2^u$, it must be 1, and this gives 
$p^{\m}+1=p+1$ 
which implies $\m=1$.
Next suppose that $p=2$ and $u>1$; now $\p$ must be odd, 
and if $u$ is
also odd then $2^{\m}=\p^u-1=(\p-1)(\p^{u-1}+\p^{u-2}+
\dots+1)$ where the
second parenthesis consists of $u$ positive odd terms and 
hence its value is
an odd number $\ge u$ and, since it divides $2^{\m}$, it 
must be $1$ which is
absurd; on the other hand, if $u$ is even then 
$2^{\m}=\p^u-1=(\p^{u/2}-1)(\p^{u/2}+1)$ and hence
$\p^{u/2}-1$ and $\p^{u/2}+1$ are both powers of 2 whose 
difference
is 2 and therefore they must be equal to 2 and 4, and this 
gives
$\m=3$ and $\p=3$ and $u=2$. Finally, if $p=2$ and $u=1$, 
then by the last
footnote we see that $\m$ must be a power of 2, and $\p$ 
must be a Fermat prime.

%As a consequence of CTT and Weak CDT we have the following.
%
%\proclaim{Special CDT}
%For any odd prime $p$, the group 
%$\operatorname{PSL}(2,p)$ is a 2-transitive but not
%3-transitive permutation group of degree $p+1$. If $p$ is 
%a Mersenne prime
%then $\operatorname{AGL}(1,p+1)$ and 
%$\operatorname{A}\Gamma\operatorname{L}(1,p+1)$ are the 
%only other 2-transitive
%but not 3-transitive permutation groups of degree $p+1$. 
%If $p$ is not a
%Mersenne prime then $\operatorname{PSL}(2,p)$ is the only 
%such permutation group.
%\endproclaim

Here is another consequence of CDT.
\proclaim{Uniqueness Theorem for Transitive Extensions}
Any two transitive extensions of a transitive permutation 
group are
isomorphic as permutation groups, with only one exception. 
\footnote{If we don't assume the given group to be 
transitive,
then there are numerous exceptions. For example every 
finite group,
in its standard representation as a regular permutation 
group,
is a transitive extension of the identity group. Since for 
increasing $l$, 
there are fewer and fewer $l$-transitive permutation 
groups, it follows
that ``most'' transitive permutation groups have no 
transitive
extensions.}
The exception is that $\operatorname{PSL}(2,7)$ and 
$\operatorname{A}\Gamma\operatorname{L}(1,8)$ have a 
common
\RM1-point stabilizer\RM; note that both these are 
\RM2-transitive but
not \RM3-transitive permutation groups of degree \RM8. 
\endproclaim

As an immediate corollary of the above theorem we have the 
following.

\proclaim{Uniqueness Theorem for Transitive Extensions of 
2-Transitive
Groups} %\newline 
Any two transitive extensions of a \RM2-transitive 
permutation group are
isomorphic as permutation groups. In particular, the Mathieu
groups $M_{11},M_{12},M_{22},M_{23}$, and $M_{24}$ are the 
unique 
transitive extensions of 
$\operatorname{PML}(2,9),M_{11},%
\operatorname{PSL}(3,4),M_{22},
\text{ and } M_{23}$ respectively.
\endproclaim

To end this review of group theory, we note that by
the {\it rank} of a transitive permutation group is meant 
the number of
orbits of its 1-point stabilizer; the lengths of these 
orbits, excluding
the obvious one point orbit, are called {\it subdegrees} 
of the group;
so the number of subdegrees is one less than the rank, and 
the sum of
the subdegrees is one less than the degree. Thus a 
2-transitive group is
simply a transitive group of rank 2. Now CT has also been 
used by
Kantor-Liebler \cite{KL}, Liebeck \cite{L}, and others, to 
give CR3
= classification of transitive groups of rank 3, which 
although much longer
than CDT, should be quite useful for Galois theory. Here 
is an amusing
sample from CR3 which does not use CT and which can be 
found in Kantor
\cite{K1}.
%\proclaim{
\thm{Sample from CR3}
For any integer $n>1$ and any prime power $q$, the groups 
$\operatorname{PSp}(2n,q)$
and $\operatorname{O}(2n+1,q)$ are the only transitive 
permutation groups of
rank \RM3
whose subdegrees are $q(q^{2n-2}-1)/(q-1)$ and $q^{2n-1}$.
A rank \RM3 transitive permutation group $G$ with 
subdegrees $q(q+1)^2$ and $q^4$
for a prime power $q>1$, is a subgroup of $\Aut 
\operatorname{PSL}(4,q)$\RM; moreover, if
$q>2$ then $G$ contains $\operatorname{PSL}(4,q)$.
%\endproclaim
\ethm

%\centerline{}
%\centerline{{\bf }}
\heading 18. A type of derivative\endheading

To continue with the calculation of Galois groups, let me 
explain how to 
throw away a root $\a=\a_1$ of a polynomial
$$
f=f(Y)=Y^n+a_1Y^{n-1}+\dots+a_n=\prod^n_{i=1}(Y-\a_i)
$$ 
by using a type of derivative. Now the coefficients 
$a_1,\dots,a_n$ belong
to a field $K$, and we want to find the polynomial
$$
f_1=f_1(Y)=\frac{f(Y)}{(Y-\a)}=Y^{n-1}+b_1Y^{n-2}+\dots+
b_{n-1}\in K(\a)[Y].
$$
To this end, first recall the three basic transformations 
of equations
described in any old book. For instance, we may quote the 
following three 
relevant articles (= sections) from Burnside-Panton's 1904 
book on 
the theory of equations \cite{BP}.
\footnote{For the last forty years I had happily assumed 
that this
Burnside of the theory of equations \cite{BP} was the same 
as the 
Burnside of the theory of groups of finite order 
\cite{Bu}. To my dismay,
at the Oxford Conference in April 1990, Peter Neumann told 
me that,
although both were named William and both obtained a D.Sc. 
from
Dublin around 1890, the equations Burnside was William 
Snow whereas the
group theory Burnside was simply William. Strangely, I 
first learnt group
theory from William Snow's book on the theory of 
equations. }
%\proclaim
\subheading{Art 31. To multiply the roots by a given 
quantity}
For any $u \neq 0$,  the polynomial $g=g(Y)$ whose roots 
are $u$ times 
the roots of $f$ is given by 
$$
g(Y)=u^nf\left(\frac Yu\right)= Y^n+c_1Y^{n-1}+\dots+
c_n=\prod^n_{i=1}(Y-u\a_i)
$$
with $c_i=u^ia_i$.
%\endproclaim
%\proclaim
\subheading{Art 32. To reciprocate the roots}
In case $\a_i\neq 0$ for $1 \leq i \leq n$, i.e., in case 
$a_n \neq 0$,
the polynomial $g=g(Y)$ whose roots are the reciprocals of 
the roots of $f$ 
is given by
$$
g(Y)=\frac{Y^n}{a_n}f\left(\frac 1Y\right)=\frac{1}{a_n}(1+
a_1Y+\dots+a_nY^n)
=\prod_{i=1}^n \left(Y-\left(\frac 1 \a_i\right)\right).
$$
%\endproclaim
%\proclaim
\subheading{Art 33. To decrease the roots by a given 
quantity}
For any $u$, the polynomial $g=g(Y)$ whose roots are $-u$ 
plus the roots
of $f$ is given by
$$
g(Y)=f(Y+u)=Y^n+c_1Y^{n-1}+\dots +c_n= \prod_{i=1}^n (Y - 
(\a_i-u))
$$
with $c_1=a_1+nu,\dots,c_n=f(u)$. 
%\endproclaim

Now in the first and the third cases provided $u\in K$, 
and in the second case
without any proviso, we have $g(Y)\in K[Y]$ and, assuming 
the roots
$\a_1,\dots,\a_n$ to be pairwise distinct, we have 
$\operatorname{Gal}(g,K)= \operatorname{Gal}(f,K)$
as permutation groups, and so for Galois theory purposes 
we may conveniently 
modify $f$ by one or more of these three transformations.
 
For example, sometimes it may be easier to compute the 
polynomial $g_1(Y)$
obtained by decreasing the roots of $f_1$ by $\a$. In view 
of what we have
just said, we get 
$\operatorname{Gal}(g_1,K(\a))= \operatorname{Gal}
(f_1,K(\a))$ and hence, assuming 
$f$ to be irreducible in $K[Y]$, for the one-point 
stabilizer $G_1$ of
$G=\operatorname{Gal}(f,K)$ we get 
$G_1=\operatorname{Gal}(g_1,K(\a))$.

Clearly $g_1$ can also be obtained by first decreasing the 
roots of $f$ by
$\a$ to get the polynomial $g=g(Y)=f(Y+\a)$,  and then 
throwing away the 
root $Y=0$ of $g$; this gives $g_1(Y)=g(Y)/Y$, and now 
remembering
that $f(\a)=0$ we get
$$
g_1(Y)=\frac{f(Y+\a)-f(\a)}{Y}.
$$
According to the calculus definition, by taking the 
``limit'' of the RHS as 
$Y$ tends to $0$, we get $f\,'(\a)$. This motivates the 
following definition
according to which $g_1$ {\it turns out to be the twisted 
$Y$-derivative of
$f$ at $\a$.}
\dfn{Definition}
For any polynomial $\h=\h(Y)$ in an indeterminate $Y$ with
coefficients in a field $L$ and for any element $\b$ in 
$L$, we call
$(\h(Y+\b)-\h(\b))/Y$ the {\it twisted} $Y$-{\it 
derivative} of $\h$ 
{\it at} $\b$. 
\enddfn

For a moment let us denote the twisted $Y$-derivative of
$\h$ at $\b$ by $\h'$. Then clearly $\h'=\h'(Y)$ is a 
polynomial in $Y$
with coefficients in $L$, and if $\h\in L$ then $\h'=0$, 
whereas: if
$\h\notin L$ then $\h'\ne 0$ and the $Y$-degree of $\h'$ 
is 1 less than the
$Y$-degree of $\h$, and the two polynomials $\h$ and $\h'$ 
have the same 
leading coefficient, and hence in particular, if $\h$ is 
monic then so in
$\h'$.

Next we note that this is $L$-{\it linear} because for any 
$\d=\d(Y)\in L[Y]$ and $\l,\m\in L$ we have
$$
(\l\d+\m\h)'=\frac{\l\d(Y+\b)+\m\h(Y+
\b)-\l\d(\b)-\m\h(\b)}{Y}
=\l\d'+\m\h'.
$$
However, the usual product rule is to be replaced by a 
{\it twisted
product rule} because by the standard trick of adding and 
subtracting the
same quantity we get
$$
\split
(\d\h)'&=\frac{\d(Y+\b)\h(Y+\b)-\d(Y+\b)\h(\b)}{Y}\\
&\quad+\frac{\d(Y+\b)\h(\b)-\d(\b)\h(\b)}{Y}\\
&=\d^*\h'+\d'\h^{\sharp}
\endsplit
$$
where $\d^*$ is the $Y$-{\it translation} of $\d$ {\it by} 
$\b$ and 
$\h^\sharp$ is the {\it evaluation} of $\h$ {\it at} $\b$, 
i.e., 
$\d^*=\d^*(Y)=\d(Y+\b)$ and
$\h^\sharp=\h(\b)$.

Finally, for any positive integer $m$ we have the {\it 
power rule}
$$(Y^m)'=Y^{m-1}+m\b Y^{m-2}+\dots+\binom mi\b^i Y^{m-i-1}+
\dots+m\b^{m-1}$$
and, in case $\CHAR L\ne 0$, for any power $q$ of $\CHAR 
L$ we have the
{\it prime power rule}
$(Y^q)'=Y^{q-1}$
and combining this with the product rule, we get the {\it 
power product rule}
$$
[Y^q\h(Y)]'=(Y+\b)^q\h'(Y)+Y^{q-1}\h(\b).
$$
The reason for explicitly mentioning $Y$ in all this is 
that there may be
other indeterminates present; for instance, if 
$\psi=\psi(X,Y)$ is a polynomial
in indeterminates $X$ and $Y$, and $\b$ is an element in a 
field which
contains $X$ as well as all the coefficients of $\psi$, 
then the twisted
$Y$-derivative of $\psi$ at $\b$ is given by
$(\psi(X,Y)-\psi(X,\b))/Y$.

Reverting to the original situation by taking $\b=\a$ and 
$L=K(\a)$ in the 
above set-up, we conclude with the following.
%\proclaim
\subheading{Summary about the twisted derivative} 
If $f=f(Y)$ is 
a nonconstant monic irreducible polynomial in an 
indeterminate $Y$ 
with coefficients in a field $K$ such that $f$ has no 
multiple root in any
overfield of $K$, and if $\a$ is a root of $f$ in some 
overfield of $K$, then
by letting $f\,'=f\,'(Y)$ to be the twisted $Y$-derivative 
of $f$ at $\a$ we have 
that the Galois group 
$\operatorname{Gal}(f\,',K(\a))$ is the one-point 
stabilizer of the
Galois group 
$\operatorname{Gal}(f,K)$.
%\endproclaim

Now without assuming $f$ to be irreducible and without any 
precondition
about multiple roots, suppose the degree of $f$ is $n>2$ 
and suppose 
for every root $\a$ of $f(Y)$ in a splitting field of $K$ 
we have 
that the twisted $Y$-derivative of $f(Y)$ at $\a$ is 
irreducible in 
$K(\a)[Y]$, then $f(Y)$ must be devoid of multiple roots; 
namely, if
$f(Y)=\prod_{i=1}^n(Y-\a_i)$ and $\a_1=\a_2$ and $f\,'(Y)$ 
is the twisted
$Y$-derivative of $f(Y)$ at $\a=\a_1$, then $f\,'(Y)$ is 
reducible in
$K(\a)[Y]$ because its degree is $n-1>1$ and it has 
$(Y-\a_2)$ as a
factor in $K(\a)[Y]$. Thus we have the following.
%\proclaim
\thm{Twisted Derivative Criterion}
If $f(Y)$ is a nonconstant monic polynomial of degree $>2$ 
in an indeterminate
$Y$ with coefficients in a field $K$ such that for every 
root $\a$ of $f(Y)$
in a splitting field of $K$ we have that the 
twisted $Y$-derivative of $f(Y)$ at $\a$ is irreducible in 
$K(\a)[Y]$, then 
$f(Y)$ has no multiple roots in any overfield of $K$.
%\endproclaim
\ethm
%\centerline{}
%\centerline{{\bf}}
\heading 19. Cycle lemma\endheading

As another tool for calculating Galois groups, let us make 
note of a
``cycle lemma''.

Let $K$ be a field and consider a monic polynomial
$$
f=f(Y)=Y^n+a_1Y^{n-1}+\dots+a_n=\prod_{i=1}^n(Y-\a_i)
$$
of degree $n$ in an indeterminate $Y$ with coefficients 
$a_1,\dots,a_n$
in $K$ having pairwise distinct roots $\a_1,\dots,\a_n$ in 
some overfield
of $K$. Now by conveniently enlarging the said overfield 
and moving
it by a $K$-isomorphism, it can be construed to contain 
any preassigned 
overfield $K^*$ of $K$, and this gives us the following
obvious but basic principle of computational Galois theory.
%\proclaim
\thm{Basic Extension Principle}
For any given overfield $K^*$ of $K$, the Galois group 
$\operatorname{Gal}(f,K^*)$, as a 
permutation group of degree $n$, acting on the roots
 $\{\a_1,\dots, \a_n\}$, 
may be regarded as a subgroup of the Galois group 
$\operatorname{Gal}(f,K)$. 
%\endproclaim
\ethm

Given any overfield $K^*$ of $K$ and any factorization
$$
f(Y)=\prod_{j=1}^me_j(Y)\quad\text{where }
e_j(Y)=Y^{n_j}+a_{j1}Y^{n_j-1}+\dots+a_{jn_j}
$$
with $a_{j1},\dots,a_{jn_j}\text{ in }K^*$
we can relabel the roots $\a_1,\dots,\a_n$ as 
$\a_{11},\dots, \a_{1,n_1},\dots$, $\a_{m1},\dots,\a_{mn_m}$ 
so that
$$
e_j=e_j(Y)=\prod_{i=1}^{n_j}(Y-\a_{ji})\quad\text{for 
}1\le j\le m
$$
and we can identify the direct product 
$S_{n_1}\times\dots\times S_{n_m}$, 
where $S_{n_j}$ is the symmetric group acting on 
$\a_{j1},\dots,\a_{jn_j}$, 
with a subgroup of the symmetric group $S_n$ acting on 
$\a_1,\dots,\a_n$.
As a second obvious but basic principle we then have the 
following.
%\proclaim
\thm{Basic Projection Principle}
For the Galois group 
$\operatorname{Gal}(f,K^*)\subset S_n$ 
we have 
$\operatorname{Gal}(f,K^*)\subset S_{n_1}\times\dots\times 
S_{n_m}$, 
and for $1\le j\le m$, the Galois group 
$\operatorname{Gal}(f,K^*)$ 
maps onto the Galois group 
$\operatorname{Gal}(e_j,K^*)\subset S_{n_j}$ 
under the natural projection $S_{n_1}\times\dots\times 
S_{n_m}\to S_{n_j}$.
%\endproclaim
\ethm

Recall that a $\n$-{\it cycle} is a permutation
$\s$\<, say in $S_n$\<, such that for some $\n$ distinct 
elements
$\a_{i_1},\dots\!,\a_{i_{\n}}$ in $\{\a_1,\dots\!,\a_n\}$ 
we have
$\s(\a_{i_1})=\a_{i_2},\dots\!,\s(\a_{i_{\n 
-1}})$$=\a_{i_{\n}},
\s(\a_{i_{\n}})=\a_{i_1}$ and $\s(\a_j)=\a_j$ for all
$j\notin \{i_1,\dots,i_{\n}\}$.
Now if $e_1(Y)$ is irreducible in $K^*[Y]$ and if either 
$n_1$ is prime or
$\operatorname{Gal}(e_1,K^*)$ is cyclic, then clearly 
$\operatorname{Gal}(e_1,K^*)$ contains an $n_1$-cycle
$\t_1$, and by the Projection Principle $\t_1$ is the 
projection of some
$\t\in\operatorname{Gal}(f,K^*)$, and if also 
$|\operatorname{Gal}(e_j,K^*)|$ and $n_1$ are
coprime for $2\le j \le m$ then upon letting $\m$ to be 
the product of
$|\operatorname{Gal}(e_2,K^*)|,\dots,|%
\operatorname{Gal}(e_m,K^*)|$ we see that
$\t^{\m}\in\operatorname{Gal}(f,K^*)$ is an $n_1$-cycle. 
Therefore in view of the
Extension Principle we get the following.

\proclaim{Cycle Prelemma}
If $|\operatorname{Gal}(e_j,K^*)|$ and $n_1$ are coprime 
for $2\le j \le m$, and
$e_1(Y)$ is irreducible in $K^*[Y]$, and either $n_1$ is 
prime or
$\operatorname{Gal}(e_1,K^*)$ is cyclic, then 
$\operatorname{Gal}(f,K)$ contains an $n_1$-cycle.
\endproclaim  

To convert the Cycle Prelemma into the Cycle Lemma, 
let $v$ be a 
(real discrete) 
{\it valuation} 
\footnote{In this paper, by a valuation we shall mean a 
real discrete
valuation.}
of $K$, i.e., $v$ is a map of $K$ onto the set of all 
integers together with
the symbol $\infty$ such that for all $a,b$ in $K$ we have
$v(a)=\infty\Leftrightarrow a=0$, and $v(ab)=v(a)+v(b)$, 
and 
$v(a+b)\ge\text{min}(v(a),v(b)).$ 
Recall that $\{a\in K:v(a)\ge 0\}$ is called 
the {\it valuation ring} of $v$, and this ring modulo the 
unique maximal ideal 
$\{a\in K:v(a)>0\}$ in it is called the {\it residue 
field} of $v$. 
Also recall that $v$ is said to be {\it trivial} on a 
subfield $k$ of $K$, or
$v$ is said to be a valuation of $K/k$, if $v(a)=0$ for 
all $0\ne a\in k$. 
Let $\hK$ be a finite 
algebraic field extension of $K$ and let 
$\hv_1,\dots,\hv_h$ be the extensions
of $v$ to $\hK$, i.e., $\hv_1,\dots,\hv_h$  are those 
valuations of $\hK$
whose valuation rings intersected with $K$ give the 
valuation ring of $v$;
we may also say that $v$ {\it splits} in $\hK$ into 
$\hv_1,\dots,\hv_h$. 
By $\OR(\hv_j:v)$ we denote the {\it reduced ramification 
exponent} 
\footnote{Also called the reduced ramification index.}
of $\hv_j$
over $v$, i.e., $\OR(\hv_j:v)$ is the unique positive 
integer such that for
all $a\in K$ we have $\hv_j(a)=\OR(\hv_j:v)v(a)$. By 
$\od(\hv_j:v)$ we
denote the {\it residue degree} of $\hv_j$ over $v$, i.e., 
$\od(\hv_j:v)$ 
is the field degree of the residue field of $\hv_j$ over 
the residue field
of $v$. Note that {\it if either $\hK/K$ is separable, or 
$v$ 
is trivial over a subfield $k$ of $K$ such that $K/k$ is 
finitely generated 
of transcendence degree $1$, then}
$$
\sum_{j=1}^h\OR(\hv_j:v)\od(\hv_j:v)=[\hK:K].
\tag$\dagger$
$$

%$$
%\aligned
%&\text{if either $\hK/K$ is separable, or $v$ is trivial 
%over a subfield
%$k$ of $K$ such that $K/k$ is}\\ 
%&\text{finitely generated of transcendence degree $1$, 
%then}
%\sum_{j=1}^h\OR(\hv_j:v)\od(\hv_j:v)=[\hK:K].
%\endaligned
%\tag$\dagger$
%$$

Also note that $\hv_j$ is {\it unramified} over $v$, or 
over $K$, means
that $\OR(\hv_j:v)=1$ and the residue field of $\hv_j$ is 
separable over the
residue field of $v$; $\hv_j$ is {\it ramified} over $v$, 
or over $K$, means 
that $\hv_j$ is not unramified over $v$; $v$ is {\it 
unramified} in $\hK$
means that $\hv_j$ is unramified over $v$
for $1\le j\le h$; and finally, $v$ is {\it ramified}
in $\hK$ means that $v$ is not unramified in $\hK$.
Now it is well known that if $f(Y)$ is irreducible in 
$K[Y]$, $\hK=K(\a_1)$,
$K^*=$ the completion of $K$ with respect to $v$, and 
$e_j(Y)$ is irreducible
in $K^*[Y]$ for $1\le j\le m$, then $h=m$ and, after a 
suitable relabelling,
$\OR(\hv_j:v)\od(\hv_j:v)=n_j$ for $1\le j\le h$; for 
instance see \S2
of \cite{A2}. Moreover, by Newton's Theorem, if the 
residue field of $v$ is an
algebraically closed field of the same characteristic as 
$K$, and if 
$n_j\not\equiv 0(\CHAR K)$ for some $j$, then for that $j$ 
the Galois group
$\operatorname{Gal}(e_j,K^*)$ is cyclic; for a proof 
of Newton's Theorem 
based on Shreedharacharya's method of completing the square,
see my new book 
on algebraic geometry 
for scientists and engineers \cite{A6}. 
Therefore by the Cycle Prelemma we get the
\proclaim{Cycle Lemma}
If $f(Y)$ is irreducible in $K[Y]$ and there exists a 
valuation $v$ of $K$
such that the residue field of $v$ is an algebraically 
closed field of the
same characteristic as $K$ and such that for the 
extensions $\hv_1,\dots,\hv_h$
of $v$ to a root field
\footnote{A {\it root field} of $f(Y)$ over $K$ is a field 
obtained by 
adjoining a root of $f(Y)$ to $K$, for instance the field 
$K(\a_1)$.} 
of $f(Y)$ over $K$ we have that $\OR(\hv_j:v)$ and 
$\OR(\hv_1:v)$ 
are coprime and $\OR(\hv_j:v)\not\equiv 0(\CHAR K)$ for 
$1<j\le h$, and 
either $\OR(\hv_1:v)$ is prime or $\OR(\hv_1:v)\not\equiv 
0(\CHAR K)$,
then the Galois group 
$\operatorname{Gal}(f,K)$ contains an $\OR(\hv_1:v)$-cycle. 
\endproclaim

The Basic Extension Principle can be refined thus.
%\proclaim
\thm{Refined Extension Principle}
Given any field extensions $K\subset K'\subset K^*$\<, 
by the Basic Extension Principle we may regard 
$\operatorname{Gal}(f,K^*)<\operatorname{Gal}(f,K')<%
\operatorname{Gal}(f,K)<S_n$, and assuming $K^*$ 
to be a finite normal extension of $K$ we have that
$\operatorname{Gal}(f,K^*)\triangleleft%
\operatorname{Gal}(f,K)$ and the factor group
$\operatorname{Gal}(f,K)/\operatorname{Gal}(f,K^*)$ is a 
homomorphic image of 
$\operatorname{Gal}(K^*,K)$.
\footnote{That is, $\operatorname{Gal}(f,K)/%
\operatorname{Gal}(f,K^*)\approx%
\operatorname{Gal}(K^*,K)/N$
for some normal subgroup $N$ of $\operatorname{Gal}(K^*,K)$.
Note that for any finite normal extension $K^*$ of a field 
$K$, without
assuming $K^*$ to be separable over $K$, the Galois group 
$\operatorname{Gal}(K^*,K)$ is
defined to be the group of all $K$-automorphisms of $K^*$.}
%\endproclaim
\ethm

To see this, let $L=K(\a_1,\dots,\a_n)$ and 
$L^*=K^*(\a_1,\dots,\a_n)$. 
Now $L$ is a (finite) Galois extension of $K$, and
given any $\s\in\operatorname{Gal}(f,K)$, we view $\s$ as 
a permutation of
$\{1,2,\dots,n\}$ such that for some (actually unique) 
$\t\in\operatorname{Gal}(L,K)$ we have 
$\t(\a_i)=\a_{\s(i)}$ for $1\le i\le n$.
Likewise, $L^*$ is a (finite) Galois extension of $K^*$, and
given any $\s^*\in\operatorname{Gal}(f,K^*)$, we view 
$\s^*$ as a permutation of
$\{1,2,\dots,n\}$ such that for some (actually unique)
$\t^*\in\operatorname{Gal}(L^*,K^*)$ we have 
$\t^*(\a_i)=\a_{\s^*(i)}$ for 
$1\le i\le n$. Obviously, $\s = \s^* \Leftrightarrow \t = 
\t^* | L$
where $\t^* | L$ denotes the restriction of $\t^*$ to $L$. 
 Thus
we get the commutative diagram
$$
\CD
\operatorname{Gal}(f,K^*)@>>>\operatorname{Gal}(f,K)\\
@VVV@VVV\\
\operatorname{Gal}(L^*,K^*)@>\d>>\operatorname{Gal}(L,K)
\endCD
$$
where the left arrow is the isomorphism $\s^* \mapsto 
\t^*$, the right
arrow is the isomorphism $\s \mapsto \t$, the top arrow is 
the inclusion
$\operatorname{Gal}(f,K^*) \subset
\operatorname{Gal}(f,K)$, and the bottom arrow $\d$ is the 
injection
$\t^* \mapsto \t^*|L$.  
Therefore our assertion is equivalent to saying that 
$\im\d$ is a normal 
subgroup of 
$\operatorname{Gal}(L,K)$ and 
$\operatorname{Gal}(L\!,K)/\im\d\approx 
\operatorname{Gal}(K^*\!,
K)/N$ for 
some normal subgroup $N$ of 
$\operatorname{Gal}(K^*\!,K)$\!. 
To prove this new version of the assertion,
let $K_0=L\cap K^*$\<, let $K^*_0$ be the maximal 
separable algebraic field
extension of $K_0$ in $K^*$\<, let 
$L^*_0=K^*_0(\a_1,\dots,\a_n)$, and let us 
depict all this in the following Hasse diagram. 
$$
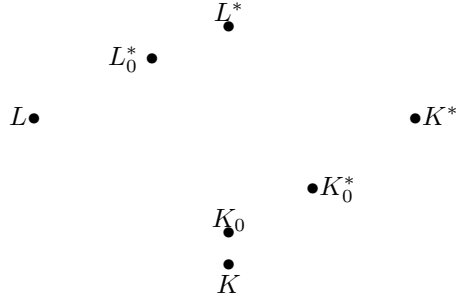
\matrix\format\c&\qquad\c&\qquad\c&\qquad\c&\qquad\c\\
&&L^*&&\\ 
\vspace{-0.1in}
&&\bullet&&\\
&\;\;L^*_0\;\bullet&&&\\
\vspace{0.15in}
L\bullet&&&&\bullet K^*\\
\vspace{0.2in}
&&&\;\bullet K^*_0&\\
&&K_0&&\\
\vspace{-0.1in}
&&\bullet&&\\
&&\bullet&&\\
\vspace{-0.05in}
&&K&&
\endmatrix
$$

Referring to the lower quadrilateral in the above diagram, 
$L^*_0/K_0$ is a 
finite Galois extension, the field $L^*_0$ 
is a compositum of the fields $L$ and 
$K^*_0$ with $L\cap K^*_0=K_0$, and the four sides of the 
said quadrilateral 
represent the finite Galois extensions $L/K_0$, 
$K^*_0/K_0$, $L^*_0/L$, and 
$L^*_0/K^*_0$; hence by the Fundamental Theorem of Galois 
Theory, the group 
$\operatorname{Gal}(L^*_0,K_0)$ is the internal direct 
product of the two normal subgroups
\footnote{That is, the intersection of the two normal 
subgroups is the identity
and they generate the whole group.}
$\operatorname{Gal}(L^*_0,L)$ and 
$\operatorname{Gal}(L^*_0,K^*_0)$, and $\l^*\mapsto\l^*|L$ 
gives a
surjection 
$\operatorname{Gal}(L^*_0,K_0)\to
\operatorname{Gal}(L,K_0)$ whose kernel is 
$\operatorname{Gal}(L^*_0,L)$
and whose restriction to 
$\operatorname{Gal}(L^*_0,K^*_0)$ is an 
isomorphism$\operatorname{Gal}(L^*_0,K^*_0)@>\h>>%
\operatorname{Gal}(L,K_0)$. By applying the Fundamental 
Theorem of Galois Theory to the left triangle in the above 
diagram, i.e., by 
noting that $L/K$, $L/K_0$, and $K_0/K$ are Galois 
extensions, we see that 
$\operatorname{Gal}(L,K_0)$ is a normal subgroup of 
$\operatorname{Gal}(L,K)$ and 
$\operatorname{Gal}(L,K)/\operatorname{Gal}(L,K_0)\approx%
\operatorname{Gal}(K_0,K)$. Upon letting
$\operatorname{Gal}(L,K_0)
\buildrel{\h^*}\over{\longrightarrow}
\operatorname{Gal}(L,K)$ 
be the natural inclusion
$\operatorname{Gal}(L,K_0)\subset\operatorname{Gal}(L,K)$, 
we see that the 
composition 
$\operatorname{Gal}(L^*_0,K^*_0)
\buildrel{\h}\over{\longrightarrow}
\operatorname{Gal}(L,K_0)
\buildrel{\h^\ast}\over{\longrightarrow}
\operatorname{Gal}(L,K)$ 
coincides with the injection
$\operatorname{Gal}(L^*_0,K^*_0)
\buildrel{\d_0}\over{\longrightarrow}
\operatorname{Gal}(L,K)$
given by $\l^*\mapsto\l^*|L$, and hence
$$
\im\d_0=\operatorname{Gal}(L,K_0)\triangleleft%
\operatorname{Gal}(L,K)
$$
and
$$
\operatorname{Gal}(L,K)/\im\d_0\approx%
\operatorname{Gal}(K_0,K).
$$
By applying the Fundamental Theorem of Galois Theory to 
the right triangle in the above diagram, i.e., by noting 
that $K^*_0/K$,
$K^*_0/K_0$, and $K_0/K$ are finite Galois extensions, we 
see that
$$
\operatorname{Gal}(K^*_0,K_0)\triangleleft%
\operatorname{Gal}(K^*_0,K)
\quad\text{and}\quad
\operatorname{Gal}(K^*_0,K)/\operatorname{Gal}(K^*_0,K_0)
\approx\operatorname{Gal}(K_0,K).
$$
Therefore, upon letting $N=\operatorname{Gal}(K^*_0,K_0)$, 
we conclude that
$$
\im\d_0\triangleleft\operatorname{Gal}(L,K)\quad\text{and}%
\quad
\operatorname{Gal}(L,K)/\im\d_0\approx%
\operatorname{Gal}(K^*_0,K)/N
$$
for some
$$
N\triangleleft\operatorname{Gal}(K^*_0,K).
$$
Finally, referring to the modified diagram obtained by 
deleting the three 
lines emanating from $K_0$ in the above diagram, 
$L/K$ and $K^*_0/K$ are finite Galois extensions; $K^*/K$ 
is a finite normal
extension; $K^*/K^*_0$ is pure inseparable; $L^*_0$ is a 
compositum of $L$ and
$K^*_0$; and $L^*$ is a compositum of $L^*_0$ and $K^*$. 
Consequently
$\l\mapsto\l|K^*_0$ gives an isomorphism 
$\operatorname{Gal}(K^*,K)\to\operatorname{Gal}(K^*_0,K)$, 
whereas $\l'\mapsto\l'|L^*_0$ 
gives an isomorphism
$$
\operatorname{Gal}(L^*,K^*)@>\d'>>%
\operatorname{Gal}(L^*_0,K^*_0)
$$
such that $\d_0(\im\d')=\im\d$; therefore by the next to 
last display we
conclude that
$$
\im\d\triangleleft\operatorname{Gal}(L,K)\quad\text{and}%
\quad
\operatorname{Gal}(L,K)/\im\d\approx%
\operatorname{Gal}(K^*,K)/N
$$
for some 
$N\triangleleft\operatorname{Gal}(K^*,K)$.

For applying to specific situations, here are some
\proclaim{Corollaries of the Refined Extension Principle}
Given any finite algebraic field extension $K'$ of $K$, by 
the Basic
Extension Principle we may regard 
$\operatorname{Gal}(f,K')<$ $\operatorname{Gal}(f,K)$ 
$<S_n$,
and then upon letting $K^*$ be a least normal
extension of $K$ containing $K'$, we have the following.

{\rm(1.1)} There exists 
$N\triangleleft\operatorname{Gal}(K^*,K)$ and
$M\triangleleft\operatorname{Gal}(f,K)$ with 
$M<\operatorname{Gal}(f,K')$ such that
$\operatorname{Gal}(f,K)/M\approx%
\operatorname{Gal}(K^*,K)/N$.

{\rm(1.2)} If $\operatorname{Gal}(K^*,K)$ is solvable, 
then there exists 
$M\triangleleft\operatorname{Gal}(f,K)$ with 
$M<\operatorname{Gal}(f,K')$ such that
$\operatorname{Gal}(f,K)/M$ is solvable.

{\rm(1.3)} If $\operatorname{Gal}(K^*,K)$ is solvable, and 
$\operatorname{Gal}(f,K)=S_n$, and 
$3\ne n\ne 4$, then $\operatorname{Gal}(f,K')=S_n\text{ or 
}A_n$.

{\rm(1.4)} If $\operatorname{Gal}(K^*,K)$ is cyclic and 
$\operatorname{Gal}(f,K)=S_n$, 
then %\newline 
$\operatorname{Gal}(f,K')=S_n\text{ or }A_n$.

{\rm(1.5)} If $\operatorname{Gal}(K^*,K)$ is cyclic of odd 
order 
and $\operatorname{Gal}(f,K)=S_n$, 
then%\newline 
$\operatorname{Gal}(f,K')=S_n$.

{\rm(1.6)} If there is no nonidentity group which is a 
homomorphic 
image of$\operatorname{Gal}
(f,K)$ as well as  $\operatorname{Gal}(K^*,K)$, then 
$\operatorname{Gal}(f,K')=\operatorname{Gal}(f,K)$.

{\rm(1.7)} If $\operatorname{Gal}(f,K)$ is a simple group 
which is not a 
homomorphic image of $\operatorname{Gal}(K^*,K)$, then
$\operatorname{Gal}(f,K')=\operatorname{Gal}(f,K)$.

{\rm(1.8)} If $\operatorname{Gal}(K^*,K)$ is solvable and 
$\operatorname{Gal}(f,K)$ is  
nonabelian simple, then %\newline 
$\operatorname{Gal}(f,K')=\operatorname{Gal}(f,K)$.

{\rm(1.9)} If $\operatorname{Gal}(K^*,K)$ is solvable, and 
$\operatorname{Gal}(f,K)=A_n$, and 
$3\ne n\ne 4$, then %\newline 
$\operatorname{Gal}(f,K')=A_n$.

{\rm(1.10)} If $\operatorname{Gal}(K^*,K)$ is cyclic, and
$\operatorname{Gal}(f,K)=A_n$, and $n=4$, then %\newline 
$\operatorname{Gal}(f,K')=A_n\text{ or }(Z_2)^2$.

{\rm(1.11)} If $\operatorname{Gal}(K^*,K)$ is cyclic, and 
$\operatorname{Gal}(f,K)=A_n$, and $n=3$, then %\newline
$\operatorname{Gal}(f,K')=A_n\text{ or }Z_1$.

{\rm(1.12)} If $\operatorname{Gal}(K^*,K)$ is cyclic of 
order 
nondivisible by $3$ and 
$\operatorname{Gal}(f,K)=A_n$, then 
$\operatorname{Gal}(f,K')=A_n$.

{\rm(1.13)} If $\operatorname{Gal}(K^*,K)$ is solvable and 
$\operatorname{Gal}(f,K)=\operatorname{PSL}(2,q)$ for a 
prime power $q>3$, 
then $\operatorname{Gal}(f,K')=\operatorname{PSL}(2,q)$.

{\rm(1.14)} If $\operatorname{Gal}(K^*,K)$ is cyclic of 
order 
nondivisible by
$\CHAR k$ and 
$\operatorname{Gal}(f,K)
=\operatorname{PSL}(2,q)$ with $1<q=n-1=$ a power of
$\CHAR k$, then 
$\operatorname{Gal}(f,K')=\operatorname{PSL}(2,q)$.
\endproclaim

(1.1) follows by taking $M=\operatorname{Gal}(f,K^*)$ in 
the Refined
Extension Principle. The implication $(1.1)\Rightarrow 
(1.2)$  
follows from the fact that a homomorphic image of a 
finite solvable group is solvable. In view of (1.1), the 
implication 
$(1.2)\Rightarrow (1.3)$ follows from the facts that if $ 
n\ge 5$ then 
$S_n$ is nonsolvable, and $A_n$ and $S_n$ are the only 
nonidentity 
normal subgroups of $S_n$, and there are no
other subgroups of $S_n$ between $A_n$ and $S_n$, whereas if
$n\le 2$ then $A_n$ and $S_n$ are the only subgroups of 
$S_n$. 
In view of (1.1), the implication $(1.3) \Rightarrow (1.4)$
follows from the fact that $S_3$ and $S_4$ are noncyclic;
$A_3$ is the only nonidentity normal subgroup of $S_3$, 
for the
factor group we have $S_3/A_3=Z_2$, and the only 
nonidentity normal
subgroup of $S_4$ other than $A_4$ is the Klein group 
$(Z_2)^2$
consisting of the four permutations (1), (12)(34), 
(13)(24), (14)(23),
and the factor group of $S_4$ by the Klein group is 
isomorphic
to $S_3$.  In view of (1.1), the implication $(1.4) 
\Rightarrow (1.5)$
follows from the fact that $S_n/A_n = Z_2$ or $Z_1$.
The implications $(1.1)\Rightarrow (1.6)\Rightarrow (1.7)$ 
are obvious. 
The implication 
$(1.7)\Rightarrow (1.8)$ follows from the fact that a 
homomorphic image of a 
finite solvable group is solvable. 
In view of (1.1), the implication $(1.8) \Rightarrow 
(1.9)$ follows from
the facts that if $n \ge 5$ then $A_n$ is nonabelian simple,
whereas if $n \leq 2$ then $A_n = Z_1$.  The implication 
$(1.1) \Rightarrow
(1.10)$ follows from the facts that $A_4$ is not cyclic, 
and the Klein
group $(Z_2)^2$ is the only nonidentity normal subgroup of 
$A_4$, and
the factor group by the Klein group is $Z_3$.  The 
implication
$(1.1) \Rightarrow (1.11)$ follows from the fact that $A_3 
= Z_3$.  In
view of (1.1), (1.10), and (1.11), the implication $(1.9) 
\Rightarrow (1.12)$
follows from the fact that $A_4/(Z_2)^2 = Z_3 = A_3$.  The 
implication
$(1.8) \Rightarrow (1.13)$ follows from the fact that 
$\operatorname{PSL}(2,q)$ is 
nonabelian simple for every prime power $q > 3$. Finally, 
in view
of (1.1), and what we have said about $S_3$ and $A_4$, the
implication $(1.13) \Rightarrow (1.14)$ follows from the 
facts
that $\operatorname{PSL}(2,2) = S_3$ and 
$\operatorname{PSL}(2,3) = A_4$.

As a consequence of the Refined Extension Principle we 
have the
%\proclaim
\thm{Substitutional Principle}
Assume that $K$ = the field $k(X)$ of rational functions 
in an
indeterminate $X$ with coefficients in a field $k$, i.e., 
$a_i=a_i(X)\in k(X)$
for $1\le i\le n$. Given any $\v(X)\in k(X)\setminus k$, let
$$
f_{\v}=f_{\v}(Y)= Y^n+a_1(\v(X))Y^{n-1}+\dots+
a_n(\v(X))\in k(X)[Y]=K[Y]
$$ 
and let $K'=k(V)$ where $V$ is an indeterminate.
Then the $k$-homomorphism $X\mapsto\v(V)$ gives an embedding
$K=k(X)\subset k(V)=K'$, and by sending $Y$ to $Y$ it 
gives an embedding
$K[Y]=k(X)[Y]\subset k(V)[Y]=K'[Y]$, and this sends $f$ to 
$f_{\v}$ with 
$X$ changed to $V$. Therefore,
$f_{\v}$ has no multiple roots in any field extension of 
$k(X)$ and,
upon letting $K^*$ = a least normal extension of $k(X)$ 
containing $k(V)$,
by the Basic Extension Principle we may regard
$\operatorname{Gal}(f,K^*)<$ $\operatorname{Gal}(f,k(V))=$ 
$\operatorname{Gal}(f_{\v},k(X))<$
$\operatorname{Gal}(f,k(X))<S_n$, and now by the Refined 
Extension Principle
$\operatorname{Gal}(f,K^*)\triangleleft%
\operatorname{Gal}(f,k(X))$ and 
$\operatorname{Gal}(f,k(X))/\operatorname{Gal}(f,K^*)%
\approx\operatorname{Gal}(K^*,k(X))/N$
for some $N\triangleleft\operatorname{Gal}(K^*,k(X))$.
\footnote{Note that by writing $\v(X)=\v'(X)/\v''(X)$, 
where $\v'(X)$ 
and $\v''(X)$ are coprime nonzero polynomials in $X$ with 
coefficients in $k$,
we have $K^*$ = a splitting field of $\v'(Y)-X\v''(Y)$ 
over $k(X)$, and
if the said splitting field coincides with a root field of
$\v'(Y)-X\v''(Y)$ over $k(X)$, then we have
$\operatorname{Gal}(f,K^*)=\operatorname{Gal}(f,k(V))=%
\operatorname{Gal}(f_{\v},k(X))$.}
%\endproclaim 
\ethm

Here are the corresponding
\proclaim{Corollaries of the Substitutional Principle}
Letting the situation be as in the Substitutional 
Principle, and remembering
that $K^*=$ a finite normal extension of $k(X)$ and 
$\operatorname{Gal}(f_{\v},k(X))<$ 
$\operatorname{Gal}(f,k(X))<S_n$, 
we have the following. 

{\rm(2.1)} There exists 
$N\triangleleft\operatorname{Gal}(K^*,k(X))$ and
$M\triangleleft\operatorname{Gal}(f,k(X))$ with 
$M<\operatorname{Gal}(f_{\v},k(X))$ such that
$\operatorname{Gal}(f,k(X))/M\approx%
\operatorname{Gal}(K^*,k(X))/N$.

{\rm(2.2)} If $\operatorname{Gal}(K^*,k(X))$ is solvable, 
then there exists 
$M\triangleleft\operatorname{Gal}(f,k(X))$ with 
$M<\operatorname{Gal}(f_{\v},k(X))$ such that
$\operatorname{Gal}(f,k(X))/M$ is solvable.

{\rm(2.3)} If $\operatorname{Gal}(K^*,k(X))$ is solvable, 
and $\operatorname{Gal}(f,k(X))=S_n$, 
and $3\ne n\ne 4$, then 
$\operatorname{Gal}(f_{\v},k(X))=S_n\text{ or }A_n$.

{\rm(2.4)} If $\operatorname{Gal}(K^*,k(X))$ is cyclic and 
$\operatorname{Gal}(f,k(X))=S_n$, 
then %\newline 
$\operatorname{Gal}(f_{\v},k(X)) =S_n$ or $A_n$.

{\rm(2.5)} If $\operatorname{Gal}(K^*,k(X))$ is cyclic of 
odd order and 
$\operatorname{Gal}(f,k(X))=S_n$, then %\newline 
$\operatorname{Gal}(f_{\v},k(X))=S_n$.

{\rm(2.6)} If there is no nonidentity group which is a 
homomorphic
image of$\operatorname{Gal}
(f,k(X))$ as well as  $\operatorname{Gal}(K^*,k(X))$, then
$\operatorname{Gal}(f_{\v},k(X))=%
\operatorname{Gal}(f,k(X))$.

{\rm(2.7)} If $\operatorname{Gal}(f,k(X))$ is a simple 
group which is not a 
homomorphic image of %\newline 
$\operatorname{Gal}(K^*,k(X))$, then
$\operatorname{Gal}(f_{\v},k(X))=%
\operatorname{Gal}(f,k(X))$.

{\rm(2.8)} If $\operatorname{Gal}(K^*,k(X))$ is solvable 
and $\operatorname{Gal}(f,k(X))$ is  
nonabelian simple, then %\newline 
$\operatorname{Gal}(f_{\v},k(X))=%
\operatorname{Gal}(f,k(X))$.

{\rm(2.9)} If $\operatorname{Gal}(K^*,k(X))$ is solvable, 
and $\operatorname{Gal}(f,k(X))=A_n$, and
$3\ne n\ne 4$, then %\newline 
$\operatorname{Gal}(f_{\v},k(X))=A_n$.

{\rm(2.10)} If $\operatorname{Gal}(K^*,k(X))$ is cyclic, and
$\operatorname{Gal}(f,k(X))=A_n$, and $n=4$, then %\newline 
$\operatorname{Gal}(f_{\v},k(X))=A_n\text{ or }(Z_2)^2$.

{\rm(2.11)} If $\operatorname{Gal}(K^*,k(X))$ is cyclic, 
and 
$\operatorname{Gal}(f,k(X))=A_n$, and $n=3$, then %\newline
$\operatorname{Gal}(f_{\v},k(X))=A_n\text{ or }Z_1$.

{\rm(2.12)} If $\operatorname{Gal}(K^*,k(X))$ is cyclic of 
order nondivisible by $3$
 and 
$\operatorname{Gal}(f,k(X)) =A_n$, then 
$\operatorname{Gal}(f_{\v},k(X))=A_n$.

{\rm(2.13)} If $\operatorname{Gal}(K^*,k(X))$ is solvable 
and 
$\operatorname{Gal}(f,k(X))=\operatorname{PSL}(2,q)$ for a 
prime power $q>3$, 
then $\operatorname{Gal}(f_{\v},k(X))=%
\operatorname{PSL}(2,q)$.

{\rm(2.14)} If $\operatorname{Gal}(K^*,k(X))$ is cyclic of 
order nondivisible 
by
$\CHAR k$ 
and$\operatorname{Gal}
(f,k(X))=\operatorname{PSL}(2,q)$ with $1<q=n-1=$ a power of
$\CHAR k$, then$\operatorname{Gal}(f_{\v},k(X))=%
\operatorname{PSL}(2,q)$.
\endproclaim

The proof of (2.1) to (2.14) follows from the above proof 
of (1.1) to
(1.14) by changing ``Refined Extension Principle'' to 
``Substitutional
Principle'' and by changing $(1.i)$ to $(2.i)$ for $1\le 
i\le 14$.

Stated in a form more suitable for applying to the specific
equations described earlier, here are some further
\proclaim{Corollaries of the Substitutional Principle}
Letting the situation be as in the Substitutional 
Principle, 
and remembering that $\operatorname{Gal}(f_{\v},k(X))<$ 
$\operatorname{Gal}(f,k(X))<S_n$, 
and assuming that $k$ is algebraically closed and 
$\v(X)=cX^r$ with 
$0\ne c\in k$ and nonzero integer $r$, we have the 
following. 

{\rm(3.1)} There exists 
$M\triangleleft\operatorname{Gal}(f,k(X))$ with 
$M<\operatorname{Gal}(f_{\v},k(X))$ such that
%\newline 
$\operatorname{Gal}(f,k(X))/M$ is a cyclic group whose 
order is nondivisible
by $\CHAR k$ but divides $r$. 

{\rm(3.2)} If $\operatorname{Gal}(f,k(X))=S_n$, then 
$\operatorname{Gal}(f_{\v},k(X))=S_n\text{ or }A_n$.

{\rm(3.3)} If $\operatorname{Gal}(f,k(X))=S_n$ and $\CHAR 
k =2$, then 
$\operatorname{Gal}(f_{\v},k(X))=S_n$.

{\rm(3.4)} If $\operatorname{Gal}(f,k(X))$ is  
nonabelian simple, then 
$\operatorname{Gal}(f_{\v},k(X))=$\newline
$\operatorname{Gal}(f,k(X))$.

{\rm(3.5)} If $\operatorname{Gal}(f,k(X))=A_n$ and $3\ne 
n\ne 4$, then 
$\operatorname{Gal}(f_{\v},k(X))=A_n$.

{\rm(3.6)} If $\operatorname{Gal}(f,k(X))=A_n$ and $n=4$, 
then 
$\operatorname{Gal}(f_{\v},k(X))=A_n\text{ or }(Z_2)^2$.

{\rm(3.7)} If $\operatorname{Gal}(f,k(X))=A_n$ and $n=3$, 
then 
$\operatorname{Gal}(f_{\v},k(X))=A_n\text{ or }Z_1$.

{\rm(3.8)} If $\operatorname{Gal}(f,k(X))=A_n$ and $\CHAR 
k=3$, then 
$\operatorname{Gal}(f_{\v},k(X))=A_n$.

{\rm(3.9)} If 
$\operatorname{Gal}(f,k(X))=\operatorname{PSL}(2,q)$ for a 
prime power $q>3$, 
then$\operatorname{Gal}(f_{\v},k(X))=%
\operatorname{PSL}(2,q)$.

{\rm(3.10)} If 
$\operatorname{Gal}(f,k(X))=\operatorname{PSL}(2,q)$ with 
$1<q=n-1=$ a power of
$\CHAR k$, then $\operatorname{Gal}(f_{\v},k(X))=%
\operatorname{PSL}(2,q)$.
\endproclaim

Namely, 
$\operatorname{Gal}(K^*,K)$, i.e., the Galois group of the 
splitting field of 
$Y^{|r|}-\left(\frac{X}{c}\right)^{r/|r|}$ over $k(X)$ is 
a cyclic group whose order is nondivisible by $\CHAR k$ 
but divides $r$,
and hence (3.1), (3.2), (3.3) follow from (2.1), (2.4), 
(2.5) respectively,
and $(3.i)$ follows from $(2.i+4)$ for $4\le i\le 10$.

Here are still some more
\proclaim{Corollaries of the Substitutional Principle}
Letting the situation be as in the Substitutional 
Principle, and remembering
that $\operatorname{Gal}(f_{\v},k(X))<$ 
$\operatorname{Gal}(f,k(X))<S_n$, 
and assuming that 
$$\v(X)=\v_d(\v_{d-1}(\dotsb(\v_1(X))\dotsb)),$$ 
where $d$
is a positive integer and 
$$\v_i(X)=\frac{\v'_i(X)}{\v''_i(X)}\in k(X)\setminus k$$ 
with coprime 
nonzero members $\v'_i(X)$ and $\v''_i(X)$ of $k[X]$ for 
$1\le i\le d$, and 
upon letting $K^*_i$ to be a splitting field of
$\v'_i(Y)-X\v''_i(Y)$ over $k(X)$ for $1\le i\le d$, we 
have the following. 

{\rm(4.1)} If for each $i$ with $1\le i\le d$ we have that 
there is no nonidentity 
group which is a homomorphic image of 
$\operatorname{Gal}(f,k(X))$ as well as  
$\operatorname{Gal}(K^*_i,k(X))$, then $%
\operatorname{Gal}(f_{\v},k(X))=%
\operatorname{Gal}(f,k(X))$.

{\rm(4.2)} If $\operatorname{Gal}(f,k(X))$ is a simple 
group which is not a homomorphic 
image of %\newline 
$\operatorname{Gal}(K^*_i,k(X))$ for any $i$ with $1\le 
i\le d$, 
then $\operatorname{Gal}(f_{\v},k(X))=%
\operatorname{Gal}(f,k(X))$.

{\rm(4.3)} If $\operatorname{Gal}(K^*_i,k(X))$ is solvable 
for each $i$ with 
$1\le i\le d$, and$\operatorname{Gal}(f,k(X))$ is 
nonabelian simple,
then $\operatorname{Gal}(f_{\v},k(X))=%
\operatorname{Gal}(f,k(X))$.

{\rm(4.4)} If $\operatorname{Gal}(K^*_i,k(X))$ is solvable 
for each $i$ with 
$1\le i\le d$, and$\operatorname{Gal}(f,k(X))=A_n$, and 
$3\ne n\ne 4$, then 
$\operatorname{Gal}(f_{\v},k(X))=A_n$.
\endproclaim

Namely, (4.1), (4.2), (4.3), and (4.4) follow by 
repeatedly applying
(2.6), (2.7), (2.8), and (2.9) respectively.

Although we shall not use it in this paper, here is a 
third basic principle
of computational Galois theory; for a proof see \S61 of 
volume I of
van der Waerden's book \cite{V}, and for some applications 
see the 1958
follow-up \cite{A4} of my 1957 paper.
\proclaim{Basic Homomorphism Principle} 
If $f(Y) \in R_v[Y]$ where $R_v$ is the valuation ring of 
a valuation
$v$ of $K$ and if $\of(Y)$ has no multiple roots in any 
overfield of
$\oK$, where $\oK$ is the residue field of $v$ and $\of(Y) 
\in \oK(Y)$
is obtained by applying the canonical epimorphism $R_v\to 
\oK$ to the 
coefficients of $f(Y)$, then, as a permutation group, the 
Galois group
$\operatorname{Gal}(\of,\oK)$ may be regarded as a 
subgroup of the Galois group 
$\operatorname{Gal}(f,K)$.
\endproclaim

To give some concrete examples of valuations, assume that 
$K=k(x)$ where
$x$ is a transcendental over a field $k$.
Then for every nonconstant irreducible $\f(x)\in k[x]$ we 
get a valuation 
$v_{\f}$ of $k(x)/k$ by taking 
$$v_{\f}\left(\frac{\p(x)\f(x)^l}{\e(x)}\right)=l$$
for any integer $l$ and any nonzero $\p(x)$ and $\e(x)$ in 
$k[x]$ which are 
nondivisible by $\f(x)$; we may call this the $\f(x)=0$ 
valuation of $k(x)/k$;
if $\f(x)=x-c$ with $c\in k$ then 
we may also call it the $x=c$ valuation of $k(x)/k$; note 
that
for any other nonconstant irreducible $\f^*(x)\in k[x]$ we 
have
$v_{\f}=v_{\f^*}\Leftrightarrow \f\text{ and }\f^*$ are 
constant multiples of
each other. In addition to these valuations, there is 
exactly one more 
valuation $v_{\infty}$ of $k(x)/k$ given by taking 
$$v_{\infty}\left(\frac{\p(x)}{\e(x)}\right)=\deg\e(x)-%
\deg\p(x)$$
for all nonzero $\p(x)$ and $\e(x)$ in $k[x]$; we may call 
this the $x=\infty$
valuation of $k(x)/k$.

Let $\l>\m\ge 0$ 
be integers, and consider the polynomial $\x(Z)Y^{\l}+
\y(Z)Y^{\m}$
in indeterminates $Y$ and $Z$ where $\x(Z)$ and $\y(Z)$ 
are nonzero coprime 
polynomials in $Z$ with coefficients in $k$. Let $y$ be an 
element in an
overfield of $K=k(x)$ such that $\x(x)y^{\l}+
\y(x)y^{\m}=0$. Now,
upon letting $\hv_1,\dots,\hv_h$  be
the extensions of $v=v_{\f}$ to $\hK=K(y)$, for $1\le j\le 
h$ we have
$$
\OR(\hv_j:v)v(\x(x)\y(x))
=\left|\OR(\hv_j:v)v\left(\frac{-\y(x)}{\x(x)}\right)\right|
=|\hv_j(y^{\l-\m})|
=(\l-\m)|\hv_j(y)|
$$
and hence
$$
\OR(\hv_j:v)
\equiv 0\left(\frac{\l-\m}{\operatorname{GCD}(\l-\m,v(%
\x(x)\y(x)))}\right)
$$
and therefore by $(\dagger)$ we get
$$
[\hK:K]
\equiv 0\left(\frac{\l-\m}{\operatorname{GCD}(\l-\m,v(%
\x(x)\y(x)))}\right).
\tag$\dagger\dagger$
$$
From this we deduce the following
\proclaim{First Irreducibility Lemma}
Let $\l>\m\ge 0$ 
be integers, and consider the polynomial $\x(Z)Y^{\l}+
\y(Z)Y^{\m}$
in indeterminates $Y$ and $Z$ where $\x(Z)$ and $\y(Z)$ 
are nonzero coprime 
polynomials in $Z$ with coefficients in a field $k$. Let 
$y$ be an element in an
overfield of $k(Z)$ such that $\x(Z)y^{\l}+\y(Z)y^{\m}=0$. 
Assume that there exists a finite number of nonconstant 
irreducible polynomials
$\f_1(Z),\dots,\f_m(Z)$ in $Z$ 
with coefficients in $k$ such that, upon letting $\n_i$
be the largest integer for which $\f_i(Z)^{\n_i}$ divides 
$\x(Z)\y(Z)$ in 
$k[Z]$, we %\newline 
have $\operatorname{GCD}(\l-\m,\n_1,\dots,\n_m)=1$. Then
$[k(y,Z):k(Z)]=\l-\m$, and the polynomial %\newline
 $\x(Z)Y^{\l}+\y(Z)Y^{\m}$
is irreducible in $k(Y)[Z]$.
\footnote{That is, $\x(Z)Y^{\l}+\y(Z)Y^{\m}$ is either a 
nonzero element
of $k(y)$ or a nonconstant irreducible polynomial in $Z$ 
with coefficients in $k(y)$.}
\endproclaim

Namely, by $(\dagger\dagger)$ we see that 
$$
[k(y,Z):k(Z)]
\equiv 0\left(\frac{\l-\m}{\operatorname{GCD}(\l-\m,\n_i)}%
\right)
\quad\text{for }1\le i\le m
$$
and clearly
$$
\frac{\l-\m}{\operatorname{GCD}(\l-\m,\n_1,\dots,\n_m)}=
\operatorname{LCM}\left(\frac{\l-\m}{\operatorname{GCD}(%
\l-\m,\n_1)},\dots,
\frac{\l-\m}{\operatorname{GCD}(\l-\m,\n_m)}\right)
$$
and hence
$$
[k(y,Z):k(Z)]
\equiv 0\left(\frac{\l-\m}{\operatorname{GCD}(\l-\m,\n_1,%
\dots,\n_m)}\right).
$$
Since $\operatorname{GCD}(\l-\m,\n_1,\dots,\n_m)=1$, we get
$$[k(y,Z):k(Z)]\equiv 0(\l-\m)$$ and since 
$y^{\l-\m}=-\y(Z)/\x(Z)$,
we conclude that $$[k(y,Z):k(Z)]=\l-\m$$ and 
the polynomial $Y^{\l-\m}+\y(Z)/\x(Z)$ is irreducible in 
$k(Z)[Y]$.
Therefore by the Gauss Lemma,
\footnote{The Gauss Lemma
says that a nonzero polynomial in indeterminates $Y$ and 
$Z$ 
with coefficients in a field $k$ is irreducible in 
$k[Y,Z]$ if and only if
it is irreducible in $k(Z)[Y]$ and, as a polynomial in 
$Y$, its coefficients
have no nonconstant common factor in $k[Z]$.}
the polynomial $\x(Z)Y^{\l-\m}+\y(Z)$ is irreducible in 
$k[Y,Z]$, and hence
again by the Gauss Lemma, 
the polynomial $\x(Z)Y^{\l-\m}+\y(Z)$ is irreducible in 
$k(Y)[Z]$, and 
therefore the polynomial $\x(Z)Y^{\l}+\y(Z)Y^{\m}$ is 
irreducible in $k(Y)[Z]$.

Let us now convert the above lemma into the following
\proclaim{Second Irreducibility Lemma}
Let $\l>\m\ge 0$, and let 
$\x_{\l}(Y,Z)$ and $\y_{\m}(Y,Z)$ be nonzero homogeneous 
polynomials of 
respective degrees $\l$ and $\m$ in $(Y,Z)$ with 
coefficients in a field $k$.
Assume that the polynomials $\x_{\l}(Y,Z)$ and 
$\y_{\m}(Y,Z)$ 
are regular in $Z$,
\footnote{That is, their degrees in $Z$ coincide with 
their degrees in
$(Y,Z)$.}
and the polynomials $\x_{\l}(1,Z)$ and $\y_{\m}(1,Z)$ have 
no nonconstant 
common factor in $k[Z]$.
Also assume that there exists a finite number of 
nonconstant irreducible 
polynomials $\f_1(Z),\dots,\f_m(Z)$ in $Z$ 
with coefficients in $k$ such that, upon letting $\n_i$ to 
be the largest 
integer for which $\f_i(Z)^{\n_i}$ divides 
$\x_{\l}(1,Z)\y_{\m}(1,Z)$ in 
$k[Z]$, we have 
$\operatorname{GCD}(\l-\m,\n_1,\dots,\n_m)=1$. Then the 
polynomial 
$\x_{\l}(Y,Z)+\y_{\m}(Y,Z)$ is irreducible in $k(Y)[Z]$.
\endproclaim

To prove this, note that $Z\mapsto YZ$ gives a 
$k(Y)$-automorphism of
$k(Y)[Z]$ which sends the polynomial $\x_{\l}(Y,Z)+
\y_{\m}(Y,Z)$ 
to the polynomial $\x_{\l}(1,Z)Y^{\l} +\y_{\m}(1,Z)Y^{\m}$.
By the First Irreducibility Lemma, the second polynomial is
irreducible in $k(Y)[Z]$ and hence so is the first.

As an illustration, let $n>t>1$ be integers such that 
$\operatorname{GCD}(n,t)=1$.
Now if $u$ is any element in an algebraic closure $\bar k$ 
of $k$ such that 
$u^n=1=u^t$, then, since $\operatorname{GCD}(n,t)=1$, we 
must have $u=1$; 
therefore $1$ is the only common root of $Z^n-1$ and 
$Z^t-1$ in $\bar k$. 
Again since $\operatorname{GCD}(n,t)=1$, either $n$ or $t$ 
is nondivisible by the 
characteristic of $k$, and hence either $Z^n-1$ or $Z^t-1$ 
is  devoid of 
multiple roots in $\bar k$. Therefore the polynomials
$(Z^n-1)/(Z-1)$ and $(Z^t-1)/(Z-1)$ have no common root in 
$\bar k$
and at least one of them has no multiple root in $\bar k$,
and hence they have no nonconstant common factor in $k[Z]$ 
and at least one
of them has no nonconstant multiple irreducible factor in 
$k[Z]$; 
since $n>t>1$, we conclude that their product has at least 
one nonconstant 
nonmultiple irreducible factor in $k[Z]$.
\footnote{That is, there exists a nonconstant irreducible 
member of $k[Z]$
which divides the said product but whose square does not 
divide the said 
product.}
By applying the $k$-automorphism $Z\mapsto Z+1$ of $k[Z]$ 
and by multiplying
the second polynomial by $-a$, where $a$ is any nonzero 
element of $k$, we see
that the polynomials $(Z+1)^n-1)/Z$ and $-a[(Z+1)^t-1]/Z$
have no nonconstant common factor in $k[Z]$ and their 
product has at 
least one nonconstant nonmultiple irreducible factor in 
$k[Z]$.
Consequently by taking $(Z+Y)^n-Y^n/Z$ and
$-a[(Z+Y)^t-Y^t]/Z$ for $\x_{\l}(Y,Z)$ and $\y_{\m}(Y,Z)$, 
and by letting $\f_1(Z),\dots,\f_m(Z)$ to be the distinct 
monic nonconstant 
irreducible factors of $\x_{\l(1,Z)}\y_{\m}(1,Z)$ in 
$k[Z]$, and by letting
$\n_i$ to be the largest integer for which 
$\f_i(Z)^{\n_i}$ divides
$\x_{\l}(1,Z)\y_{\m}(1,Z)$ in $k[Z]$, we see that 
$\x_{\l}(Y,Z)$ and $\y_{\m}(Y,Z)$ in $k[Z]$ are nonzero 
homogeneous polynomials
of respective degrees $\l=n-1$ and $\m=t-1$ in $(Y,Z)$ 
with coefficients in
$k$ such that the polynomials
$\x_{\l}(1,Z)$ and $\y_{\m}(1,Z)$ have no nonconstant 
common factor in $k[Z]$, 
and for some $i$ we have $\n_i=1$, and hence 
trivially we have 
$\operatorname{GCD}(\l-\m,\n_1,\dots,\n_m)=1$. Therefore 
by the
Second Irreducibility Lemma we get the following
\proclaim{Third Irreducibility Lemma}
Let $n>t>1$ be integers such that 
$\operatorname{GCD}(n,t)=1$. Let $Y$ and $Z$ be 
indeterminates over a field $k$, and let $0\ne a\in k$. 
Then the polynomial
$((Z+Y)^n-Y^n/Z-a[(Z+Y)^t-Y^t]/Z$ is irreducible in
$k(Y)[Z]$.
\endproclaim

%\centerline{}
%\centerline{{\bf }}
\heading 20. The tilde polynomial and borrowing 
cycles\endheading

Let us now apply the twisted derivative method to the 
polynomial
$$
\tF_{n,s}=\tF_{n,s}(X,Y)=Y^n-aY^t+X^s
$$ 
with $1\le t<n\text{ and }\operatorname{GCD}(n,t)=1$,
where $s$ is a positive integer,
$a$ is a nonzero element in a field $k$ of characteristic 
$p$, 
\footnote{Here $p$ may or may not be zero. Later on we 
shall specialize to the
case of $p\ne 0$; with the further assumption that 
$n\equiv 0(p)$ and
$s\equiv 0 (t)$, the polynomial $\tF_{n,s}$ reduces to the
polynomial $\tF_{n}$ considered in \S11.}
and $X$ and $Y$ are indeterminates over $k$.
We want to calculate the Galois group 
$\tG_{n,s}=\operatorname{Gal}(\tF_{n,s},k(X))$.
\footnote{ In a moment we shall show that $\tF_{n,s}$ has 
no multiple roots
in any overfield of $k(X)$, and hence the Galois group 
$\operatorname{Gal}(\tF_{n,s},k(X))$
is defined.}

Let $\ok$ be an algebraic closure of $k$, and let
$$
d=\frac{n-t}{e}\quad\text{ where }e=\cases 
1&\text{if }p=0,\\
\text{max $p^\m$ with }n-t\equiv 0(p^{\m})&\text{if }p\ne 0.
\endcases
\tag{i}
$$
Now 
$$
Y^n-aY^t=Y^t\prod_{i=1}^{d}(Y-\b_i)^e
\tag{ii}
$$
with pairwise distinct nonzero elements $\b_1,\dots,\b_{d}$ 
in  $\bar k$
and $\operatorname{GCD}(e,t)=1$, and hence by the First 
Irreducibility Lemma we see that 
$\tF_{n,s}$ is irreducible in $\bar k(X)[Y]$, and 
therefore $\tF_{n,s}$ 
{\it is irreducible in} $k(X)[Y]$. 
\footnote{For $s=1$ this also follows by noting that 
$\tF_{n,1}$
is monic of degree 1 in $X$.}

Let us put $x=X$ and let $y$ be a root of $\tF_{n,s}$ in 
an overfield of 
$k(x)$. Then
$$
x^s=ay^t-y^n
\tag{iii}
$$
and hence $y$ is transcendental over $k$. 
Let $E(y,Z)$ be the polynomial in an indeterminate $Z$ 
such that $E(y,Y)$
is the twisted $Y$-derivative of $\tF_{n,s}(x,Y)$ at $y$. 
Then
$$
E(y,Z)=\frac{(Z+y)^n-y^n}{Z}-\frac{a[(Z+y)^t-y^t]}{Z}
$$
and hence by the Third Irreducibility Lemma we see that 
$E(y,Z)$ {\it is
irreducible in} $k(y)[Z]$.

For a moment suppose that $s=1<t$; then by (iii) we have 
$x\in k(y)$, and
hence $E(y,Z)$ is irreducible in $k(x,y)[Z]$, and  
therefore by the Twisted
Derivative Criterion we see that  $\tF_{n,1}(X,Y)$ 
has no multiple roots in any overfield of $k(X)$, and  so
$\operatorname{Gal}
(\tF_{n,1},k(X))$ is defined; since $E(y,Z)$ is 
irreducible in $k(x,y)[Z]$,
it also follows that 
$\operatorname{Gal}(\tF_{n,1},k(X))$ is $2$-transitive. 
Thus, 
{\it if $s=1<t$ then $\tF_{n,1}$ is devoid of multiple 
roots and 
$\tG_{n,1}$ is $2$-transitive}.

Reverting to general $s$ but still assuming $t>1$, 
since $\tF_{n,s}(X,Y)= \tF_{n,1}(X^s,Y)$
and $\tF_{n,1}(X,Y)$ has no multiple roots in any 
overfield of $k(X)$, it 
follows that $\tF_{n,s}(X,Y)$ has no multiple roots in any 
overfield of $k(X)$. 
Thus, {\it if $t>1$ then $\tF_{n,s}$ is devoid of multiple 
roots and 
hence $\tG_{n,s}$ is defined}.

To give a direct proof of $\tF_{n,s}$ being devoid of 
multiple roots, let us
calculate its $Y$-discriminant. 

So first 
recall that the $Y$-discriminant of a monic polynomial 
$f=f(Y)$ of degree
$n>0$ in $Y$ with coefficients in a field $K$ is denoted 
by $\text{Disc}_Y(f)$
and is defined by putting
$$
\text{Disc}_Y(f)=\text{Res}_Y(f,f_Y)=\text{the 
$Y$-resultant of $f$ and $f_Y$},
$$
where $f_Y$ is the (ordinary) $Y$-derivative of $f$. Upon 
letting $m$ be 
the degree of $f_Y$, we note that 
$\text{Res}_Y(f,f_Y)$ is the determinant of a certain $n+
m$ by $n+m$ matrix;
also note that $m=n-1\Leftrightarrow n\not\equiv 0(\CHAR 
K)$; 
finally note that if $f_Y=0$, i.e.,
if $\CHAR K=p\ne 0$ and $f\in K[Y^p]$, then we take 
Disc$_Y(f)=0$.
For calculational purposes, upon letting 
$f=\prod_{i=1}(Y-\a_i)$ 
we observe that
$$\text{Disc}_Y(f)=\prod_{i=1}^nf_Y(\a_i)$$
and upon assuming $f_Y\ne 0$ and upon letting 
$f_Y=\e\prod_{i=1}^m(Y-\e_i)$
we observe that
$$
\text{Disc}_Y(f)=(-1)^{nm}\e^n\prod_{i=1}^mf(\e_i).
$$
Finally, upon letting
$$
\text{Disc}^*_Y(f)=\text{the modified $Y$-discriminant of }
F=\prod\nolimits_{i<j}(\a_i-\a_j)^2
$$
we note that
\footnote{On pages 82--87 of volume I of van der Waerden's 
book \cite{V},
it seems to be wrongly asserted that Disc$^*_Y(f)=$ 
Disc$_Y(f)$. Some authors 
call Disc$^*_Y(f)$ the discriminant of $f$.}
$$
\text{Disc}^*_Y(f)=(-1)^{n(n-1)/2}\text{Disc}_Y(f).
$$
Let us record the following well-known 
\proclaim{Discriminant Criterion}
If $\CHAR K\ne 2$ and $f$ is devoid of multiple roots, 
then\RM:
$\operatorname{Gal}(f,K)\subset 
A_n\Leftrightarrow\text{Disc}^*_Y(f)$ 
is a square in $K$.
\endproclaim

Now for the (ordinary) $Y$-derivative of $\tF_{n,s}$  we 
have
$$
\align
(\tF_{n,s})_{_Y}
&=nY^{n-1}-taY^{t-1}\\
&=\cases
-taY^{t-1}\ne 0&\text{ if }n\equiv 0(p),\\
nY^{t-1}\prod_{i=1}^{n-t}(Y-b\w_i)\ne 0&\text{ if 
}n\not\equiv 0(p),
\endcases
\endalign
$$
with 
$b,\w_1,\dots,\w_{n-t}$ in $\ok$ 
such that $nb^{n-t}=ta$ and
$\prod_{i=1}^{n-t}(Y-\w_i)=Y^{n-t}-1$
and hence
$$
\text{Disc}_Y(\tF_{n,s})
=\cases\!\!\!(-1)^{n(t-1)}(-ta)^n(X^s)^{t-1}&\text{ if 
}n\equiv 0(p),\\
\!(-1)^{n(n-1)}n^n(X^s)^{t-1}\prod_{i=1}^{n-t}
[X^s-n^{-1}(n-t)ab^t\w_i^t]&\text{ if }n\not\equiv 0(p),
\endcases
$$
and therefore, because $\operatorname{GCD}(n-t,t)=1$, we get
$$
\text{Disc}_Y(\tF_{n,s})
=\cases
(-1)^{nt}t^na^nX^{s(t-1)}&\text{ if }n\equiv 0(p),\\
X^{s(t-1)}[n^nX^{s(n-t)}-(n-t)^{n-t}t^ta^n]
&\text{ if }n\not\equiv 0(p),
\endcases
\tag{iv}
$$
and, observing that if $n\equiv 0(p)$ then in $k$ we have 
$n^n=0$ and
$-(n-t)^{n-t}t^t=(-1)^{n-t+1}t^n$ = $(-1)^{nt}t^n$,
\footnote{Note that if $n$ and $t$ are both odd then $n-t+
1$ and $nt$ are both
odd, whereas if one out of $n$ and $t$ is odd and the 
other even then
$n-t+1$ and $nt$ are both even, and finally, since 
$\operatorname{GCD}(n,t)=1$, the
remaining possibility of $n$ and $t$ both being even 
cannot occur.}
we conclude that in both the cases we have
$$
\text{Disc}_Y(\tF_{n,s})=
n^nX^{s(n-1)}-(n-t)^{n-t}t^ta^nX^{s(t-1)}.
$$
Thus always $\text{Disc}_Y(\tF_{n,s})\ne 0$ and hence
$\tF_{n,s}$ is devoid of multiple roots.

If $s=1$ and $n-t\not\equiv 0(p)$ then by (i), (ii), (iii) 
we see that the
valuation $x=0$ of $\ok(x)/\ok$ splits in 
$\ok(x,y)=\ok(y)$ into the 
$n-t+1=d+1$ valuations $y=0,y=\b_1,\dots,y=\b_d$ 
with reduced ramification exponents $t,1,\dots,1$, and 
hence if also 
$t\not\equiv 0(p)$ then by the Cycle Lemma we can find a 
$t$-cycle in
$\operatorname{Gal}(\tF_{n,1},\ok(X))$, and therefore, 
since by the Basic Extension Principle we have
$\operatorname{Gal}(\tF_{n,1},\ok(X))<\operatorname{Gal}(%
\tF_{n,1},k(X))$, we get a 
$t$-cycle in
$\operatorname{Gal}(\tF_{n,1},k(X))$. Thus, {\it if 
$n-t\not\equiv 0(p)$ and $t\not\equiv 0(p)$
then $\tG_{n,1}$ contains a $t$-cycle}.

Let us note the following Corollary of an 1871 Theorem of 
Jordan \cite{J2};
a proof can also be found in 13.3 on page 34 of Wielandt 
\cite{Wi}, and
in 4.4 on page 171 of Volume I of Huppert \cite{HB}.

\proclaim{Jordan's Corollary}
$A_n$ and $S_n$ are the only primitive permutation groups 
of degree $n$ 
containing a $3$-cycle.
\endproclaim

Let us also note the following two  1892 Theorems of 
Marggraff \cite{Mar};
these are given as Theorems 13.5 and 13.8 on pages 35 and 
38 of Wielandt
\cite{Wi} respectively.
\footnote{Apparently, Marggraff's Second Theorem can also 
be found in
\cite{J2}; see Neumann \cite{N}.}

\proclaim{Marggraff's First Theorem}
$A_n$ and $S_n$ are the only primitive permutation groups 
of degree $n$ 
having a $(n-\n)$-point stabilizer, with $1<\n<n/2$, which 
acts 
transitively on the remaining $\n$ symbols.
\endproclaim

\proclaim{Marggraff's Second Theorem}
A primitive permutation group of degree $n$ containing a 
$\n$-cycle,
with $1<\n<n$, is $(n-\n+1)$-transitive. 
\endproclaim

For a moment suppose that $n-t\not\equiv 0(p)$ and 
$1<t\not\equiv 0(p)$.
\footnote{Note that if $n\equiv 0(p)$ then automatically 
$n-t\not\equiv 0(p)$
and $t\not\equiv 0(p)$ because by assumption %\newline 
$\operatorname{GCD}(n,t)=1$.}
Now, in view of what we have proved above, $\tG_{n,1}$ is 
a 2-transitive
permutation group containing a $t$-cycle and hence if 
$t=2$ then obviously
\footnote{A 2-cycle is simply a transposition, and a 
2-transitive 
permutation group containing a transposition must contain 
all transpositions 
and hence must be the symmetric group.} 
$\tG_{n,1}=S_n$; if $t=3$ then by Jordan's Corollary 
\footnote{A 2-transitive permutation group is 
automatically primitive.}
we have $\tG_{n,1}=A_n\text{ or }S_n$; if $t<\frac n2$ 
then by Marggraff's
First Theorem we have $\tG_{n,1}=A_n\text{ or }S_n$; if 
$t<n-4$ then by 
Marggraff's Second Theorem and CTT we have 
$\tG_{n,1}=A_n\text{ or }S_n$; 
if $t=n-4$ 
then automatically
$12\ne n\ne 24$ and hence by 
Marggraff's Second Theorem and CTT we have 
$\tG_{n,1}=A_n\text{ or }S_n$; 
and if $t=n-3$ 
then again automatically
$12\ne n\ne 24$ and
hence if also
$11\ne n\ne 23$, then by 
Marggraff's Second Theorem and CTT we have 
$\tG_{n,1}=A_n\text{ or }S_n$. 
In all these cases, 
if $t$ is even then we must have $\tG_{n,1}=S_n$ because 
$A_n$ cannot have a 
cycle of even length, and if $t$ is odd and $p\ne 2$ then, 
in view of the 
Discriminant Criterion, by (iv) we can unambiguously 
decide between $A_n$ and 
$S_n$; in particular, if $k$ is algebraically closed and 
$n\equiv 0(p)$ and 
$p\ne 2$ and $t$ is odd then we get $A_n$ because in that 
case 
Disc$^*_Y(\tF_{n,1})$ is a square in $k(x)$. 
Again, in all these cases,
assuming $k$ to be algebraically closed, in view of 
Corollaries (3.2) to (3.8)
of the Substitutional Principle, we see that (1) $t\text{ 
is even }\Rightarrow
\tG_{n,1}=S_n\Rightarrow\tG_{n,s}=A_n\text{ or }S_n$; (2)
$t\text{ is odd and }n=4\Rightarrow\tG_{n,1}=A_n\text{ or 
}S_n\text{ with }n=4
\Rightarrow\tG_{n,s} =(Z_2)^2\text{ or }A_n\text{ or 
}S_n$; (3)
$t\text{ is odd and }n>4\text{ and }p=2\Rightarrow\tG_{n,1}
=A_n\text{ or }S_n\text{ with }n>4 \Rightarrow\tG_{n,s} 
=A_n\text{ or }S_n$;
(4) $t\text{ is odd and }n=4\text{ and }p\ne 2\text{ and 
Disc}^*_Y(\tF_{n,1})
\text{ is a square in }
k(x)\Rightarrow\tG_{n,1}=A_n\text{ with }n=4
\Rightarrow\tG_{n,s} =(Z_2)^2\text{ or }A_n$;
(5) $t\text{ is odd and }n>4\text{ and }p\ne 2\text{ and 
Disc}^*_Y(\tF_{n,1})
\text{ is a square in }k(x)\Rightarrow\tG_{n,1}
=A_n$ with $n>4\Rightarrow\tG_{n,s} =A_n$; and finally
(6) $t\text{ is odd and }p\ne 2\text{ and 
Disc}^*_Y(\tF_{n,1})$
is not a square in $k(x)\Rightarrow\tG_{n,1}
=S_n\Rightarrow\tG_{n,s} =A_n\text{ or }S_n$.

Now for a moment suppose that $n\equiv 0(p)$. Then 
$t\not\equiv 0(p)$ 
and hence by (iv)
we see that $x=0$ and $x=\infty$ are the only valuations 
of $k(x)/k$ which are
possibly ramified in $k(x,y)$. Moreover, if $s\equiv 0(t)$ 
then, in view of
(i), (ii), (iii), either by direct reasoning we see that 
the valuation $x=0$
of $k(x)/k$ is unramified in $k(x,y)$, or alternatively, 
first, upon letting
$x^*=x^s$ we see that the valuation $x^*=0$ of $k(x^*)/k$ 
splits in
$k(x^*,y)=k(y)$ into several valuations out of which $y=0$ 
is the only 
valuation which is possibly ramified over $k(x^*)$ and for 
it the reduced
ramification index is $t$ and the residue degree is $1$, 
and now upon 
letting
$$
x'=\cases
x&\text{ if }p=0,\\
x^{s'}\text{ where }s'=
\text{max }p^{\lambda}
\text{ with }s\equiv 0(p^{\lambda})
&\text{ if }p\ne 0, 
\endcases
$$
we see that $x\prime^{s/s'}=x^*$ and 
$s/s'$ is a positive integer which is divisible by $t$
but not divisible by $p$ and hence by MRT
\footnote{MRT=Abhyankar's Lemma = pages 181--186 of 
\cite{A5}.}
we see that the valuation $x'=0$ of $k(x')$ is unramified 
in $k(x',y)$,
and finally from this we deduce that the valuation $x=0$ 
of $k(x)/k$
is unramified in $k(x,y)$.
\footnote{This last deduction follows from the easy to 
prove fact which says 
that if a valuation $v$ of a field $K$ is unramified in a 
finite separable 
algebraic field extension $L$ of $K$ then the unique 
extension of $v$ to a 
finite purely inseparable field extension $K'$ of $K$ is 
unramified in the 
compositum of $L$ and $K'$.}

Thus, if $n\equiv 0(p)$ and $s\equiv 0(t)$, then
$x=\infty$ is the only valuation of $k(x)/k$ which is 
possibly ramified
in $k(x,y)$, and hence 
\footnote{Say by Proposition 1 of \cite{A1}.}
$x=\infty$ is the only valuation of $k(x)/k$ which is 
possibly ramified
in a least Galois extension of $k(x)$ containing $k(x,y)$. 
It only 
remains to note that, by Result 4 on page 841 of the 1957 
paper 
\cite{A3}, as a consequence of the genus formula, every 
member of the algebraic 
fundamental group of the affine line over an algebraically 
closed field 
$k$ of characteristic $p\ne 0$ is a quasi $p$-group; 
therefore in our case, 
if $\tG_{n,s}=A_n\text{ or }S_n$ and $p\ne 2$ then we must 
have 
$\tG_{n,s}=A_n$, because clearly $S_n$ is not a quasi 
$p$-group for any 
$n\ge p>2$;
likewise, if $\tG_{n,s}=(Z_2)^2\text{ or }A_n\text{ or 
}S_n$ with
$n=4$ and $p=2$, then we must have
$\tG_{n,s}=(Z_2)^2\text{ or }S_n$  because $A_4$ is not a 
quasi
$2$-group.

Let us put all this together in the following
%\proclaim
\subheading{Summary about the tilde polynomial}
Let $k$ be an algebraically closed field $k$ of 
characteristic $p\ne 0$, and
consider the polynomial $\tF_n=Y^n-aY^t+X^s$ in 
indeterminates $X$ and
$Y$ over $k$, where $a$ is a nonzero element of $k$, and 
$n,s,t$ are positive 
integers with $\operatorname{GCD}(n,t)=1$ and $t<n\equiv 
0(p)$ and 
$s\equiv 0(t)$. Then $\tF_n$ gives an unramified covering 
of the affine line
over $k$, and for the Galois group 
$\tG_n=\operatorname{Gal}(\tF_n,k(X))$ we have the
following.

(II.1) If $1<t<4$ and $p\neq 2$, then $\tG_n=A_n$. 

(II.2) If $1<t<n-3$ and $p\neq 2$, then $\tG_n=A_n$. 

(II.3) If $1<t=n-3$ and $p\neq 2$ and $11\neq p \neq 23$, 
then $\tG_n=A_n$. 

(II.4) If $1<t<4<n$ and $p=2$, then $\tG_n=A_n$ or $S_n$.

(II.5) If $1<t<n-3$ and $p=2$, then $\tG_n=A_n$ or $S_n$.

Here CT was not used in the proofs of (II.1) and (II.4). 
Moreover, 
\footnote{In view of the above discussion, without using CT 
we see that if 
$1<t<4=n$ and $p=2$, then $\tG_n=(Z_2)^2$ or $S_n$, and 
using CT we 
see that
if $1<t=n-3$ and $\tG_n\ne A_n$, then
either
$n=p=11$ and $\tG_n=M_{11}$, or 
$n=p=23$ and $\tG_n=M_{23}$.}
the 
following special cases of (II.2) and (II.5) were proved 
without using CT.

(II.2*) If $1<t<n/2$ and $p\ne 2$, then $\tG_n=A_n$.

(II.5*) If $1<t<n/2$ and $p=2$, then $\tG_n=A_n$ or $S_n$.
%\endproclaim
\rem{Note}
In the next section we shall consider the unramified 
covering
of the affine line given by the polynomial
$\oF_n=Y^n-XY^t+1\text{ with }n=p+t\text{ and }t\not\equiv 
0(p)$ which was
introduced in \S11; there the calculation of the Galois 
group of $\oF_n$
will be based on the fact that it contains a $p$-cycle 
because the
valuation $X=\infty$ splits into the valuations $Y=0$ and 
$Y=\infty$
with reduced ramification exponents $t$ and $p$ 
respectively. 
\footnote{Actually, by reciprocating the roots of $\oF_n$ 
we shall
get the polynomial $Y^n-XY^p+1$ for which the valuation 
$X=\infty$ splits
into the valuations $Y=0$ and $Y=\infty$ with reduced 
ramification
exponents $p$ and $t$ respectively.}
In the present section we considered the unramified covering
of the affine line given by the polynomial
$\tF_n=Y^n-aY^t+X^s\text{ with }n\equiv 0(p)\text{ and 
}\operatorname{GCD}(n,t)=1
\text{ and }s\equiv 0(t)$ which was also introduced in 
\S11; 
contrary to the Galois group of the  polynomial $\oF_n$,
it is not easy to find any cycle in the Galois group of the 
polynomial $\tF_n$ because now the valuation $X = \infty$ 
has the
valuation $Y = \infty$ as the only extension  and for it 
the 
reduced ramification exponent is $n$ which is a multiple 
of $p$; so what we
did was to ``borrow'' a $t$-cycle by going down to the 
subfield $k(X^s)$
of $k(X)$, i.e., by embedding the Galois group over $k(X)$ 
as a  subgroup of
the Galois group over $k(X^s)$ where the latter does 
contain a $t$-cycle.
This method of borrowing cycles can already be found in 
Hilbert's 1892
paper \cite{H} 
where he ``borrows'' a transposition to embed $A_n$ in 
$S_n$, thereby
constructing $A_n$ coverings of the rational number field 
$\bq$. 
\endrem 
%\centerline{}
%\centerline{{\bf }}
\heading 21. The bar polynomial\endheading

Let us now turn to the polynomial
$$
\oF_{n,q}=Y^n-XY^t+1
$$
with $n=q+t$
and positive integer $t\not\equiv 0(p)$
mentioned in \S11. This is a polynomial in indeterminates 
$X$ and $Y$ 
with coefficients in a field $k$ of characteristic $p\ne 
0$, and $q$ is a 
positive power of $p$.
%\footnote{Note that, since $n=q+t$, we have: 
%$t\not\equiv 
%0(p)\Leftrightarrow\operatorname{GCD}(n,t)=1$.}

The (ordinary) $Y$-derivative of $\oF_{n,q}$ is given by
$$
(\oF_{n,q})_{_Y}=n(Y^{n-1}-XY^{t-1})
$$
and so we get 
$$
\oF_{n,q}=\frac Yn (\oF_{n,q})_Y+1
$$ 
and hence for the $Y$-discriminant of $\oF_{n,q}$ we have
$$
\text{Disc}_Y(\oF_{n,q})=n^n=\text{ a nonzero element of }k.
$$
Therefore this gives an unramified covering of the affine 
line over $k$.
At any rate, $\oF_{n,q}$, as a polynomial in $Y$, is 
devoid of multiple
roots, and so we can talk about its Galois group over 
$k(X)$. We are
interested in calculating this Galois group 
$\oG_{n,q}=\operatorname{Gal}(\oF_{n,q},k(X))$.

Now $\oF_{n,q}$, as a polynomial in $X$\<, is linear and 
in it the coefficient of 
$X$ has no common factor with the terms devoid of $X$. 
Therefore $\oF_{n,q}$
is irreducible in $k[Y][X]$ and hence in $k(X)[Y]$. 
Consequently, $\oG_{n,q}$,
{\it as a permutation group of degree $n$, is transitive}.

Regarding $\oF_{n,q}$ as a polynomial in $Y$ and 
reciprocating its roots we 
get the polynomial
$$
\uh=\uh(Y)=Y^{q+t}-xY^q+1,
$$
where we have put $x=X$. Let $y$ be a root of this 
polynomial in some 
overfield of $k(x)$. Then solving $\uh(y)=0$ for $x$ we 
get the equation
$x=y^t+y^{-q}$.
The $q=p$ case of this equation was really the starting 
point of this paper 
and it was originally obtained by taking $h=c_1=1$ in 
Proposition 1 of the 
1957 paper \cite{A3}. The case of general $q$ can also be 
obtained by taking
$h=q/p\text{ and }c_{q/p}=1\text{ and }c_i=0\text{ for } 
1\le i<q/p$ in that
Proposition. 
By the said Proposition, or directly by looking at the above
equation, we see that (1) the simple transcendental 
extension $k(y)$ of $k$ is
a separable algebraic field extension of the simple 
transcendental extension 
$k(x)$ of $k$ with field
degree $[k(y):k(x)]=q+t$; (2) $v_{\infty}:x=\infty$ is the 
only valuation of 
$k(x)/k$ which is ramified in $k(y)$; (3) $v_{\infty}$ 
splits in $k(y)$ into 
the valuations $w_0:y=0$ and $w_{\infty}:y=\infty$; and 
(4) for the residue 
degrees and reduced ramification exponents we have 
$\od(w_0:v_{\infty})=1$ and $\OR(w_0:v_{\infty})=q$, and 
$\od(w_{\infty}:v_{\infty})=1$ and 
$\OR(w_{\infty}:v_{\infty})=t$.        
Now $\oG_{n,q}=\operatorname{Gal}(\oF_{n,q},k(X))=%
\operatorname{Gal}(\uh,k(x))$, 
and hence in view of 
(3) and (4), by the Cycle Lemma we see that if $k$ is 
algebraically
closed and $q=p$ then $\oG_{n,p}$ contains a $p$-cycle; by 
the Basic 
Extension Principle, as a permutation group, 
$\operatorname{Gal}(\oF_{n,q},\bar k(X))$ is
a subgroup of 
$\operatorname{Gal}(\oF_{n,q},k(X))$ where $\bar k$ is an 
algebraic
closure of $k$; therefore, without assuming $k$ to be 
algebraically
closed, we have that {\it if $q=p$ then $\oG_{n,p}$ 
contains a $p$-cycle.}

Let prime denote the twisted $Y$-derivative at $y$. Then
$$
\align 
\uh'(Y)&=[Y^q(Y^t-x)]'+(1)'\quad  \text{(by linearity)}\\
&=[Y^q(Y^t-x)]'\quad  \text{(because constant}' = 0)\\
&=(Y+y)^q(Y^t-x)'+Y^{q-1}(y^t-x) \quad \text{(by power 
product rule)}\\
&=(Y+y)^q(Y^t)'+Y^{q-1}(y^t-x) \quad \text{(because 
constant}' = 0)\\
&=(Y+y)^q(Y^t)'-y^{-q}Y^{q-1} \quad \text{(because }x=y^t+
y^{-q})
\endalign
$$
and hence by the definition of $(Y^t)'$ we get
$$
\uh'(Y)=[(Y+y)^q]\frac{(Y+y)^t-y^t}{Y}-y^{-q}Y^{q-1}
$$
and therefore
$$
\uh'(0)=ty^{q+t-1}\ne 0.
$$

Let $\ud(Y)$ be the polynomial obtained by reciprocating 
the roots of 
$\uh'(Y)$. Then $\ud(Y)$ is a monic polynomial of degree 
$q+t-1$ in $Y$ with 
coefficients in $k(y)$ and we have
$$
\align
\ud(Y)=&\left(\frac{Y^{q+t-1}}{ty^{q+
t-1}}\right)\uh'\left(\frac 1Y\right)\\
=&\left[\left(\frac{Y^q}
{y^q}\right)\left(\frac 1Y +y\right)^q\right]
\frac{\left(\frac{Y^t}{y^t}\right)\left[\left(\frac 1Y +
y\right)^t-y^t\right]}
{\left(\frac{tY}{y}\right)\left(\frac 1Y\right)}
-\left(\frac{Y^{q+t-1}}{ty^{q+t-1}}\right)
\left(\frac{y^{-q}}{Y^{q-1}}\right)
\\=&[(Y+y^{-1})^q]\frac{y[(Y+
y^{-1})^t-Y^t]}{t}-\frac{Y^t}{ty^{2q+t-1}}.
\endalign
$$

Let $\ul(Z)$ be the polynomial obtained by multiplying the 
roots of $\ud(Z)$ by $ty$. Then 
$\ul(Z)$ is a monic polynomial of degree $q+t-1$ in an 
indeterminate $Z$ 
with coefficients in $k(y)$ and we have
$$
\align
\ul(Z)=&(ty)^{q+t-1}\ud\left(\frac{Z}{ty}\right)\\
=&\left[(ty)^q\left(\frac{Z}{ty}+\frac 1y\right)^q\right]
\frac{(ty)^{t-1}y\left[\left(\frac{Z}{ty}+\frac 1y\right)^t-
\left(\frac{Z}{ty}\right)^t\right]}{t}
-\frac{(ty)^{q+t-1}\left(\frac{Z}{ty}\right)^t}{ty^{2q+
t-1}}\\
=&[(Z+t)^q]\frac{(Z+t)^t-Z^t}{t^2}-
\frac{t^{q-2}Z^t}{y^{q+t}}
\endalign
$$
and hence
$$
\ul(Z)=(Z+t)^q\g(Z)-t^{q-2}y^{-q-t}Z^t,
$$
where $\g(Z)$ is the monic polynomial of degree $t-1$ in 
$Z$ with coefficients
in the prime field $k_p\subset k$ given by
$$
\g(Z)=t^{-2}[(Z+t)^t-Z^t]
=\sum_{i=0}^{t-1}\g_iZ^{t-1-i},
$$
where $\g_i=t^{i-1}\bigl(\stack{t}{i+1}\bigr)\in k_p$
with $\g_0=1$.
Now $\g(0)=t^{t-2}$ and $\g(-t)=(-1)^{t+1}t^{t-2}$ and hence
$\g(0)\ne 0\ne \g(-t)$.

Let $z$ be a root of $\ul(Z)$ in some overfield of $k(y)$. 
Then solving
$\ul(z)=0$ for $y$ we get
$$
y^{q+t}=\frac{t^{q-2}z^t}{(z+t)^q\g(z)}\tag"$(\ast)$"
$$
and so in particular we see that $z$ is transcendental 
over $k$.
Consider the valuation $\k:z=0$ of $k(z)/k$ and let $\z$ 
be an extension of
it to $k(y,z)$. Since $t\ne 0\ne \g(0)$, by the above 
equation we see that
$$
(q+t)\z(y)=\z(y^{q+t})=\z\left(\frac{t^{q-2}z^t}{(z+
t)^q\g(z)}\right)
=\OR(\z:\k)\k\left(\frac{t^{q-2}z^t}{(z+t)^q\g(z)}\right)
=\OR(\z:\k)t.
$$
Since $q+t$ and $t$ are coprime, we must have 
$\OR(\z:\k)\ge q+t$. Therefore
$[k(y,z):k(z)]\ge q+t$ and hence by the above equation for 
$y^{q+t}$ 
we get $[k(y,z):k(z)]= q+t$.

Since $[k(y,z):k(z)]= q+t$, by the above equation for 
$y^{q+t}$ we see that
the polynomial $(z+t)^q\g(z)Y^{q+t}-t^{q-2}z^t$ is 
irreducible in $k(z)[Y]$.
Since the nonzero polynomials %\newline 
$(Z+t)^q\g(Z)$ and $t^{q-2}Z^t$ in $Z$ with
coefficients in $k$ have no nonconstant common factor in 
$k[Z]$, we conclude 
that the polynomial 
$$(Z+t)^q\g(Z)Y^{q+t}-t^{q-2}Z^t$$
is irreducible in 
$k[Y,Z]$. Consequently the polynomial $(Z+t)^q\g(Z)y^{q+
t}-t^{q-2}Z^t$
is irreducible in $k(y)[Z]$, and hence {\it the 
polynomial} $\ul(Z)$
{\it is irreducible in} $k(y)[Z]$. 
\footnote{This also follows from the First Irreducibility 
Lemma. In fact, the
above argument is a special case of the argument used in 
the proof of the said
lemma.}
Therefore {\it the polynomials}
$\ul(Y)$ {\it and} $\ud(Y)$ {\it are irreducible in} 
$k(y)[Y]$,
and hence 
$\operatorname{Gal}(\ul,k(y))= 
\operatorname{Gal}(\ud,k(y))$ = {\it a transitive
permutation group of degree} $n-1$. Since this group is the
one-point stabilizer of $\oG_{n,q}= 
\operatorname{Gal}(\oF_{n,q},k(x))$,
we conclude that: $\oG_{n,q}$ {\it is} 2-{\it transitive}.

As we have noted above, if $q=p$ then
$\oG_{n,p}$ contains a $p$-cycle, and hence by a 1873 
Theorem of Jordan
we see that if also $t>2$, then $\oG_{n,p}=A_n$ or $S_n$. 
In Result 4 on page 
841 of the 1957 paper \cite{A3} we have noted that, as a 
consequence of the 
genus formula, every member of the algebraic fundamental 
group
of the affine line over an algebraically closed
ground field of characteristic $p\ne 0$ is a quasi 
$p$-group.
Also obviously, for $n>1$ and $p\ne 2$, the symmetric 
group $S_n$ is 
not a quasi $p$-group. Therefore {\it if $k$ is 
algebraically
closed and $q=p\ne 2$ and $n=p+t$ with $t\not\equiv 0(p)$ 
and $t>2$, then
$\oG_{n,p}$ = the alternating group $A_n$.}

Here is the 1873 Theorem of Jordan \cite{J3} we spoke of; 
proofs can also 
be found in 13.9 on page 39 of Wielandt \cite{Wi} and in 
3.7 on page 331 of
volume III of Huppert-Blackburn \cite{HB}.

\proclaim{Jordan's Theorem}
If a primitive permutation group $G$ of degree $n=p+t$ 
contains a
$p$-cycle, where $p$ is prime and $t>2$, then $G=A_n\text{ 
or }S_n$.
\footnote{A 2-transitive permutation group is 
automatically primitive.}
\endproclaim

In the case $p=2$, a $p$-cycle is simply a transposition, 
and obviously
$S_n$ is the only 2-transitive 
permutation group of degree $n$ which contains a 
transposition. 
\footnote{A 2-transitive permutation group which contains 
a transposition,
must automatically 
contain all transpositions, and hence must be the 
symmetric group.}
Therefore {\it if $q=p=2$
and $n=p+t$ with $t\not\equiv 0(p)$ and $t>0$, then 
$\oG_{n,p}=$ the
symmetric group $S_n$.}

To throw away a second root of $\oF_{n,q}$, or more 
precisely to throw away
a root of $\ul(Z)$, this time around let prime denote the 
twisted 
$Z$-derivative at $z$. Then
$$
\align
&\ul'(Z)=[(Z+t)^q\g(Z)]' - 
t^{q-2}y^{-q-t}[Z^t]'\quad\text{(by linearity)}\\
&\quad =(Z+z+t)^q\g'(Z)+Z^{q-1}\g(z)-t^{q-2}y^{-q-t}[Z^t]'
\quad\text{(by power product rule)}\\
&\quad =(Z+z+t)^q\g'(Z)+Z^{q-1}\g(z)-\g(z)(z+
t)^qz^{-t}[Z^t]'
\quad\text{(by (*))}\\
&\quad =(Z+z+t)^q\bar\g(Z)+\g(z)Z^{q-1}-\g(z)(z+
t)^qz^{-t}\bar\r_{t-1}(Z),
\endalign
$$
where $\bar\r_j(Z)$ is the monic polynomial of degree 
$j\ge 0$ in $Z$
with coefficients in the prime field $k_p\subset k$ given by
$$
\bar\r_j(Z)=\frac{(Z+z)^{j+1}-z^{j+1}}{Z}
$$
and
$$
\bar\g(Z)=\sum_{j=0}^{t-2}\g_{t-2-j}\bar\r_j(Z).
$$

Let $\uf(W)$ be the monic polynomial of degree $q+t-2$ in 
an indeterminate $W$
with coefficients in $k(z)$ obtained by multiplying the 
roots of $\ul'(W)$
by $z$. Then
$$
\align
\uf(W)=&z^{q+t-2}\ul'\left(\frac Wz\right)\\
=&\left[z^q\left(\frac Wz +z+t\right)^q\right]
\left[z^{t-2}\bar\g\left(\frac Wz \right)\right]
+\g(z)z^{t-1}\left[z^{q-1}\left(\frac Wz 
\right)^{q-1}\right]\\
&-\g(z)(z+t)^qz^{q-t-1}
\left[z^{t-1}\bar\r_{t-1}\left(\frac Wz \right)\right]
\endalign
$$
and hence
$$
\uf(W)=(W+z^2+tz)^q\hat\g(W)+\g(z)z^{t-1}W^{q-1}-
\g(z)(z+t)^qz^{q-t-1}\hat\r_{t-1}(W)
$$
where $\hat\r_j(W)$ is the monic polynomial
of degree $j\ge 0$ in $W$ with coefficients in $k[z]$ 
given by
$$
\hat\r_j(W)=\frac{(W+z^2)^{j+1}-z^{2j+2}}{W}
$$
and 
$$
\hat\g(W)=\sum_{j=0}^{t-2}z^{t-2-j}\g_{t-2-j}\hat\r_j(W).
$$
Note that if $t=1$ then $\hat\g(W)=0$,
\footnote{Because in that case the summation in the above 
expression for
$\hat\g(W)$ is empty.}
whereas if $t>1$ then $\hat\g(W)$ is a monic polynomial of 
degree $t-2$
in $W$ with coefficients in $k[z]$, and hence, in 
particular, if $t=2$
then $\hat\g(W)=1$.

For a moment suppose that $t=1$; 
then $\g(Z)$ is a monic polynomial of degree 
$t-1=0$ in $Z$ and hence $\g(Z)=1$ and therefore $\g(z)=1$; 
also as noted above $\hat\g(W)=0$; 
finally $\hat\r_{t-1}(W)$ is a monic polynomial of degree 
$t-1=0$ in $W$ and 
hence $\hat\r_{t-1}(W)=1$. Thus,
$$
t=1\Rightarrow\left\{\aligned 
&\ul(Z)=(Z+1)^q-y^{-q-1}Z\quad\text{and }\;\\
&\uf(W)=W^{q-1}-(z+1)^qz^{q-2}.
\endaligned\right.\tag1*
$$ 

Next, for a moment suppose that $t=2$; 
then 
$$\g(Z)=2^{-2}[(Z+2)^2-Z^2]=Z+1$$
and hence $\g(z)=z+1$; 
also as noted above $\hat\g(W)=1$ and clearly $(W+z^2+
tz)^q=W^q+(z+2)^qz^q$; 
finally $\hat\r_{t-1}(W)$ is the monic polynomial of 
degree 1 in $W$ 
given by $\hat\r_{t-1}(W)=W+2z^2$ and hence 
$$
\split
\uf(0)=&(z+2)^qz^q-(z+1)(z+2)^qz^{q-3}(2z^2)\\
=&(z+2)^qz^{q-1}(z-2z-2)
=-(z+2)^{q+1}z^{q-1}.
\endsplit
$$
Thus,
$$
t=2\Rightarrow\left\{\aligned 
&\ul(Z)=(Z+2)^q(Z+1)-2^{q-2}y^{-q-2}Z^2\quad\text{and }\\
&\uf(W)=W^q+(z+1)zW^{q-1}-(z+1)(z+2)^qz^{q-3}W\\
&\phantom{\uf(W)=}-(z+2)^{q+1}z^{q-1}.
\endaligned\right.\tag2*
$$ 

In the case of general $t$, by looking at the equation 
$(\ast)$, we get the
following ``ramification diagram'' for the field extensions 
$k(x)\subset k(y)\subset k(y,z)$, where the square 
bracketed numbers are the 
reduced ramification exponents.
$$
\smalldisplay
\spreadmatrixlines{-0.1in}
\matrix \format \l&\l&\l\\
  &   &\Bigl |\rightarrow \z_i:\psi_i(z)=0[1] \text{ for 
}1\le i\le t^*-1\\
  &   &\Bigl |\\
  &\Bigl |\rightarrow w_{\infty}:y= \infty 
[t]\text{---}&\Bigl |\\
  &\Bigl |&\Bigl |\\
  &\Bigl |&\Bigl |\rightarrow \z_{t^*}:\psi_{t^*}(z)=0[q] \\
  &\Bigl |\\
v_{\infty}:x=\infty\text{---} &\Bigl |\\
  &\Bigl |\\
  &\Bigl |&\Bigl |\rightarrow \z_0:z=0[t] \\
  &\Bigl |&\Bigl |\\
  &\Bigl |\rightarrow w_0:y=0[q]\,\text{---\kern -.16667em --}  &\Bigl |\\
  &    &\Bigl |\\
  &    &\Bigl |\rightarrow \z_{\infty,j}:z=\infty\text{ 
and }
\e_j(y^{\frac{q+
t}{t'}}z^{\frac{q-1}{t'}})=0[\frac{q-1}{t'}]\\ 
\vspace{0.1in}  &    & \quad\ \text{ for } 1 \leq j \leq 
t''. \\
\endmatrix
$$

To explain the top half of the extreme right hand side in 
the above diagram,
first note that, since $t\not\equiv 0(p)$, the 
$Z$-discriminant of 
$(Z+t^{-1})^t-t^{-t}\in k[Z]$ is a nonzero constant; now 
$0$ is a root of this
monic polynomial and by reciprocating the roots of  
$((Z+t^{-1})^t-t^{-t}/Z$ we get $\g(Z)$; therefore $\g(Z)$ 
is 
free from multiple factors; since $\g(-t)\ne 0$, we 
conclude that 
$\g(Z)=\psi_1(Z)\psi_2(Z)\dotsb\psi_{t^*-1}(Z)$, where, 
upon letting 
$\psi_{t^*}(Z)=Z+t$, we have that 
$\psi_1(Z),\psi_2(Z),\dots,\psi_{t^*}(Z)$ 
are pairwise distinct nonconstant monic 
irreducible polynomials in $Z$ with coefficients in $k$; 
let $\k_i$ be the
valuation of $k(z)/k$ given by $\psi_i(z)$; since 
$\operatorname{GCD}(q+t,q)=1$,
by $(\ast)$ we
see that $\k_i$ has a unique extension $\z_i$ to $k(y,z)$, 
and
$\z_1,\z_2,\dots,\z_{t^*}$ are exactly all the extensions 
of $w_{\infty}$ to 
$k(y,z)$, and we have $\bar r(\z_i:\k_i)=q+t$ for $1\le 
i\le t^*$,
$\bar r(\z_i:w_{\infty})=1$ for $1\le i\le t^*-1$,
and $\bar r(\z_{t^*}:w_{\infty})=q$. 

Turning to the bottom half, since $\operatorname{GCD}(q+
t,t)=1$, by $(\ast)$ we see that the
valuation $\k_0:z=0$ of $k(z)/k$ has a unique extension 
$\z_0$ to $k(y,z)$;
since $t\not\equiv 0(p)$, upon letting 
$t'=\operatorname{GCD}(q+t,q-1)$ we have
$Z^{t'}-t^{q-2}=\e_1(Z)\e_2(Z)\dots\e_{t''}(Z)$, where
$\e_1(Z),\e_2(Z),\dots,\e_{t''}(Z)$ are pairwise distinct 
nonconstant monic
irreducible polynomials in $Z$ with coefficients in $k$; 
now by $(\ast)$
we see that the valuation $\k_{\infty}:z=\infty$ of 
$k(z)/k$ has $t''$
extensions 
$\z_{\infty,1},\z_{\infty,2},\dots,\z_{\infty,t''}$ to 
$k(y,z)$
which are characterized by saying that
$\z_{\infty,j}(\e_j(y^{(q+t)/t'}z^{(q-1)/t'}))>0$ for
$1\le j \le t''$, and moreover
$\z_0,\z_{\infty,1}, \z_{\infty,2},\dots,\z_{\infty,t''}$ 
are exactly all the extensions of $w_0$ to $k(y,z)$, and we 
have$\bar r(\z_0:\k_0)=q+t$, $\bar r(\z_0:w_0)=t$, 
$\bar r(\z_{\infty,j}:\k_{\infty})=(q+t)/t'$ for $1\le 
j\le t
''$, and $\bar r(\z_{\infty,j}:w_0)=(q-1)/t'$ for $1\le 
j\le t''$.

Finally note that, since $v_{\infty}$ is the only 
valuation of $k(x)$ which
is ramified in $k(y)$, no valuation of $k(y)$, other than 
$w_0$ and 
$w_{\infty}$, can be ramified in $k(y,z)$, and no 
valuation of $k(y,z)$,
other than 
$\z_1,\z_2,\dots,\z_{t^*},\z_0,\z_{\infty,1}, \z_{\infty,2},
\dots,\z_{\infty,t''}$,
can be ramified in the splitting field of $\uf(W)$ over 
$k(y,z)$.

Now, as noted above, if $t=1$ then $\g(Z)=1$, whereas if 
$t=2$ then 
$\g(Z)=Z+1$, and hence, in connection with the top half of 
the above diagram,
we have that
$$
t=1\Rightarrow \;y^{q+1}=\frac{z}{(z+1)^q},
\text{ and }t^*=1\text{ and }\psi_1(z)=z+1\tag{$1'$}
$$
whereas
$$
t=2\Rightarrow \;
y^{q+2}=\frac{2^{q-2}z^2}{(z+2)^q(z+1)},
\ t^*=2,\ \psi_1(z)=z+1,\ \text{and}
\ \psi_2(z)=z+2.\tag{$2'$}
$$

Concerning the bottom half of the above diagram we note that
$$
\split
t=2\Rightarrow&
\frac{\operatorname{LCM}(\bar r(\z_0:w_0),\bar 
r(\z_{\infty,1}:w_0),\dots,
\bar r(\z_{\infty,t''}:w_0))}{\bar r(\z_0:w_0)}\\
&\qquad =
\cases \frac{q-1}{6}&\text{ if }q\equiv 1(3),\\
\frac{q-1}{2}&\text{ if }q\not\equiv 1(3).\endcases
\endsplit
\tag"$(2'')$"
$$

Recall that a permutation group is said to be {\it 
semiregular}
if its stabilizer, at any point in the permuted set, is the
identity; in our terminology this is equivalent to 
1-antitransitive.
By analogy, let us say that a univariate monic polynomial, 
with coefficients
in some field, is {\it semiregular} over that field if by 
adjoining any
one of its roots we get the splitting field. Note that 
then, assuming the
roots to be distinct, the polynomial is semiregular if and 
only if its
Galois group is semiregular.

For a moment suppose that $t=1$ and $k$ contains all the 
$(q-1)\text{th}$ roots of $1$. Then by $(1^*)$ we see that 
$\uf(W)$ is 
semiregular over $k(y,z)$ and hence its Galois group is 
1-antitransitive; 
since this group is the 2-point stabilizer of the 
2-transitive permutation
group $\oG_{n,q}$, we conclude that $\oG_{n,q}$ is a (2,3) 
group.
Now if $p=2$ then $\z_1((z+1)^qz^{q-2})=q(q+1)$ and 
$\operatorname{GCD}(q-1,q(q+1))=1$, and 
therefore by $(1^*)$ we see that $\uf(W)$ is irreducible 
in $k(y,z)[W]$ and 
its Galois group is a cyclic group of order $q-1$, and 
hence in particular the 
said Galois group is a (1,1) group; consequently, 
$\oG_{n,q}$ is a (3,3) group.
On the other hand, if $p\ne 2$ then by $(1')$ we have 
$(z+1)^qz^{q-2}=z^{q-1}y^{-(q+1)}$ and hence by $(1^*)$ we 
get
$$
\uf(W)=(W^{(q-1)/2}-z^{(q-1)/2}y^{-(q+1)/2})
(W^{(q-1)/2}+z^{(q-1)/2}y^{-(q+1)/2})
$$
and also 
$$
\z_1(z^{(q-1)/2}y^{-(q+1)/2})
=\frac 12\z_1((z+1)^qz^{q-2})=\frac{q(q+1)}{2}$$
and
$$
\operatorname{GCD}
\left(\frac{q-1}{2},\frac{q(q+1)}{2}\right)=1,
$$
and hence the above two factors
of $\uf(W)$ are irreducible in $k(y,z)[W]$ and the Galois 
group of
$\uf(W)$ over $k(y,z)$  is a cyclic group of order 
$\frac{q-1}{2}$.
Thus, {\it if $t=1$ and $k$ contains all the 
$(q-1)\text{th}$ roots of $1$, 
then $\oG_{n,q}$ is a} (2,3) {\it group, and moreover if 
also 
$p=2$ then $\oG_{n,q}$ is a} (3,3) {\it group and its 
order is $(q+1)q(q-1)$, 
whereas if $p\ne 2$ then $\oG_{n,q}$ is not a} (3,3) {\it 
group and its order 
is $\frac{(q+1)q(q-1)}{2}$}.
 
Therefore by an obvious corollary of the 
Zassenhaus-Feit-Suzuki Theorem
we see that {\it if $t=1$ and $k$
contains all the $(q-1)\text{th}$ roots of \RM1\RM, then}
$\oG_{n,q}=\operatorname{PSL}(2,q)$ {\it with the possible 
exception that for $q=7$ we may
have} 
$\oG_{n,q}=\operatorname{A}\Gamma\operatorname{L}(1,8)$.

The said corollary may be formulated in the  following 
manner. For deducing 
this corollary from the Zassenhaus-Feit-Suzuki Theorem, 
the only thing
we need to check is that if the (degree, order) pair 
$(2^l,2^l(2^l-1)l)$
of a Feit Group $\operatorname{A} \Gamma 
\operatorname{L}(1,2^l)$, where $l$ is a prime, equals
$(q+1, \frac{(q+1)q(q-1)}{2})$ then $l= 
(2^l-2)/2=2^{l-1}-1$ and
by direct calculation we see that this is not possible for 
$l=2$, but
is possible for $l=3$, and is never possible for $l \geq 
4$ because
then $2^{l-1}-1=1+2+2^2+ \dots + 2^{l-2} \geq 1+2(l-2)=l+
(l-3) > l$.

\proclaim{Corollary of the Zassenhaus-Feit-Suzuki-Theorem}
If $p=2$, i.e., if $q$ is a positive power of $2$, then 
$\operatorname{PSL}(2,q)$ is
the only $(3,3)$ group of degree $q+1$.  If $p \neq 2$, 
i.e., if $q$
is a positive power of an odd prime $p$, then 
$\operatorname{PSL}(2,q)$ is the only
$(2,3)$ group of degree $q+1$ and 
order $\frac{(q+1)q(q-1)}{2}$ with the exception
that, in the case $q=7$, the group $\operatorname{A} 
\Gamma \operatorname{L}(1,8)$ also satisfies this
description.
\endproclaim

Now, on the one hand, by Result 4 on page 841 of \cite{A3} 
we know that,
in case $k$ is algebraically closed,  
$\oG_{n,q}$ is a quasi $p$-group, and, on the other hand,
for the group $\operatorname{A} \Gamma 
\operatorname{L}(1,8) = \operatorname{GF}(8) \rtimes 
\Gamma \operatorname{L}(1,8)$
we have $\Gamma \operatorname{L}(1,8) = 
\operatorname{GL}(1,8) \rtimes \Aut 
(\operatorname{GF}(8))$ and hence
$\Aut (\operatorname{GF}(8)) = Z_3$ is a homomorphic image 
of $\operatorname{A} \Gamma \operatorname{L}(1,8)$ and
therefore $\operatorname{A} \Gamma \operatorname{L}(1,8)$ 
is not a quasi $7$-group.  
\footnote{A quasi $p$-group can be characterized as a 
finite group having
no homomorphic image whose order is prime to $p$ and 
greater than 1.}
Consequently, in case $k$ is algebraically closed, the 
said exception cannot 
occur and hence {\it if $k$ is algebraically closed and 
$t=1$, then} 
$\oG_{n,q}=\operatorname{PSL}(2,q) = 
\operatorname{PSL}(2,n-1).$

Before proceeding further, let us record the following

\proclaim{Fourth Irreducibility Lemma}
Remember that $p$ is a prime number, and $t$ is a  
positive integer.
\footnote{In this lemma, $t$ is allowed to be divisible by 
$p$.}
For $0\le i\le t$, let $f_i(Y)$ and $g_i(Y)$ be monic 
polynomials of
degree $p+t-i$ in $Y$ with coefficients in a field $K_i$, 
such that
$g_i(Y)$ is obtained from $f_i(Y)$ by making one or more 
of the three 
operations of multiplying the roots by a nonzero quantity 
in $K_i$, 
reciprocating the roots, 
\footnote{In case of reciprocating the roots, we of course 
assume that zero 
is not a root.}
and decreasing the roots by a quantity in $K_i$. Assume 
that, for 
$1\le i\le t$, the polynomial $f_i(Y)$ is obtained by 
throwing away a root 
$\a_i$ of $g_{i-1}(Y)$, and for the field $K_i$ we have 
$K_i=K_{i-1}(\a_i)$. 
Also assume that $f_0(Y)$ has no multiple roots.
\footnote{In any field extension of $K_0$.}
Finally assume that there exist valuations $\hat u^{(i)}$ 
and $u^{(i)}$ 
of $K_i$ for $0\le i\le t$ such that $\hat 
u^{(0)}=u^{(0)}$\RM; 
for $1\le i\le t$, the valuation $\hat u^{(i)}$ is an 
extension of 
$\hat u^{(i-1)}$ to $K_i$\RM; for $1\le i\le t$, the 
valuation $u^{(i)}$ is an 
extension of $u^{(i-1)}$ to $K_i$\RM; for every $i$ with 
$1\le i\le t$ we have 
$\bar r(\hat u^{(i)}:\hat u^{(i-1)})\not\equiv 0(p)$\RM; 
and for some $j$ with 
$1\le j\le t$ we have $\bar r(u^{(j)}:u^{(j-1)})\equiv 
0(p)$. 

Let $u'$ be an extension of $u^{(t)}$ to a splitting field 
$K'$
of $g_t(Y)$ over $K_t$, and let 
$r^*=\prod_{i=1}^t\bar r(u^{(i)}:u^{(i-1)})$ 
and $r'=\bar r(u':u^{(t)})$. 
Let $u_1^{(t)}=u^{(t)},u_2^{(t)},\dots,u_{\d}^{(t)}$ be 
all the extensions
of $u^{(t-1)}$ to $K_t$, and let 
$$\ds r''=\frac{\operatorname{LCM}(\bar 
r(u_1^{(t)}:u^{(t-1)}),
\dots,r(u_{\d}^{(t)}:u^{(t-1)}))} {\bar 
r(u^{(t)}:u^{(t-1)})}.$$  
Finally let $\bar u^{(0)},\dots,\bar u^{(t)},\bar u'$ be 
any valuations of
$K_0,\dots,K_t,K'$ respectively, with $\bar 
u^{(0)}=u^{(0)}$, such that
$\bar u^{(i)}$ is an extension of $\bar u^{(i-1)}$ for 
$1\le i\le t$, and
$\bar u'$ is an extension of $\bar u^{(t)}$ for $1\le i\le 
t$, and let
$\bar r^*=\prod_{i=1}^t\bar r(\bar u^{(i)}:\bar 
u^{(i-1)})$ and
$\bar r'=\bar r(\bar u':\bar u^{(t)})$.

Then the polynomial $g_t(Y)$ is irreducible in $K_t[Y]$, 
and we have
$|\operatorname{Gal}(g_t,K_t)| \equiv 0(r')$ and $r'\equiv 
0(r'')$ and 
$\bar r^*\bar r'=r^*r'$. Moreover, if $g_{t-1}(Y)$ is 
irreducible in
$K_{t-1}[Y]$, and the residue field of $u^{(0)}$
is an algebraically closed field of the same 
characteristic as $K_0$, and
$\bar r(u_j:u^{(t-1)})\not\equiv 0(\CHAR K_0)$ for $1\le 
j\le\d$, then
$r'=r''$.
\endproclaim

To see this, take an extension $\hat u'$ of $\hat u^{(t)}$ 
to $K'$.
By assumption 
$$
\bar r(u^{(j)}:u^{(j-1)})\equiv 0(p)\quad \text{for some $j$
with $1\le j\le t$}
$$
and clearly $\bar r(u':u^{(0)})=r^*r'$, and hence 
we must have $\bar r(u':u^{(0)})\equiv 0(p)$. Now both the 
valuations
$\hat u'$ and $u'$ are extensions of $\hat 
u^{(0)}=u^{(0)}$ to $K'$
which is a Galois extension of $K_0$, 
\footnote{Because it is a splitting field of $f_0(Y)$ over 
$K_0$.}
and hence $\bar r(\hat u':\hat u^{(0)})=\bar 
r(u':u^{(0)})$, 
and therefore $\bar r(\hat u':\hat u^{(0)})\equiv 0(p)$. 
By assumption $\bar r(\hat u^{(i)}:\hat 
u^{(i-1)})\not\equiv 0(p)$ 
for $1\le i\le t$ and clearly $\bar r(\hat u':\hat u^{(0)})=
\bar r(\hat u^{(1)}:\hat u^{(0)})\dots \bar r(\hat 
u^{(t)}:\hat u^{(t-1)}) 
\bar r(\hat u':\hat u^{(t)})$, and hence we must have
$\bar r(\hat u':\hat u^{(t)})\equiv 0(p)$. 
By assumption $K'$ is a splitting field of $g_t(Y)$ over 
$K_t$ and hence
$|\operatorname{Gal}(g_t,K_t)|$ must be divisible by $\bar 
r(\hat u':\hat u^{(t)})$,
and therefore $|\operatorname{Gal}(g_t,K_t)|\equiv 0(p)$; 
since $g_t(Y)$ is a monic 
polynomial of degree $p$ in $Y$ with coefficients in 
$K_t$, we conclude that 
$g_t(Y)$ is irreducible in $K_t[Y]$. Since $K'$ is a 
splitting field of 
$g_t(Y)$ over $K_t$, we also get 
$|\operatorname{Gal}(g_t,K_t)|\equiv 0(r')$.
Now $K'$ is a Galois extension of $K_{t-1}$, and hence 
$r'\equiv 0(r'')$.
Since $K'$ is a Galois extension of $K_0$, and $\bar u'$ 
and $u'$ are
extensions of $\bar u^{(0)}=u^{(0)}$ to $K'$, we also get
$\bar r^*\bar r'=\bar r(\bar u':\bar u^{(0)})=\bar 
r(u':u^{(0)})=r^*r'$.
If $g_{t-1}(Y)$ is irreducible in $K_{t-1}[Y]$ then $K'$ 
is a least Galois
extension of $K_{t-1}$ containing $K_t$, and hence the 
last assertion
follows from Proposition 7 on page 845 of \cite{A3}.

Note that this lemma gives an alternative proof of the 
irreducibility of
$\ul(Y)$ in case of $q=p$ and $t=1$. In fact we have the 
following

\proclaim{Corollary of the Fourth Irreducibility Lemma}
If $k$ is algebraically closed and $F(Y)$ is a monic 
irreducible polynomial of 
degree $p+1$ in $Y$ with coefficients in $k(X)$ such that 
no valuation of 
$k(X)/k$, other than the valuation given by $X=\infty$, is 
ramified in the
root field of $F(Y)$ over $k(X)$, then the Galois group 
$\operatorname{Gal}(F,k(X))$ is $2$-transitive.
\endproclaim

To deduce this from the lemma, it suffices to note that, 
in view of
Proposition 6 on page 835 of \cite{A3}, in the root field 
of $F(Y)$ over
$k(X)$, the valuation $X=\infty$ must split into two 
valuations with
reduced ramification exponents $p$ and $1$ respectively.
\footnote{
Alternatively, in view of Proposition 6 on page 835 of 
\cite{A3}, 
this Corollary can be deduced from the fact that if $G$ is 
a transitive
permutation group of degree $p+1$ of order divisible by 
$p$ then $G$ must
contain a $p$-cycle and hence it must be $2$-transitive.}

More interestingly, in view of the ramification diagram 
and implication
$(2'')$, by the above lemma we see that {\it if $q=p$ and 
$n=p+t$
with $t=2\not\equiv 0(p)$, then
$\uf(W)$ is irreducible in $k(y,z)[W]$ and so $\oG_{n,p}$ 
is} 
3-{\it transitive, and the order of its} 3-{\it point 
stabilizer}
$\operatorname{Gal}(\uf,k(y,z))$ {\it is divisible by 
$\t$, where $\t=\frac{p-1}{6}$ 
or $\t=\frac{p-1}{2}$ according as $p\equiv 1(3)$ or 
$p\not\equiv 1(3)$, 
and hence the order of $\oG_{n,p}$ is divisible by $(p+
2)(p+1)p\t$, and
moreover the reduced ramification exponent of any 
extension of $\z_0$
\RM(resp.\ $\z_1,\z_2,\z_{\infty,1},\dots,\z_{\infty,t''})$
to a splitting field of $\uf(W)$ over $k(y,z)$ is 
divisible by $\t$
\RM(resp.\ $p\t,\t,1,\dots,1)$, and, in case $k$ is 
algebraically closed,
the reduced ramification exponent of any extension of $\z_0$
\RM(resp.\ $\z_1,\z_2,\z_{\infty,1},\dots,\z_{\infty,t''})$
to a splitting field of $\uf(W)$ over $k(y,z)$ equals $\t$
\RM(resp.\ $p\t,\t,1,\dots,1)$.}
\footnote{In the present case, by $(2')$ we have $t^*=2$.}

Therefore by the following Corollary of CTT we see that 
{\it if 
$q=p$ and $n=p+t$ with $t=2\not\equiv 0(p)$, and if 
$\oG_{n,p}$ is
neither equal to $A_n$ nor equal to $S_n$},
\footnote{That is if, as a permutation group, it is not 
isomorphic to either
of these.}
{\it then either $p=7$ and 
$\oG_{n,p}=\operatorname{PSL}(2,8)$, or $p=7$ and 
$\oG_{n,p}=\operatorname{P}\ug\operatorname{L}(2,8)$, or
$p=31$ and 
$\oG_{n,p}=\operatorname{P}\ug\operatorname{L}(2,32)$.} 

\proclaim{Corollary of CTT}
Let $G$ be a \RM3-transitive permutation group of degree 
$n=p+2$ with
an odd prime $p$, such that $G$ is neither equal to $A_n$ 
nor equal to
$S_n$. Then $p=2^\m-1$ for some prime $\m$,
\footnote{Recall that a prime $p$ is called a Mersenne 
prime if it is of
the form $2^{\m}-1$ for some positive integer $\m$, and 
then automatically 
$\m$ is itself a prime number.}
and either $G=\operatorname{PSL}(2,2^{\m})$ or 
$G=\operatorname{P}\ug\operatorname{L}(2,2^{\m})$. 
Moreover, if the order 
of $G$ is divisible by $(p+2)(p+1)p\t$, where 
$\t=\frac{p-1}{6}$
or $\frac{p-1}{2}$ according as $p\equiv 1(3)$ or 
$p\not\equiv 1(3)$, then
either $p=7$ and $G=\operatorname{PSL}(2,8)$, or $p=7$ and 
$G=\operatorname{P}\ug\operatorname{L}(2,8)$, or 
$p=31$ and $G=\operatorname{P}\ug\operatorname{L}(2,32)$. 

\endproclaim

To deduce 
the above Corollary 
%this 
from CTT, by inspection we see that no group listed in 
items (3) to (8) of CTT has degree of the form $p+2$ with 
an odd
prime $p$. Moreover the degree of every group listed in 
items (1) and (2)
of CTT is of the form $\p^{\m}+1$, with prime $\p$ and 
positive integer $\m$,
and if we have  $\p^{\m}+1=p+2$ then, in case of odd $\p$, 
we would get
the contradiction $p=\p^{\m}-1 \equiv 0(2)$.  Therefore 
by CTT we conclude that $p=2^{\m}-1$ for
some prime $\m$, and $G$ is a group between 
$\operatorname{PSL}(2,2^{\m})$ and 
$\operatorname{P} \ug \operatorname{L}(2,2^{\m})$. Since 
$\m$ is prime and $\operatorname{PSL}(2,2^{\m})$ is a 
normal
subgroup of 
$\operatorname{P}\ug\operatorname{L}(2,2^{\m})$ of index 
$\m$, there are no groups
strictly between $\operatorname{PSL}(2,2^{\m})$ and 
$\operatorname{P} \ug \operatorname{L}(2,2^{\m})$. 
Therefore $G=\operatorname{PSL}(2,2^{\m})$ or 
$G=\operatorname{P}\ug\operatorname{L}(2,2^{\m})$. 
Now if the order of $G$ is divisible by $(p+2)(p+1)p\t$, 
then
$(p+2)(p+
1)p\m=|\operatorname{P}\ug\operatorname{L}(2,2^{\m})|\ge 
|G|\ge(p+2)(p+1)p\t$ and hence
$\m\ge\t\ge\frac{p-1}{6}=(2^{\m}-2)/6$ and therefore 
$2^{\m-1}-1\le 3\m$;
on the other hand, if $\m>5$ then 
$2^{\m-1}-1=(1+2+2^2)+(2^3+2^4+\dots+2^{\m-2})\ge 7+
8(\m-4)=3\m+5(\m-5)>3\m$.
Therefore if the order of $G$ is divisible by $(p+2)(p+
1)p\t$, then we must
have $\m=2$ or 3 or 5, i.e., $p=3$ or 7 or 31; now the 
case $p=3$ is ruled 
out because then $\operatorname{PSL}(2,2^{\m})=A_n$ and 
$\operatorname{P}\ug\operatorname{L}(2,2^{\m})=S_n$ with
$\m=2$ and $n=5$; also, in case of $p=31$, we cannot have 
$G=\operatorname{PSL}(2,32)$
because 
$|\operatorname{PSL}(2,32)|=33\cdot 32\cdot 31<33\cdot 
32\cdot 31\cdot 5=(p+2)(p+1)p\t$.  

Once again 
by Result 4 on page 841 of \cite{A3} we know that, in case 
$k$ is algebraically
closed, the group $\oG_{n,q}$ is a quasi $p$-group, and 
obviously, in the case
$p\ne 2$, the group $S_n$ is not a quasi $p$-group because 
it has $Z_2$ 
as a homomorphic image, and likewise, in the case of a 
Mersenne prime
$p=2^{\m}-1$, with a prime $\m$, the group 
$\operatorname{P}\ug\operatorname{L}(2,2^{\m})$ is not a
quasi $p$-group because it has $Z_{\m}$ as a homomorphic 
image.
Therefore by the above italicized conclusion we see that 
{\it
if $k$ is algebraically closed and $q=p$ and $n=p+t$ with 
$t=2\not\equiv 0(p)$, then in the case $p\ne 7$ we have 
$\oG_{n,p}=A_n$,
whereas in the case $p=7$ we have
$\oG_{n,p}=A_n$ or $\operatorname{PSL}(2,n-1)$.}

Now let $\w$ be an element in an overfield of $k(y,z)$ 
such that
$$
\uf(\w)=0,\tag"$(\ast\ast)$"
$$
and let $\uc(U)$ be the monic polynomial of degree $q-1$ 
in an indeterminate 
$U$ with 
coefficients in $k(y,z,\w)$ obtained by throwing away the 
root $\w$ of $\uf(U)$.
Actually, $\uf(U)\in k(z)[U]$ and hence $\uc(U)\in 
k(z,\w)[U]$.

To compute $\uc(U)$, let prime stand for the twisted 
$U$-derivative at $\w$; 
then we have
$\uc(U)=\uf'(U)$, and by the prime power rule we have 
$(U^q)'=U^{q-1}$,
and by the power rule we have $(U^0)'=0$ and $(U)'=1$, and 
by direct 
calculation we have
$$
\align
(U^{q-1})'
=&\frac{(U+\w)^{q-1}-\w^{q-1}}{U}\\
=&\frac{\w^{q-1}\left[-1+\left(1+\frac U\w\right)^q
\left(1+\frac U\w\right)^{-1}\right]}{U}\\
=&\frac{\w^{q-1}\left[-1+\left(1+\frac{U^q}{\w^q}\right)
\left(1-\frac U\w+\frac{U^2}{\w^2}-\frac{U^3}{\w^3}+\dotsb 
\right)\right]}{U}\\
=&\frac{\w^{q-1}\left[-1+\left(1-\frac{U}{\w}
+\frac{U^2}{\w^2}-\dots +
\frac{U^{q-1}}{\w^{q-1}}\right)\right]}{U}\\
=&U^{q-2}-\w U^{q-3}+w^2U^{q-4}-\dots +\w^{q-3}U-w^{q-2},
\endalign
$$
and therefore, in view of $(2^*)$, by linearity we see that
$$
t=2\Rightarrow \left\{\aligned \uc(U)
=&U^{q-1}+(z+1)z[U^{q-2}-\w U^{q-3}+\w^2U^{q-4}-\dots+
\w^{q-3}U]\\
&-(z+1)z[\w^{q-2}+(z+2)^qz^{q-4}].
\endaligned\right.\tag$2^{**}$
$$

For a moment suppose that $k$ is algebraically closed and 
$q=p=7$ and $t=2$.
Since $A_9$ is 4-transitive and $\operatorname{PSL}(2,8)$ 
is sharply 3-transitive, by
the above italicized assertion it follows that 
$\oG_{n,p}\ne A_n\Leftrightarrow 
\oG_{n,p}=\operatorname{PSL}(2,n-1)\Leftrightarrow
\uc(U)$ is reducible in 
$k(y,z,\w)[U]\Leftrightarrow\uc(U)$ completely factors
(into linear factors) in $k(y,z,\w)[U]$. 
To convert the question of reducibility of $\uc(U)$ over 
the field $k(y,z,\w)$
to its reducibility over the simpler field $k(z,\w)$, 
first note that
$k(z,\w)\subset k(y,z,\w)$ and hence $\uc(U)$ is reducible 
in 
$k(z,\w)[U]\Rightarrow\uc(U)$ is reducible in 
$k(y,z,\w)[U]$; 
next, by $(2')$ we see that
$k(y,z,\w)$ is a cyclic extension of $k(z,\w)$ of degree 1 
or 3 or 9,
and hence $\uc(U)$ completely factors in 
$k(y,z,\w)[U]\Rightarrow$ the
degrees of the irreducible factors of $\uc(U)$ in 
$k(z,\w)[U]$ are
(3,3) or (3,1,1,1) or (1,1,1,1,1,1). {\it Therefore if $k$ 
is algebraically 
closed and $q=p=7$ and $n=p+t$ with $t=2$ then 
$\oG_{n,p}\ne A_n\Leftrightarrow 
\oG_{n,p}=\operatorname{PSL}(2,n-1)\Leftrightarrow
\uc(U)$ is reducible in 
$k(y,z,\w)[U]\Leftrightarrow\uc(U)$ completely factors
\RM(into linear factors\/\RM) in 
$k(y,z,\w)[U]\Leftrightarrow\uc(U)$ is reducible in 
$k(z,\w)[U]\Leftrightarrow$ the degrees of the irreducible 
factors of 
$\uc(U)$ in $k(z,\w)[U]$ are $(3,3)$ or $(3,1,1,1)$ or 
$(1,1,1,1,1,1)$.}

To further simplify the question of reducibility of 
$\uc(U)$, we proceed
to eliminate the last two possibilities of $(3,1,1,1)$ and 
$(1,1,1,1,1,1)$
for the degrees of the irreducible factors of $\uc(U)$ in 
$k(z,\w)$.
Afterwards, we shall also indicate how the question of 
reducibility of
$\uc(U)$ in $k(z,\w)[U]$ can be converted to its 
reducibility in
$k_p(z,\w)[U]$.

Recall that $\g(Z)$ is a monic polynomial of degree $t-1$ 
in $Z$ with 
coefficients in the prime field $k_p\subset k$. Let
$$
y^*=\frac {(z+t)^q\g(z)}{z^t}
$$
and
$$
\ul^*(Z)=(Z+t)^q\g(Z)-y^*Z^t\in k_p(y^*)[Z].
$$
Now the element $y^*$ is transcendental over $k$, the 
polynomial $\ul^*(Z)$ 
is irreducible in $k(y^*)[Z]$, the element $z$ is a root 
of $\ul^*(Z)$, 
the polynomial $\ul'(Z)\in k_p(z)[Z]$ is obtained by 
throwing away the root 
$z$ of $\ul^*(Z)$, the polynomial $\uf(W)\in k_p(z)[W]$ is 
obtained by
multiplying the roots of $\ul'(W)$ by $z$, the element 
$\w$ is a root of
$\uf(W)$, the polynomial $\uc(U)\in k_p(z,\w)[U]$ is 
obtained by throwing
away the root $\w$ of $\uf(U)$, and finally {\it if $q=p$ 
and $t=2$ then the
polynomial $\uf(W)$ is irreducible in $k(z)[W]$ 
and hence} 
$\operatorname{Gal}(\ul^*,k(y^*))$ {\it and} 
$\operatorname{Gal}(\ul^*,k_p(y^*))$
{\it are} $2$-{\it transitive permutation groups of 
degree} $p+1$.

For $t=2$ we have $\g(Z)=Z+1$ and so
$$
t=2\Rightarrow\;
y^*=\frac {(z+2)^q(z+1)}{z^2}\;
\text{ and }\;
\ul^*(Z)=(Z+2)^q(Z+1)-y^*Z^2.\tag$2^{***}$
$$

If $q=p$ and $t=2$ then either directly by the above 
expression for $y^*$, 
or indirectly by the 
formulae $(\ast)$ and $(2')$ and the ramification diagram 
of the 
field extension $k(y)\subset k(y,z)$, we see that 
$w^*_0:y^*=0$
and $w^*_{\infty}:y^*=\infty$ are the only valuations of 
$k(y^*)/k$ which are
ramified in $k(z)$, the valuation $w^*_0$ splits into the
valuations $\k_1:z+1=0$ and $\k_2:z+2=0$ of $k(z)/k$ with 
reduced ramification
exponents $\bar r(\k_1:w^*_0)=1$ and $\bar 
r(\k_2:w^*_0)=p$, and
the valuation $w^*_{\infty}$ splits into the valuations 
$\k_0:z=0$ and 
$\k_{\infty}:z=\infty$ of $k(z)/k$ with reduced 
ramification exponents 
$\bar r(\k_0:w^*_{\infty})=2$ and $\bar 
r(\k_{\infty}:w^*_{\infty})=p-1$.
Thus, in case of $q=p$ and $t=2$ we have the special 
ramification diagram
for the field extension $k(y^*)\subset k(z)$ where the 
square bracketed
numbers are the reduced ramification exponents:
$$
\spreadmatrixlines{-0.1in}
\matrix \format \l&\l&\l\\
  &   &\Bigl |\rightarrow \k_1:z+1=0[1] \\
  &   &\Bigl |\\
  & w^*_0:y^*=0\;\text{---\!\!\!\!\!--}&\Bigl |\\
  &&\Bigl |\\
  &&\Bigl |\rightarrow \k_2:z+2=0[p] \\
\vspace{0.5pc}
  &\\
  &\\
  &\\
  &&\Bigl |\rightarrow \k_0:z=0[2] \\
  &&\Bigl |\\
  & w^*_{\infty}:y^*={\infty}\,\text{---}  &\Bigl |\\
  &    &\Bigl |\\
  &    &\Bigl |\rightarrow \k_{\infty}:z=\infty[p-1]. \\
\endmatrix
$$

For a moment suppose that $q=p$ and $t=2$. Since the 
degree of $\uf(W)$ is 
$p$, by the above diagram we see that $\k_1$ has a unique 
extension
$\l_1$ to $k(z,\w)$, for $\l_1$ we have $\bar 
r(\l_1:\k_1)=p$, and the
reduced ramification exponent of any extension of $\l_1$ 
to a splitting field
of $\uc(U)$ over $k(z,\w)$ equals the reduced ramification 
exponent of
any extension of $\k_2$ to the said splitting field. By 
direct calculation
with the expression of $\uf(W)$ given in $(2^*)$, it can 
be shown that
$\k_2$ has 3 extensions to $k(z,\w)$ with reduced 
ramification exponents
$(1,\frac{p-1}{2},\frac{p-1}{2})$. Therefore, say by the 
Fourth Irreducibility
Lemma, the reduced ramification exponent of any extension 
of $\l_1$ to
the said splitting field must be divisible by 
$\frac{p-1}{2}$.
Consequently the degree of some irreducible factor of 
$\uc(U)$ in
$k(z,\w)[U]$ must be $\ge \frac{p-1}{2}$.  
The said direct calculation will not be given here,
\footnote{But it will be given in my forthcoming paper 
\cite{A7}.}
and hence, for deducing the last consequence, at least in 
the case when $k$ is 
algebraically closed, let us indirectly argue  in the 
following manner.

So for a moment suppose that $k$ is algebraically closed 
and $q=p$ and $t=2$. 
Now in view of MRT,
\footnote{MRT = Abhyankar's Lemma =  Pages 181--186 of 
\cite{A5}.}
by implication $(2')$ and by the paragraphs ``Finally note 
that...'' and
``More interestingly...'', we see that $\k_1,\k_2,\k_0$ 
are the only 
valuations of $k(z)/k$ which are ramified in a splitting 
field
of $\uf(W)$ over $k(z)$, and the reduced ramification 
exponent of any
extension of $\k_1$ (resp.\ $\k_2,\k_0$) to the said 
splitting field is
$\frac{p(p-1)}{2}$ (resp.\ $\frac{p-1}{2},\frac{p-1}{2}$). 
Since the degree of 
$\uf(W)$ is $p$, it follows that $\k_1$ has a unique 
extension
$\l_1$ to $k(z,\w)$, for $\l_1$ we have $\bar 
r(\l_1:\k_1)=p$, and the
reduced ramification index of any extension of $\l_1$ to a 
splitting field
of $\uc(U)$ over $k(z,\w)$ is $\frac{p-1}{2}$.
Therefore the degree of some irreducible factor of 
$\uc(U)$ in
$k(z,\w)[U]$ must be $\ge \frac{p-1}{2}$.  

Finally for a moment suppose that $k$ is algebraically 
closed and
$q=p=7$ and $t=2$. Then by the above paragraph and the 
paragraph following
$(2^{**})$ we see that $\uc(U)$ is reducible
in $k(z,\w)[U]\Leftrightarrow$ the degrees of the 
irreducible factors of
$\uc(U)$ in $k(z,\w)[U]$ are (3,3) or (3,1,1,1).
Thus, $\uc(U)$ is reducible in 
$k(z,\w)[U]\Leftrightarrow\operatorname{Gal}(\ul^*,k(y^*))$
is 2-transitive but neither 3-transitive nor sharply 
2-transitive.
Therefore, since $\operatorname{AGL}(1,8)$ is sharply 
2-transitive, by the Special CDT
we conclude that $\uc(U)$ is reducible in 
$k(z,\w)[U]\Leftrightarrow
\operatorname{Gal}(\ul^*,k(y^*))=\operatorname{PSL}(2,7)$
 or $\operatorname{A}\Gamma\operatorname{L}(1,8)$;
by the Zassenhaus-Feit-Suzuki Theorem, both of these are 
(2,3) groups,
and hence $\uc(U)$ is reducible in 
$k(z,\w)[U]\Leftrightarrow$ the degrees of 
the irreducible factors of $\uc(U)$ in $k(z,\w)[U]$ are 
(3,3) 
and they have a common splitting field over $k(z,\w)$ with 
Galois group $Z_3$.
It follows that $\uc(U)$ is reducible in 
$k_p(z,\w)[U]\Leftrightarrow$
the degrees of the irreducible factors of $\uc(U)$ in 
$k_p(z,\w)[U]$ are 
$(3,3)\Leftrightarrow\operatorname{Gal}(\ul^*,k_p(y^*))$ 
is 2-transitive but
neither 3-transitive nor sharply 2-transitive 
$\Leftrightarrow\operatorname{Gal}(\ul^*,k_p(y^*))$ is not 
3-transitive; therefore
by CTT, Special CDT, and the Zassenhaus-Feit-Suzuki
Theorem, we see that
$\uc(U)$ is reducible in $k_p(z,\w)[U]\Leftrightarrow$$%
\operatorname{Gal}(\ul^*,k_p(y^*))=\operatorname{PSL}(2,7)%
\text{ or }\operatorname{A}\Gamma\operatorname{L}(1,8)%
\Leftrightarrow
\operatorname{Gal}(\ul^*,k_p(y^*))\not\in\{S_8,A_8,$$%
\operatorname{AGL}(3,2),\operatorname{PGL}(2,7)\}$.

Let us close this long section with a
%\proclaim
\subheading{Summary about the bar polynomial}
Let $k$ be an algebraically closed field of characteristic 
$p\ne 0$.
Let $q$ be a positive power of $p$ and let $n=q+t$ where 
$t$ is a positive
integer with $t\not\equiv 0(p)$. Let $\oF_{n,q}=Y^n-XY^t+
1\in k[X,Y]$ and let
$\oG_{n,q}=\operatorname{Gal}(\oF_{n,q},k(X))$.
Without using CT we have shown that 
if $t=1$ then $\oG_{n,q}=\operatorname{PSL}(2,q)=%
\operatorname{PSL}(2,n-1)$, 
whereas if $q=p$ and $t>2\ne p$ then $\oG_{n,q}=A_n$, 
and finally if $q=p$ and $t>2=p$ then $\oG_{n,q}=S_n$.
\footnote{Since $\operatorname{PSL}(2,2)=S_3$, it follows 
that if 
$q=p=2$ then $\oG_{n,q}=S_n$.}
Using CT we have shown that if $q=p\ne 7$ and $t=2$ then 
$\oG_{n,q}=A_n$.
Using CT, and referring to 
$(2^*),(2'),(\ast\ast),(2^{**}),(2^{***})$ for
notation and remembering that $z$ is transcendental over 
$k$, we have also 
shown that if $q=p=7$ and $t=2$ then
$\oG_{n,q}\ne A_n\Leftrightarrow\oG_{n,q}=%
\operatorname{PSL}(2,n-1)\Leftrightarrow\uc(U)$
is reducible in $k(y,z,\w)[U]\Leftrightarrow\uc(U)$ 
completely factors 
(into linear factors) in $k(y,z,\w)[U]\Leftrightarrow\uc(U)$
is reducible in $k(z,\w)[U]\Leftrightarrow$ the degrees of 
the irreducible
factors of $\uc(U)$ in $k(z,\w)[U]$ are $(3,3)$ 
and they have a common splitting field over $k(z,\w)$ with 
Galois group 
$Z_3\Leftrightarrow\operatorname{Gal}(\ul^*,k(y^*))=%
\operatorname{A}\Gamma\operatorname{L}(1,8)
\Leftrightarrow\uc(U)$ is reducible in 
$k_p(z,\w)[U]\Leftrightarrow$
the degrees of the irreducible factors of $\uc(U)$ in 
$k_p(z,\w)[U]$ are 
$(3,3)\Leftrightarrow\operatorname{Gal}(\ul^*,k_p(y^*))=%
\operatorname{A}\Gamma\operatorname{L}(1,8)$.
%\endproclaim
\rem{Remark} To establish the above long chain of 
equivalences, we still need to
prove that
$$
\oG_{9,7}=\operatorname{PSL}(2,8)\Rightarrow%
\operatorname{Gal}(\ul^*,k(y^*))\ne
\operatorname{PSL}(2,7)\tag$'$
$$
and
$$
\align
&\operatorname{Gal}(\ul^*,k(y^*))\tag$''$\\
&\qquad=\operatorname{A}\Gamma\operatorname{L}(1,8)
 \Rightarrow\operatorname{Gal}(\ul^*,k_p(y^*))
 \not\in\{S_8,A_8,\operatorname{AGL}(3,2),\\
\shoveright{\operatorname{PGL}(2,7),%
\operatorname{PSL}(2,7)\}.}
\endalign
$$
Now 
$\operatorname{Gal}(\ul,k(y))$ and 
$\operatorname{AGL}(1,8)$ are the 1-point stabilizers of
$\oG_{9,7}$ and $\operatorname{PSL}(2,8)$ respectively, 
and hence to prove $(')$ it
suffices to show that
$$
\operatorname{Gal}(\ul,k(y))=\operatorname{AGL}(1,8)%
\Rightarrow\operatorname{Gal}(\ul^*,k(y^*))\ne
\operatorname{PSL}(2,7)\tag$'''$
$$
Let $f$ be a nonconstant univariate monic polynomial with 
coefficients
in a field $K$ having no multiple roots in any field 
extension of $K$, and
let $K^*$ be a (finite) Galois extension of $K$. By $(2')$ 
and $(2^{***})$
we see that, in the case $q=p=7$ and $t=2$, the field 
$k(y)$ is a Galois 
extension of the field $k(y^*)$ and hence, by taking 
$f=\ul=\ul^*$ and 
$(K,K^*)=(k(y^*),k(y))$, implication $(''')$ follows from 
the implication
$$
\operatorname{Gal}(f,K^*)=\operatorname{AGL}(1,8)\Rightarrow
\operatorname{Gal}(f,K)\ne 
\operatorname{PSL}(2,7).\tag$\sharp$
$$
Also clearly for some finite algebraic field extension 
$k^*$ of $k_p$
we have that 
$\operatorname{Gal}(\ul^*,k^*(y^*))= 
\operatorname{Gal}(\ul^*,k(y^*))$ and $k^*(y^*)$ is a
Galois extension of $k_p(y^*)$ and hence, by taking
$f=\ul^*$ and $(K,K^*)=(k_p(y^*),k^*(y^*))$, implication 
$('')$ follows
from the implication
$$
\split
&\operatorname{Gal}(f,K^*)=\operatorname{A}\Gamma%
\operatorname{L}(1,8)\\
&\qquad\Rightarrow\operatorname{Gal}(f,K)
\not\in\{S_8,A_8,\operatorname{AGL}(3,2),%
\operatorname{PGL}(2,7),\operatorname{PSL}(2,7)\}.
\endsplit
\tag$\sharp\sharp$
$$
\endrem
Now $\operatorname{AGL}(1,8)$ is a nonidentity solvable 
group and $\operatorname{PSL}(2,7)$ is a nonabelian
simple group, and hence $\operatorname{AGL}(1,8)$ is not 
isomorphic to a normal subgroup of 
$\operatorname{PSL}(2,7)$, and therefore $(\sharp)$ 
follows from the Refined
Extension Principle. Likewise, 
$\operatorname{A}\Gamma \operatorname{L}(1,8)$ is a 
nonidentity solvable group but all the nonidentity normal
subgroups of $\{S_8,A_8$, $\operatorname{AGL}(3,2)$, 
$\operatorname{PGL}(2,7)$, $\operatorname{PSL}(2,7)\}$ are 
nonsolvable 
because they respectively contain the nonabelian simple 
group 
$\{A_8,A_8$, $\operatorname{PSL}(3,2)$, 
$\operatorname{PSL}(2,7)$, $\operatorname{PSL}(2,7)\}$ as 
a normal 
subgroup, and hence $\operatorname{A}\Gamma 
\operatorname{L}(1,8)$ cannot be isomorphic to a normal 
subgroup of any one of the groups 
$\{S_8,A_8$, $\operatorname{AGL}(3,2)$, 
$\operatorname{PGL}(2,7),\operatorname{PSL}(2,7)\}$, and 
therefore 
$(\sharp\sharp)$ also follows from the Refined Extension 
Principle. 
%\endrem
\rem{Note}
Assuming $q=p=7$ and $t=2$, in my forthcoming paper 
\cite{A7}, first
by resolving the singularities of the curve $\uf(\w)=0$ in 
the $(z,\w)$-plane
I calculate its genus in terms of infinitely near 
singularities, and
then, seeing that the genus is 2  and hence the curve is 
hyperelliptic, by
means of adjoints I express it as a double covering of the 
line, and
finally, using the resulting square-root parametrization 
of the curve I 
show that the polynomial $\uc(U)$ factors in $k(z,\w)[U]$ 
into
two factors of degree 3, and so I conclude that 
$\oG_{9,7}=\operatorname{PSL}(2,8)$.
Initially I used the square-root parametrization to take 
the ``norm''
of $\uc(U)$ which is a monic polynomial 
$\uc^{\sharp}(T,U)$ of degree 12 in $U$ 
with coefficients which are polynomials in $T$ over the 
prime field 
$k_7$. I used REDUCE and MACSYMA to calculate 
$\uc^{\sharp}(T,U)$; the
largest coefficient degree in $T$ turned out to be 216, 
and it took 10 
pages to print out the exact expressions of the 
coefficients. From what I 
have said above, it follows that 
$\uc(U)$ is reducible in 
$k(z,\w)[U]\Leftrightarrow\uc^{\sharp}(T,U)$ is 
reducible in $k_7[T,U]\Leftrightarrow\uc^{\sharp}(T,U)$ 
factors 
in $k_7[T,U]$ into two polynomials of degree 3 in $U$. 
Although the last two 
tasks are finitistic in nature, the computer algebra 
packages REDUCE
and MACSYMA refused to factor bivariate polynomials over a 
finite field!
So I turned back to hand calculation in a hyperelliptic 
function field, using
MACSYMA only for checking ordinary polynomial operations !!
\endrem
%\centerline{} 
%\centerline{{\bf }}
\heading 22. The roof polynomial and decreasing 
induction\endheading

Let $X$ and $Y$ be indeterminates over an algebraically  
closed field $k$
of characteristic $p\ne 0$, let $0\ne a\in k$, let $s$ and 
$t$ be positive
integers with $t\not\equiv 0(p)$, and consider the 
polynomial
$$
\hF_n=Y^n-aX^sY^t+1\quad \text{with }n=p+t
$$
mentioned in \S11. This  corresponds to the $q=p$ case of 
the more
general polynomial
$$
\hF_{n,q}=Y^n-aX^sY^t+1
$$ 
with $n=q+t$ and $q=$ a positive
power of $p$
which was also mentioned in \S11. The polynomials $\hF_n$ 
and $\hF_{n,q}$
in turn may be obtained by substituting $aX^s$ for $X$ in 
the polynomials
$$
\oF_{n,q}=Y^n-XY^t+1
$$
with $n=q+t\text{ and }q=\text{a positive
power of }p$
and
$$
\oF_{n}=Y^n-XY^t+1\quad \text{with }n=p+t.
$$
As noted in the beginning of \S21, 
the $Y$-discriminant of the polynomials $\oF_n$ and 
$\oF_{n,q}$ is the
nonzero element $n^n$ of $k$, and hence 
the $Y$-discriminant of the polynomials $\hF_n$ and 
$\hF_{n,q}$ is also the
nonzero element $n^n$ of $k$, and therefore each one of 
the four polynomials 
$\oF_{n},\oF_{n,q},\hF_{n},\hF_{n,q}$ gives an unramified 
covering of the 
affine line over $k$.
\footnote{That is, say in view of Proposition 1 of 
\cite{A1}, $X=\infty$ is
the only valuation of $k(X)/k$ which is possibly ramified 
in the splitting
field of $\oF_{n}$ (resp.\ $\oF_{n,q},\hF_{n},\hF_{n,q}$) 
over $K(X)$.}

In the beginning of \S21, we noted that
the polynomials $\oF_{n}$ and $\oF_{n,q}$ are irreducible 
in $k(X)[Y]$, and
considering the Galois groups 
$\oG_{n}=\operatorname{Gal}(\oF_{n},k(X))$ and
$\oG_{n,q} =\operatorname{Gal}(\oF_{n,q},k(X))$, as 
summarized in the Summary and the 
Note at the end of that section, in the rest of that 
section we 
proved the following.
\footnote{As said earlier, my work on this paper started 
when, in September
1988, Serre told me that he could prove (I.1*). As he now 
tells me, his
proof, which also applies to (III.1*), uses ``descending 
Galois theory''
which is different from my method which he calls the 
``ascending''
method.
With Serre's kind permission,
his letter to me,
dated 15 November 1990,
describing his ``descending''
proof, is appended herewith.}

(I.1*) If $t=1$, then $\oG_n=\operatorname{PSL}(2,p)=%
\operatorname{PSL}(2,n-1)$.

(I.2*) If $t=2$ and $p=7$, then $\oG_n=%
\operatorname{PSL}(2,8)=\operatorname{PSL}(2,n-1)$.

(I.3*) If $t=2$ and $p\neq 7$, then $\oG_n=A_n$.

(I.4*) If $t>2$ and $p\neq 2$, then $\oG_n=A_n$.

(I.5*) If $p=2$, then $\oG_n=S_n$.

(III.1*) If $t=1$, then $\oG_{n,q}=%
\operatorname{PSL}(2,q)=\operatorname{PSL}(2,n-1)$.

As noted in the said Summary and Note, on the one hand, 
CT was not used in the proofs of (I.1*), (I.4*), (I.5*), 
and (III.1*), and 
on the other hand, the proof (I.2*) was complete only 
modulo the 
reducibility of the polynomial $\uc(U)$ which will be 
established in
\cite{A7}. 
Actually, CT was used only in the sense that, assuming 
$t=2$, we first 
showed $\oG_n$ to be a $3$-transitive permutation group of
degree $n$ and from this by CTT we deduced that 
$\oG_n=A_n$ or $S_n$.
Now the $3$-transitivity tells us that $|\oG_n|\ge 
n(n-1)(n-2)$ and
hence if $p=3$ then obviously $\oG_n=A_n$ or $S_n$. 
Likewise, for $p=5$
we need not invoke CT because again in that case, say in 
view of the following
elementary theorem which occurs as item IV on page 148 of 
Carmichael's
book \cite{Ca}, 
\footnote{The statement of Charmichael's Theorem given on 
page 154 of Volume I 
of Huppert \cite{HB} says that $A_n$ and $S_n$ are the 
only $l$-transitive 
groups of degree $n$ for which $l>\frac n3$. This seems 
incorrect beacause
for $(l,n)=(4,11)$ we have $4>\frac{11}{3}$ but the 
Mathieu group $M_{11}$ 
is a 4-transitive group of degree 11, and for 
$(l,n)=(5,12)$ 
we have $5>\frac{12}{3}$ but the Mathieu group $M_{12}$ is 
a 
5-transitive group of degree 12.}
the $3$-transitivity directly tells us that
$\oG_n=A_n$ or $S_n$.

\proclaim{Carmichael's Theorem}
$A_n$ and $S_n$ are the only $l$-transitive groups of 
degree $n$ for
which $l>[\frac n3+1]$ where $[\frac n3+1]$ is the  
greatest integer not exceeding $\frac n3+1$.
\endproclaim

By Corollaries (3.5) to (3.10) of the Substitutional 
Principle,
the Galois groups 
$\hG_{n,q}=\operatorname{Gal}(\hF_{n,q},k(X))$ and
$\hG_{n}=\operatorname{Gal}(\hF_{n},k(X))$ have the same 
description as the above
description of the Galois groups 
$\oG_{n}=\operatorname{Gal}(\oF_{n},k(X))$ and 
$\oG_{n,q}=\operatorname{Gal}(\oF_{n,q},k(X))$.
Let us summarize this in the following
%\proclaim
\subheading{Summary about the roof polynomial}
Let $k$ be an algebraically closed field of characteristic 
$p\ne 0$, let
$0\ne a\in k$, and let $s$ and $t$ be positive integers 
with 
$t\not\equiv 0(p)$. Then the polynomials  
$\hF_n=Y^n-aX^sY^t+1 \text{ with }n=p+t$ and
$\hF_{n,q}=Y^n-aX^sY^t+1 \text{ with }n=q+t
\text{ and }q=\text{a positive power of }p$, give 
unramified coverings of the
affine line over $k$,
and for their Galois groups 
$\hG_n=\operatorname{Gal}(\hF_n,k(X))$ and
$\hG_{n,q}=\operatorname{Gal}(\hF_{n,q},k(X))$ we have the 
following.

(I.1) If $t=1$, then $\hG_n=\operatorname{PSL}(2,p)=%
\operatorname{PSL}(2,n-1)$.

(I.2) If $t=2$ and $p=7$, then $\hG_n=%
\operatorname{PSL}(2,8)=\operatorname{PSL}(2,n-1)$.

(I.3) If $t=2$ and $p\neq 7$, then $\hG_n=A_n$.

(I.4) If $t>2$ and $p\neq 2$, then $\hG_n=A_n$.

(I.5) If $p=2$, then $\hG_n=S_n$.

(III.1) If $t=1$, then $\hG_{n,q}=\operatorname{PSL}(2,q)=%
\operatorname{PSL}(2,n-1)$.

Here, CT is not used in the proofs of (I.1), (I.4), (I.5) 
and (III.1);
likewise, it is not used in the proof of (I.3) for $p<7$.
Moreover, referring to the polynomial $\uc(U)$ obtained by 
taking
$q=p=7$ in item (2**) of the \S21, the proof of (I.2) is 
complete only modulo the reducibility of $\uc(U)$ in 
$k(z,\w)$ to 
be established in \cite{A7}.
%\endproclaim

Now what does the ``decreasing induction'' in the title of 
this section
refer to? Roughly speaking, it says that if we can find an 
unramified
covering of the affine line with Galois group $A_n$ then 
we can find one
with Galois group $A_{n-1}$. If this were so without any 
qualification,
then an $A_n$ covering for large $n$ would yield an $A_n$ 
covering for all
smaller $n$. But there is a qualification!  More 
precisely, we have
the following
%\proclaim
\subheading{Method of decreasing induction} 
Assuming that we are working over an algebraically closed 
field $k$ of 
characteristic $p\ne 0$, let there be given an irreducible 
$n$-fold 
unramified covering of the affine line by the affine line, 
\footnote{That is, a covering of the projective line by 
the projective
line, which is unramified over the affine line.}
such that the point 
at infinity splits exactly into two points, say the origin 
and the point at 
infinity, with ramification exponents $p$ and 
$t=n-p\not\equiv 0(p)$, and let 
$G$ be the Galois group of the given covering.
\footnote{That is, $G$ is the Galois group of the 
associated least
Galois covering.}
Now if $G=A_n$ with $n>5$ then we can find an unramified 
covering of the 
affine line with Galois group $A_{n-1}$.  More generally, 
without any 
condition on $n$ or $G$, we can find an unramified 
covering of the affine 
line with Galois group $G^*\triangleleft G_1=$ the 
$1$-point stabilizer
of $G$ such that $G_1/G^*$ is cyclic of order nondivisible 
by $p$.
%\endproclaim

To prove this informally, say the original line is the 
$X$-axis,
and the covering line is the $Y$-axis, and let $C$ be 
the least Galois covering of the $X$-axis ``containing'' 
(or, ``dominating'') the $Y$-axis. Then the Galois
group of $C$ over the $Y$-axis is the one-point stabilizer 
$G_1$ of $G$. 
Moreover, the origin of the $Y$-axis is tamely
ramified in $C$, say with reduced ramification exponent 
$t^*$, 
\footnote{By the argument on pages 843--845 of \cite{A3} 
we see that $t^*$
divides LCM$(t,(p-1)!)$ and hence the origin of the 
$Y$-axis is
tamely ramified in $C$.}
and the point at infinity of 
the $Y$-axis is the only other point which is possibly 
ramified in $C$. 
Now considering the $Y^*$-axis with $a^*Y^{*r}=Y$ with 
$0\ne a^*\in k$ and 
$r\equiv 0(t^*)$, and passing to
the ``compositum'' $C^*$ of $C$ and the $Y^*$-axis, by MRT 
\footnote{MRT = pages 181--186 of \cite{A5}.}
we see that
$C^*$ is an unramified covering of the affine line = the 
$Y^*$-axis, and 
upon letting $G^*$ to be the Galois group of $C^*$ over the
$Y^*$-axis, by Corollary (3.1) of the Substitutional 
Principle we see that
$G^*\triangleleft G_1=$ the $1$-point stabilizer
of $G$, and $G_1/G^*$ is cyclic of order nondivisible by 
$p$.

To apply this method to the polynomial $\oF_n$, from the 
initial material 
\footnote{Up to formula $(\ast)$.}
of \S21 we recall that
\footnote{Until further notice 
$q=p$ and $n=p+t$ with positive integer $t\not\equiv 0(p)$.}
by putting $X=x$ and reciprocating 
the roots of $\oF_n$ we get the polynomial
$$
\uh(Y)=Y^{p+t}+xY^p+1
$$
and by letting $y$ be a root of $\uh(Y)$ in an overfield 
of $k(x)$ we have
$k(x,y)=k(y)$, and by letting $\uh'(Y)$ be the twisted 
$Y$-derivative of 
$\uh(Y)$ at $y$, and letting $\ud(Y)$ be the polynomial 
obtained by 
reciprocating the roots of $\uh'(Y)$, and letting $\ul(Z)$ 
be the polynomial 
obtained by multiplying the roots of $\ud(Z)$ by $ty$, we 
get
$$
\ul(Z)=\g(Z)(Z+t)^p-t^{p-2}y^{-p-t}Z^t,
$$
where
$$
\g(Z)=t^{-2}[(Z+t)^t-Z^t].
$$
Now $\ul(Z)$ is irreducible in $k(y)[Z]$ and has no 
multiple roots in
any overfield of $k(y)$, and upon taking 1-point 
stabilizers, by
(I.3*), (I.4*), and (I.5*) we get the following where CT 
is used only in
the $p>5$ case of (I.3**).

(I.3**) If $t=2$ and $p\ne 7$, then 
$\operatorname{Gal}(\ul(Z),k(y))=A_{n-1}$.

(I.4**) If $t>2$ and $p\ne 2$, then 
$\operatorname{Gal}(\ul(Z),k(y))=A_{n-1}$.

(I.5**) If $p=2$, then 
$\operatorname{Gal}(\ul(Z),k(y))=S_{n-1}$.

Remembering that $s$ is any positive integer and $a$ is 
any nonzero element 
of $k$, we let
$$
\ul_{s,a}(Z)=\g(Z)(Z+t)^p-aY^{-s}Z^t
$$
and we note that if $s\equiv 0(p+t)$ then $\ul_{s,a}(Z)$ 
can be obtained 
from $\ul(Z)$ by replacing $y$ by $a'Y^{s/(p+t)}$, where 
$a'\in k$ is 
such that $t^{p-2}a\prime^{-p-t}=a$. In view of the LCM Theorem,
\footnote{That is, Proposition 7 on page 845 of \cite{A3}.}
by the first ramification diagram 
\footnote{Together with its explanation in the three 
paragraphs following it.}
in \S21, we see that LCM$(p-1,t)$ is divisible by the 
reduced ramification exponent of every extension
of the valuation $y=0$ of $k(y)/k$ to a splitting field of 
$\ul(Z)$ over
$k(y)$, and no valuation of $k(y)/k$, other than the 
valuations
$y=0$ and $y=\infty$, is ramified in the said splitting 
field. Therefore,
if $s\equiv 0((p+t)\operatorname{LCM}(p-1,t))$ then 
by MRT we see that no valuation
of $k(Y)/k$, other than the valuation $Y=\infty$, is 
ramified in a splitting
field of $\ul_{s,a}(Z)$ over $k(Y)$, and obviously 
$\ul_{s,a}(Z)$ 
has no multiple roots in any overfield of $k(Y)$, and in 
view
of Corollaries (3.3), (3.5), and (3.8) of the Substitution 
Principle, by 
(I.3**), (I.4**), and (I.5**) we see that $(3')$ if $t=2$ 
and $p\ne 7$, then 
$\operatorname{Gal}(\ul_{s,a}(Z),k(Y))=A_{n-1}$, whereas 
$(4')$ if $t>2$ and
$p\ne 2$, then 
$\operatorname{Gal}(\ul_{s,a}(Z),k(Y))=A_{n-1}$, and 
finally $(5')$ if
$p=2$, then 
$\operatorname{Gal}(\ul_{s,a}(Z),k(Y))=S_{n-1}$. Thus we 
have the following
where CT is used only in the $p>5$ case of (I.$3'$).

(I.$0'$) If $s\equiv 0((p+t)\operatorname{LCM}(p-1,t))$ 
then no valuation
of $k(Y)/k$, other than the valuation $Y=\infty$, is 
ramified in a splitting
field of $\ul_{s,a}(Z)$ over $k(Y)$, and $\ul_{s,a}(Z)$ 
has no multiple roots in any overfield of $k(Y)$.

(I.$3'$) If $t=2$ and $p\ne 7$ and $s\equiv 0((p+
t)\operatorname{LCM}(p-1,t))$, 
then$\operatorname{Gal}(\ul_{s,a}(Z),k(Y))=A_{n-1}$.

(I.$4'$) If $t>2$ and $p\ne 2$ and $s\equiv 0((p+
t)\operatorname{LCM}(p-1,t))$, 
then$\operatorname{Gal}(\ul_{s,a}(Z),k(Y))=A_{n-1}$.

(I.$5'$) If $p=2$ and $s\equiv 0((p+
t)\operatorname{LCM}(p-1,t))$, then 
$\operatorname{Gal}(\ul_{s,a}(Z),k(Y))=S_{n-1}$.

Actually, the above four assertions remain valid if we 
replace the assumption
that $s\equiv 0((p+t)\operatorname{LCM}(p-1,t))$ by the 
weaker assumption that
$s\equiv 0(\operatorname{LCM}(p-1,t))$. To see this, first 
note that $\ul(Z)$ can be 
obtained from $\ul_{1,1}(Z)$ by replacing $Y$ by 
$t^{2-p}y^{p+t}$; now since
no valuation of $k(y)/k$, other than the valuations $y=0$ 
and $y=\infty$,  
is ramified in a splitting field of $\ul(Z)$ over $k(y)$, 
it follows that
no valuation of $k(Y)/k$, other than the valuations $Y=0$ 
and $Y=\infty$,  
is ramified in a splitting field of $\ul_{1,1}(Z)$ over 
$k(Y)$; moreover
since $\ul(Z)$ is irreducible in $k(y)[Z]$ and has no 
multiple roots in any
overfield of $k(y)$, it follows that $\ul_{1,1}(Z)$ is 
irreducible in 
$k(Y)[Z]$ and has no multiple roots in any overfield of 
$k(Y)$; finally,
in view of Corollary (3.1) the Substitutional Principle, 
by (I.3**) to
(I.5**) we see that (3***) if $t=2$ and $p\ne 7$, then 
$\operatorname{Gal}(\ul_{1,1}(Z),k(Y))=A_{n-1}$ or 
$S_{n-1}$, and (4***) if
$t>2$ and $p\ne 2$, then again 
$\operatorname{Gal}(\ul_{1,1}(Z),k(Y))=A_{n-1}$ or 
$S_{n-1}$, 
and (5***) if $p=2$, then 
$\operatorname{Gal}(\ul_{1,1}(Z),k(Y))=S_{n-1}$. Again, 
$\ul_{s,a}(Z)$ can be obtained from $\ul_{1,1}(Z)$ by 
substituting
$aY^s$ for $Y$ and hence, in view of Corollaries (3.2) to 
(3.8) of the
Substitutional Principle, we get the following where CT is 
used only in
the $p>5$ case of (I.$3''$). 

(I.$0''$) No valuation of $k(Y)/k$, other than the 
valuations $Y=0$ and 
$Y=\infty$, is ramified in a splitting field of 
$\ul_{s,a}(Z)$ over $k(Y)$, 
and $\ul_{s,a}(Z)$ is irreducible
\footnote{The irreducibility of $\ul_{s,a}$ follows from 
the First 
Irreducibility Lemma. Alternatively, it follows from the 
fact that a 
polynomial is irreducible if and only if its Galois group 
is transitive.}
in $k(Y)[Z]$ and has no multiple roots in any overfield of 
$k(Y)$.

(I.$3''$) If $t=2$ and $p\ne 7$, then 
$\operatorname{Gal}(\ul_{s,a}(Z),k(Y))=A_{n-1}$
or $S_{n-1}$.

(I.$4''$) If $t>2$ and $p\ne 2$, then 
$\operatorname{Gal}(\ul_{s,a}(Z),k(Y))=A_{n-1}$ 
or $S_{n-1}$.

(I.$5''$) If $p=2$, then 
$\operatorname{Gal}(\ul_{s,a}(Z),k(Y))=S_{n-1}$.

Now $\g(0)\ne 0\ne t$ and, upon letting $Z'$ to be a root 
of $\ul_{1,1}(Z)$
in an overfield of $k(Y)$, we have
$$
Y=\frac{Z\prime^t}{\g(Z')(Z'+t)^p}
$$
and hence the valuation $Y=0$ of $k(Y)/k$ splits in 
$k(Y,Z')=k(Z')$ into
the valuations $Z'=0$ and $Z'=\infty$ with reduced 
ramification exponents
$t\not\equiv 0(p)$ and $p-1\not\equiv 0(p)$ respectively, 
and therefore by the
LCM Theorem we see that LCM$(p-1,t)$ is divisible by the 
reduced ramification exponent of every extension of the 
valuation $Y=0$ of 
$k(Y)/k$ to a splitting field of $\ul(Z)$ over $k(Y)$, and 
hence by MRT
we see that if $s\equiv 0(\operatorname{LCM}(p-1,t))$ then 
no valuation
of $k(Y)/k$, other than the valuation $Y=\infty$, is 
ramified in a splitting
field of $\ul_{s,a}(Z)$ over $k(Y)$, and therefore by 
Result 4 on page 841
of \cite{A3} we know that 
$\operatorname{Gal}(\ul_{s,a}(Z),k(Y))$ is a quasi 
$p$-group.
Thus, in view of (I.$0''$), (I.$3''$), (I.$4''$), and 
(I.$5''$), we conclude 
with the following where CT is used only in the $p>5$ case 
of (IV.$2'$).

(IV.$0'$) If $s\equiv 0(\operatorname{LCM}(p-1,t))$, then 
no valuation
of $k(Y)/k$, other than the valuation $Y=\infty$, is 
ramified in a splitting
field of $\ul_{s,a}(Z)$ over $k(Y)$, and $\ul_{s,a}(Z)$ is 
irreducible
in $k(Y)[Z]$ and has no multiple roots in any overfield of 
$k(Y)$.

(IV.$1'$) If $t>2$ and $p\ne 2$ and $s\equiv 
0(\operatorname{LCM}(p-1,t))$,
then$\operatorname{Gal}(\ul_{s,a}(Z),k(Y))=A_{n-1}$.

(IV.$2'$) If $t=2$ and $p\ne 7$ and $s\equiv 
0(\operatorname{LCM}(p-1,t))$,
then$\operatorname{Gal}(\ul_{s,a}(Z),k(Y))=A_{n-1}$.

(IV.$4'$) If $p=2$ and $s\equiv 
0(\operatorname{LCM}(p-1,t))$, then 
$\operatorname{Gal}(\ul_{s,a}(Z),k(Y))=S_{n-1}$.

As we have said, CT is used only when $t=2$; now in the 
case $t=2$ we have 
$\g(Z)=Z+1$ and hence $\ul_{s,a}(Z)=(Z+1)(Z+
2)^p-aY^{-s}Z^t$; therefore 
by (I.$0''$) we see that if $t=2$ then 
$\ul_{s,a}(Z)=(Z+1)(Z+2)^p-aY^{-s}Z^t$ is unramified 
outside $Y=0$ and 
$Y=\infty$. For getting hold of a variation of (I.$3'''$) 
without CT, let
us reprove the said unramifiedness in a more general 
context.
So let
$$
E(Z)=(Z+1)(Z+b)^p-YZ^{\t}
$$
with $1<{\t}\le p\text{ and }0\ne b\in k$.
Now for the (ordinary) $Z$-derivative we have
$$
E_Z(Z)=(Z+b)^p-{\t}YZ^{{\t}-1}
$$
and hence
$$
E(Z)=(Z+1)E_Z(Z)+E^*(Z),
$$
where
$$
E^*(Z)=({\t}-1)YZ^{{\t}-1}\left(Z+\frac{{\t}}{{\t}-1}\right)
$$
and therefore for the $Z$-discriminant we have
$$
\align
\text{Disc}_Z(E(Z))
&=\text{Res}_Z(E(Z),E_Z(Z))\\
&=\text{Res}_Z(E^*(Z),E_Z(Z))\\
&=({\t}-1)^pY^pE_Z(0)^{{\t}-1}E_Z\left(\frac{-{\t}}{{%
\t}-1}\right)\\
&=({\t}-1)^pY^pb^{p({\t}-1)}E_Z\left(\frac{-{\t}}{{\t}-1}%
\right)\\
&=({\t}-1)^pY^pb^{p({\t}-1)}(-{\t}Y)\left(\frac{-{\t}}{{%
\t}-1}\right)^{{\t}-1}
\quad\text{ if }b=\frac{{\t}}{{\t}-1}\\
&=(-{\t})^{\t}({\t}-1)^{p-{\t}+1}b^{p({\t}-1)}Y^{p+1}
\quad\text{ if }b=\frac{{\t}}{{\t}-1}\\
\endalign
$$
and hence if $b=\frac{{\t}}{{\t}-1}$ then no valuation of 
$k(Y)/k$, other than
the valuations $Y=0$ and $Y=\infty$, is ramified in a 
splitting field of
$E(Z)$ over $k(Y)$, and $E(Z)$ has no repeated roots in 
any overfield
of $k(Y)$. Therefore, upon remembering that $a$ is any 
nonzero element of
$k$ and $s$ and $t$ are any positive integers with 
$t\not\equiv 0(p)$,  
and upon letting
$$
E_{s,a}(Z)=(Z+1)(Z+b)^p-aY^{-s}Z^t
$$
we get the following.
\footnote{For a while there will be no reference to $n$.}

(IV.0**) If $1<t<p$ and $b=\frac{t}{t-1}$ then no 
valuation of $k(Y)/k$, 
other than the valuations $Y=0$ and $Y=\infty$, is 
ramified in a splitting field
of $E_{s,a}(Z)$ over $k(Y)$, and $E_{s,a}(Z)$ is 
irreducible 
\footnote{The irreducibility of $E_{s,a,b}$ follows from 
the First 
Irreducibility Lemma.}
in $k(Y)[Z]$ and has no multiple roots in any overfield of 
$k(Y)$.

Before proceeding further, let us make note of the following
\proclaim{Alternate Corollary of the Fourth Irreducibility 
Lemma}
If $F(Z)$ is a monic irreducible polynomial of degree $p+
1$ in $Z$ with
coefficients in $k(Y)$ such that some valuation of 
$k(Y)/Y$, say the
valuation $Y=\infty$, splits in a root field of $F(Z)$ 
over $k(Y)$ into
two valuations with reduced ramification exponents $p$ and 
$1$, then
the Galois group $\operatorname{Gal}(F(Z),k(Y))$ is 
$2$-transitive.
\endproclaim

For a moment suppose that $1<t<p$ and $b=\frac{t}{t-1}$. 
Now upon letting
$Z^*$ to be a root of $E_{1,1}(Z)$ in some overfield of 
$k(Y)$, we have
$$
Y=\frac{Z^{*t}}{(Z^*+1)(Z^*+b)^p}
$$
and hence the valuation $Y=\infty$ of $k(Y)/k$ splits in 
$k(Y,Z^*)=k(Z^*)$
into two valuations with reduced ramification exponents 
$p$ and $1$ and
therefore, by the above Alternate Corollary, 
$\operatorname{Gal}(E_{1,1}(Z),k(Y))$ is
$2$-transitive. By the above equation for $Y$ we also see 
that the valuation
$Y=0$ of $k(Y)/k$ splits in $k(Z^*)$ into the valuations 
$Z^*=0$ and
$Z^*=\infty$ with reduced ramification exponents 
$t\not\equiv 0(p)$ and
$p+1-t\not\equiv 0(p)$ respectively. Therefore on the one 
hand, by
the LCM Theorem and MRT we see that if $s\equiv 0(t(p+1-t))$
then no valuation of $k(Y)/k$, other than the valuation 
$Y=\infty$, is
ramified in a splitting field of $E_{s,a}(Z)$ over $k(Y)$, 
and hence by
Result 4 on page 841 of \cite{A3} we see that 
$\operatorname{Gal}(E_{s,a}(Z),k(Y))$ is
a quasi $p$-group, and so if we already knew that 
$\operatorname{Gal}(E_{s,a}(Z),k(Y))=A_{p+1}$ or $S_{p+1}$ 
then we would conclude that
$\operatorname{Gal}(E_{s,a}(Z),k(Y))=A_{p+1}$. On the 
other hand, by the Cycle Lemma we see 
that if $\operatorname{GCD}(p+1-t,t)=1$ then 
$\operatorname{Gal}(E_{1,1}(Z),k(Y))$ contains a $t$-cycle, 
and hence by Marggraff's First Theorem we see that if also 
$t<\frac{p+1}{2}$ 
then 
$\operatorname{Gal}(E_{1,1}(Z),k(Y)=A_{p+1}$ or $S_{p+1}$.
Note that if $p>5$ and $t$ is an odd prime factor of 
$\frac{(p-1)(p-3)}{4}$
then either $t$ is an odd prime factor of $\frac{p-1}{2}$ 
or $t$ is an odd 
prime factor of $\frac{p-3}{2}$, and in both the cases 
$1<t<\frac{p+1}{2}$
and $\operatorname{GCD}(p+1,t)=1$. Also note that, in the 
case
$p>5$, there do exist odd
prime factors of $\frac{(p-1)(p-3)}{4}$, because otherwise 
we would have
$p-1=2^m$ with $m>2$ and $p-3=2^{m'}$ with $m'>1$, and 
this would give
$2^{m'-1}=\frac{p-3}{2}=\frac{p-1-2}{2}=2^{m-1}-1$ which 
would be a
contradiction since $2^{m'-1}$ is even and $2^{m-1}-1$ is 
odd.
Thus, without using CT, we have proved the following.
\footnote{This was inspired by discussions with Walter 
Feit.}

(IV.0*) If $p>5$ then $t$ can be chosen so that 
$1<t<\frac{p+1}{2}$ 
and$\operatorname{GCD}(p+1,t)=1$, and for any such $t$, 
upon 
assuming $b=\frac{t}{t-1}$
and $s\equiv0(t(p+1-t))$, we have that no valuation of 
$k(Y)/k$, other
than the valuation $Y=\infty$, is ramified in a splitting 
field $E_{s,a}(Z)$
over $k(Y)$, and $E_{s,a}(Z)$ is irreducible in $k(Y)[Z]$ 
and has no
multiple roots in any overfield of $k(Y)$.

(IV.3*) If $p>5$ then $t$ can be chosen so that 
$1<t<\frac{p+1}{2}$
and$\operatorname{GCD}
(p+1,t)=1$, and for any such $t$, upon assuming 
$b=\frac{t}{t-1}$
and $s\equiv
 0(t(p+1-t))$, we have 
$\operatorname{Gal}(E_{s,a}(Z),k(Y))=A_{p+1}$.

The above two results with $(Y,Z)$ changed to 
$(X,Y)$, together with the 
results (IV.$0'$), (IV.$1'$), (IV.$2'$), and (IV.$4'$) 
with $(n,Y,Z)$ 
changed
to$(n+1,X,Y)$, 
\footnote{Note that with these changes, if $n+
1-p=t\not\equiv 0(p)$ then 
$\g(Y)=\frac
{(Y+n+1)^{n+1-p}-Y^{n+1-p}}{(n+1)^{2}}$.}
may be summarized in the following
%\proclaim
\subheading{Summary about the primed roof polynomial}
Let $k$ be an algebraically closed field of characteristic 
$p\ne 0$, let
$a,b$ be nonzero elements in $k$, let $n,s,t$ be positive 
integers with 
$n+1\not\equiv 0(p)$ and $n>t\not\equiv 0(p)$, and 
consider the monic 
polynomial of degree $n$ in $Y$ with coefficients in 
$k(X)$ given by
$$
\hF'_n=h(Y)(Y+b)^p-aX^{-s}Y^t\quad\text{with }0\neq b\in k,
$$
where $h(Y)$ is the monic polynomial of degree $n-p$ in 
$Y$ with coefficients 
in $k$ given by
$$
h(Y)=\frac{(Y+n+1)^{n+1-p}-Y^{n+1-p}}{(n+1)^2}.
$$
Then in the following cases $\hF'_n$ is irreducible in 
$k(Y)[Z]$, has no
multiple roots in any overfield of $k(Y)$, and gives an 
unramified covering of 
$L_k$ with the indicated Galois group 
$\hG'_n=\operatorname{Gal}(\hF'_n,k(X))$.

(IV.1) If $n+1-p=t>2\neq p$ and $b=t$ and $s\equiv 0(p-1)$ 
and $s\equiv 0(t)$, 
then $\hG'_n=A_n$.

(IV.2) If $n+1-p=t=2$ and $p\neq 7$ and $b=t$ and $s\equiv 
0(p-1)$, 
then $\hG'_n=A_n$.

(IV.3) If $n=p+1$ and $p>5$, then $t$ can be chosen so 
that $1<t<\frac{p+1}{2}$
and $\operatorname{GCD}(p+1,t)=1$, and for any such $t$, 
upon assuming
$b=\frac{t}{t-1}$ and $s\equiv 0(t(p+1-t))$, we have 
$\hG'_n=A_n$. 

(IV.4) If $n+1-p=t$ and $p=2$ and $b=t$ and $s\equiv 0(t)$, 
then $\hG'_n=S_n$.

Here CT is used only in the $p>5$ case of (IV.2).
%\endproclaim

Referring to the summaries about the tilde polynomial and 
the roof
polynomial and the primed roof polynomial, we have 
established the following
four corollaries and we have arranged a proof of the First 
and the Second 
Corollaries independent of CT.

\proclaim{First Corollary} For any $n\geq p >2$,
there exists an unramified covering of the affine line in 
characteristic
$p$ whose Galois group is $A_n$.
\endproclaim

\proclaim{Second Corollary} For any $n\geq p =2$,
there exists an unramified covering of the affine line in 
characteristic
$p$ whose Galois group is $S_n$.
\endproclaim

\proclaim{Third Corollary} Unramified coverings of the 
affine line in 
characteristic $p$ with a few more Galois groups have been 
constructed.
\endproclaim 

\proclaim{Fourth Corollary} 
Let $G$ be a quasi $p$-group. Assume that $G$ has a 
subgroup $H$ of index 
$p+1$ such that $H$ does not contain any nonidentity 
normal subgroup of $G$.
\footnote{Equivalently, $G$ is isomorphic to a transitive 
permutation group
of degree $p+1$.}
Also assume that $p$ is not a Mersenne prime and $p$ is 
different from
$11$ and $23$. Then there exists an unramified
covering of the affine line in characteristic $p$ having 
$G$ as Galois group.
\endproclaim

It only remains to note that, in view of the said three 
Summaries, the 
above Forth Corollary follows from CTT, Special CDT, and 
the Corollary of the
Fourth Irreducibility Lemma given in \S21.
Moreover, given any $n\ge p>2$,
as a definite alternative for
getting an $(A_n)$-covering as
asserted in the First Corollary without
$CT$: if $p+2<n\nequiv0(p)$ then use
(I.4) with $t=n-p$; 
if $n=p+2$ and $p<7$ then use
(I.3) with $t=2$; if $n=p+2$
and $p\ge7$ then use (IV.1) with
$t=3$ and$s=t(p-1)$; if
$n=p+1$ and $p<7$ then use (IV.2)
with $t=2$ and $s=p-1$; if $n=p+1$ and $p\ge7$
then use (IV.3) with $t=$ the smallest odd
prime factor of $\frac{(p-1)(p-3)}{2}$ and with $s=t(p+
1-t)$,
[note that in case of $n=8$ and $p=7$
this gives $t=3$ and $s=15]$; if
$n\equiv0(p)$ then use (II.1)
with $t=2$. Likewise, given any
$n\ge p=2$, as a definite alternative for
getting an $(S_n)$-covering as asserted
in the Second Corollary without $CT$:
if $n\nequiv 0(p)$ then use (I.5) with
$t=n-p$; if $n\equiv0(p)$ then use (IV.4)
with $t=n+1-p$ and $s=t$.

The above cited three summaries are transcribed from my 
e-mail message to
Serre dated 28 August 1989. This was only one out of the 
nearly a hundred
e-mail and s-mail messages which flashed back and forth 
between him and
me in the two year period September 1988 to September 
1990. Indeed it has
been a tremendous pleasure working with him. So once again 
my hearty 
MERCI MON AMI to Serre.

%\centerline{}
%\centerline{{\bf REFERENCES}}

%\let\vfill=\vfill%
%\redefine\vfill#1{\let\vfill=\vfill}
%\collect
%\enddocument

\vskip1pc

\heading Appendix by J.-P. Serre$^*$\endheading

\noindent Harvard, November 15, 1990
\vskip1pc
\noindent Dear Abhyankar,
\vskip.5pc
\par Here is my original proof that 
$\operatorname{PSL}_2(\bold{F}_q)$
occurs for the equation $Y^{q+1}-XY+1=0$.\footnote""{$^*$The author
expresses his appreciation to J. P. Serre for permission to
include the following letter in this paper.\hfill}
\par
I use ``descending Galois theory,'' i.e., I construct {\it a 
priori\/}
the Galois covering one wants. This is different from your 
``ascending''
method; in particular, I don't need any characterization of 
$\operatorname{PGL}_2(\bold{F}_q)$, or 
$\operatorname{PSL}_2(\bold{F}_q)$,
as a permutation group on $q+1$ letters.
\subheading{Notation}
\par
$p$ is a prime; $q$ is a power $p^e$ of $p$.
\par
$G=\operatorname{PGL}_2(\bold{F}_q)$, i.e.,
the quotient of $\operatorname{GL}_2(\bold{F}_q)$
by its center $\bold{F}^\ast_q$.
\par
$k$ is an algebraically closed field of characteristic $p$;
all the curves I consider are over $k$.
\subheading{Preliminary construction}
I start from the obvious fact that $G$ acts in a natural way
(by ``fractional linear transformations'') on the projective
line $\bold{P}_1$. In algebraic terms this means that $G$
acts
on $k(t)$ by $t\mapsto (at+b)/(ct+d)$. The quotient curve
$\bold{P}'_1=\bold{P}_1/G$ is of course (L\"uroth's theorem!)
a projective line. Equivalently, the field of invariants of 
$G$
in $k(t)$ is a purely transcendental field $k(x)$.
\par
The first computational problem which arises is to write $x$
explicitly. To do so, let us call $(u,v)$ the homogeneous 
coordinates
on $\bold{P}_1$, so that $t=v/u$. The invariant theory of 
$k[u,v]$
with respect to the action of $G$ has been done long ago by 
Dickson. The
basic covariants are the following homogeneous polynomials:
$$
A(u,v)=uv^q-vu^q,\qquad
 B(u,v)=(uv^{q^2}-vu^{q^2})/A(u,v).
$$
They are of degree $q+1$ and $q^2-q$ respectively. Hence the
ratio
$$
x=B(u,v)^{q+1}/A(u,v)^{q^2-q}
$$
is invariant by $G$. Its expression in terms of $t$ is easy 
to find:
if we write $A(u,v)=u^{q+1}a(t)$, $B(u,v)=u^{q^2-q}b(t)$, we
 have
$$
a(t)=t^q-t,\qquad b(t)=(t^{q^2}-t)/(t^q-t)=a(t)^{q-1}+1
$$
and
$$
x=b(t)^{q+1}/a(t)^{q^2-q}=(a(t)^{q-1}+1)^{q+1}/a(t)^{q^2-q}.
\tag"$(\ast)$"
$$
This shows that $x$ is a rational function of $t$ of degree 
$q(q^2-1)$.
Since it is invariant by $G$, which has order $q(q^2-1)$, 
Galois theory shows
it {\it generates\/} the field of the $G$-invariant elements
of $k(t)$.
Hence we have found our parameter for 
$\bold{P}'_1=\bold{P}_1/G$.
\subheading{Ramification}
It is necessary to study the ramification in the Galois 
extension
$k(t)/k(x)$ thus constructed. This amounts to looking for the
fixed
points of the action of $G$ on the projective line 
$\bold{P}_1$. It is easy
to see that these fixed points make up two orbits. Namely:
\par
(a) The $\bold{F}_q$-rational points of $\bold{P}_1$. This 
orbit has $q+1$
elements. The stabilizer of an element is a triangular 
subgroup (``Borel
subgroup'') of order $q(q-1)$. Since that order is divisible 
by $p$, there
is wild ramification.
\par
The point of $\bold{P}^\prime_1$ corresponding to this orbit 
is
$x=\infty$.
\par
(b) The ``quadratic'' points, i.e., the 
$\bold{F}_{q^2}$-rational points of
$\bold{P}_1$ which are not rational over $\bold{F}_q$. There 
are $q^2-q$
of them. The stabilizer of such a point is a cyclic group of 
order $q+1$
(``nonsplit Cartan subgroup''). Since that order is prime to 
$p$, the
ramification at such a point is tame.
\par
The point of $\bold{P}^\prime_1$ corresponding to this orbit 
is
$x=0$.
\par
Hence we see that the covering of $\bold{P}^\prime_1$ we get 
in this way is
{\it ramified both at\/} $0$ {\it and\/} $\infty$, {\it and 
nowhere else\/}.
The next step is thus:
\subheading{Getting rid of the ramification at $0$ using 
Abhyankar's lemma}
We consider the cyclic extension $k(X)$ of $k(x)$ defined by 
the
equation $X^{q+1}=x$. By making a base change to that 
extension
(i.e., by considering $k(t,X)/k(X))$ we get rid of the 
ramification
at $0$. Only the ramification at $\infty$ remains. Of course, 
one has
to see what the new Galois group is. There are two cases:
\par
(i) $p=2$. The extensions $k(X)/k(x)$
and $k(t)/k(x)$ are disjoint. Hence the new
Galois group is equal to the old one, namely
$G=\operatorname{PGL}_2(\bold{F}_q)$, which
happens to be equal to $\operatorname{PSL}_2
(\bold{F}_q)$.
\par
(ii) $p\ne2$. The extensions $k(X)/k(x)$ and
$k(t)/k(x)$ have a quadratic extension in common,
namely $k(x^{1/2})$. Hence the new Galois group
is $G'=\operatorname{PSL}_2(\bold{F}_q)$.
\par
In both cases, one thus gets {\it a Galois extension of
$k(X)$ with Galois group $\operatorname{PSL}_2(\bold{F}_q)$
which is ramified only at $X=\infty$}.
\par
It remains to see that this extension is the same as the one
you get by the equation $Y^{q+1}-XY+1=0$.
\subheading{An equation for the degrees $q+1$ extension}
I go back to the $k(t)/k(x)$ extension with Galois group 
$G$.
Let $H$ be the triangular subgroup of $G$, of index $q+1$.
The fixed field $k(t)^H$ is an extension of degree $q+1$
of $k(x)$ and we want to find a generator $y$ for that field 
(which
is also a purely transcendental field, of course). We may 
assume
that $H$ is the group of transformations $t\mapsto at+b$. 
This 
shows that the polynomial
$$
y=a(t)^{q-1}=(t^q-t)^{q-1}
$$
is invariant by $H$. Since its degree is $q(q-1)=|H|$, the 
same
argument as above shows that $k(y)$ is equal to $k(t)^H$.
\par
I now write the equation of degree $q+1$ relating $y$ to $x$. 
This
is easy, since by construction, we have $x=(y+1)^{q+1}/y^q$,
cf.\ $(\ast)$ above. We thus get the equation
$$
(y+1)^{q+1}-xy^q=0.
$$
\par
But I want to work on $k(X)$, with $x=X^{q+1}$. We have:
$$
(y+1)^{q+1}-X^{q+1}y^q=0.
$$
Let me put $Y=(y+1)/yX$, i.e., $y=1/(XY-1)$. The above 
equation
becomes
$$
Y^{q+1}-XY+1=0,
$$
and we are done.
\par
This is the proof I found in 1988
when I started thinking about your
problem. The first part is natural
enough---and could indeed be applied
to other groups, \`a la Nori. The second
part (the search for the degree $q+1$
equation) is not; it looks like a happy
coincidence, and I would not have found
it if I had not known in advance your
polynomial $Y^{q+1}-XY+1$.
\par
With best regards,
\par
\quad Yours
$$
\text{J-P. Serre}
$$
\par
\noindent PS\<$\quad$\<The determination of the
invariants of $G$ in $k(t)$ is not new.
I am almost certain to have seen it in print very
long ago, as an elementary exercise in Galois
theory. (Indeed: see Lang's {\it Algebra\/},
2nd ed., p.\ 349, exercise 33, and also P. Rivoire,
{\it Ann. Inst. Fourier\/} {\bf6} (1955--1956),
pp.\ 121--124.)


\Refs

\ref \key A1
\by S. S. Abhyankar
\paper On the ramification of algebraic functions
\jour Amer. J. Math.
\vol 77
\yr 1955
\pages 572--592
\endref

\ref \key A2
\bysame
\paper Local uniformization on algebraic surfaces over 
ground fields
of characteristic $p\ne 0$
\jour Ann. of Math.
\vol 63
\yr 1956
\pages 491--526
\endref

\ref \key A3
\bysame
\paper Coverings of algebraic curves
\jour Amer. J. Math.
\vol 79
\yr 1957
\pages 825--856
\endref

\ref \key A4
\bysame
\paper On the ramification of algebraic functions part 
{\rm II}\RM: Unaffected
equations for characteristic two
\jour Trans. Amer. Math. Soc.
\vol 89 \rm(1958)
\year 1958
\pages 310--324
\endref

\ref \key A5
\bysame
\paper Tame coverings and fundamental groups of algebraic 
varieties,
Part {\rm III:} Some other sets of conditions for the 
fundamental group to
be abelian
\jour Amer. J. Math.
\vol 82
\yr 1960
\pages 179--190
\endref

\ref \key A6
\bysame
\book Algebraic geometry for scientists and engineers
\publ Amer. Math. Soc., Providence, RI
\yr 1990
\endref

\ref \key A7
\bysame
\paper Square-root parametrization of plane curves \rm(to 
appear)
\endref

\ref \key As
\by M. Aschbacher
\book Finite groups
\publ Cambridge University Press, Cambridge
\yr 1986
\endref

\ref \key B
\by R. Brauer
\paper On the structure of groups of finite order
\jour Proc. Internat. Congress Math.
\vol 1
\yr 1954
\pages 209--217
\endref

\ref \key Bu
\by W. Burnside
\book Theory of groups of finite order
\publ Cambridge University Press, Cambridge
\yr 1911
\endref
 
\ref \key BP
\by W. S. Burnside and A. W. Panton
\book Theory of equations
\publ Vols. I and II, Dublin, Hodges, Figgs and Co., London
\yr 1904
\endref
 
\ref \key C
\by P. J. Cameron
\paper Finite permutation groups and finite simple groups
\jour Bull. London Math. Soc.
\vol 13
\yr 1981
\pages 1--22
\endref

\ref \key Ca
\by R. D. Carmichael
\book Groups of finite order
\publ Dover, New York
\yr 1956
\endref
 
\ref \key Ch
\by C. Chevalley
\paper Sur certains groupes simples
\jour T\^ohoku Math. J.
\vol 7
\yr 1955
\pages 14--66
\endref

\ref \key CKS
\by C. W. Curtis, W. M. Kantor, and G. M. Seitz
\paper The \RM2-transitive permutation representations of 
the finite
Chevalley groups 
\jour Trans. Amer. Math. Soc.
\vol 218
\yr 1976
\pages 1--59
\endref

\ref \key D1
\by L. E. Dickson
\book Linear groups with an exposition of the Galois field 
theory
\publ Teubner, Leipzig
\yr 1958
\endref

\ref \key D2
\bysame
\paper On finite algebras 
\jour Nachr. Akad. Wiss. G\"ottingen 
\vol 
\yr 1905
\pages 358--393 
\endref

\ref \key F
\by W. Feit
\paper On a class of doubly transitive permutation groups
\jour Illinois J. Math.
\vol 4
\yr 1960
\pages 170--186
\endref

\ref \key Fr
\by F. G. Frobenius
\paper \"Uber aufl\"osbare Gruppen \rm IV
\jour Sitzungsberichte der K\"oniglich Preu\ss ischen 
Akademie der Wissenchaften zu Berlin
\vol 
\yr 1901
\pages 1216--1230; reprinted as No. 63 on pages 189-203 of 
Band III
of the ``Gesammelte Abhandlungen'' of Frobenius published 
by Springer-Verlag,
New York, 1968 
\endref

\ref \key FT
\by W. Feit and J. Thompson
\paper Solvability of groups of odd order
\jour Pacific J. Math.
\vol 13
\yr 1963
\pages 775--1029
\endref

\ref \key G1
\by D. Gorenstein
\book Finite groups
\publ Chelsea Publishing Company, New York
\yr 1980
\endref

\ref \key G2
\bysame
\book Finite simple groups
\publ Plenum Press, New York
\yr 1983
\endref

\ref \key G3
\bysame
\paper Classifying the finite simple groups
\jour Bull. Amer. Math. Soc.
\vol 14
\yr 1986
\pages 1--98
\endref

\ref \key G4
\bysame
\paper The classification of the finite simple groups, a 
personal 
journey\,\RM:
The early years\,, {\rm A Century of Mathematics in 
America, Part I}
\jour Hist. of Math., vol. 1, Amer. Math. Soc., 
Providence, RI,
1988, pp. 447--476
\endref

\ref \key H
\by D. Hilbert
\paper \"Uber die Irreduzibilit\"at ganzzahligen 
Koefficienten 
\jour Crelle J.
\vol 110
\yr 1892
\pages 104--129
\endref

\ref \key HB
\by B. Huppert and N. Blackburn
\book Finite groups \rm I, II, III
\publ Springer-Verlag, New York
\yr 1982
\endref

\ref \key J1
\by C. Jordan
\book Trait\'e des substitutions et des \'equations 
alg\'ebriques
\publ Gauthier-Vill., Paris
\yr 1870
\endref

\ref \key J2
\bysame
\paper Th\'eor\`emes sur les groupes primitifs
\jour J. Math. Pures Appl.
\vol 16
\pages 383--408
\endref

\ref \key J3
\bysame
\paper Sur la limite de transitivit\'e des groupes non 
altern\'es
\jour Bull. Soc. Math. France
\vol 1
\year 1873
\pages 40--71
\endref

\ref \key K1
\by W. M. Kantor
\paper Rank \RM3 characterizations of classical geometries
\jour J. Algebra
\vol 36
\yr 1975
\pages 309--313 
\endref

\ref \key K2
\bysame
\paper Homogeneous designs and geometric lattices
\jour J. Combin. Theory Ser. A
\vol 38
\yr 1985
\pages 66--74 
\endref

\ref \key KL
\by W. M. Kantor and R. L. Liebler
\paper The rank \RM3 permutation representations of the 
finite classical
groups
\jour Trans. Amer. Math. Soc.
\vol 271
\yr 1982
\pages 1--71 
\endref

\ref \key L
\by M. W. Liebeck
\paper The affine permutation groups of rank three
\jour Proc. London Math. Soc.
\vol 54
\yr 1987
\pages 477--516 
\endref

\ref \key Mar
\by B. Marggraff
\book \"Uber primitive Gruppen mit transitiven
Untergruppen geringeren Grades
\publ Dissertation, Giessen
\yr 1892
\endref

\ref \key Mas
\by G. Mason
\paper The classification of quasi-thin simple groups. 
\rm(To appear)
\endref

\ref \key Mat
\by E. Mathieu
\paper M\'emoire sur l'\'etude des fonctions de plusieurs 
quantit\'es,
sur la mani\`ere de les former, et sur les substitutions 
qui les laissent 
invariables
\jour J. Math  Pures Appl.,
\vol 18
\yr 1861
\pages 241--323
\endref

\ref \key Mo
\by E. H. Moore
\paper International Mathematics Congress, Chicago \rm1893
\jour New York,
1896, p. 210
\endref

\ref \key N
\by P. M. Neumann
\paper Some primitive permutation groups
\jour Proc. London Math. Soc.
\vol 50
\yr 1985
\pages 265--281
\endref

\ref \key O
\by M. E. O'Nan
\paper Normal structure of the one-point stabilizer of a 
doubly
transitive permutation group {\rm I} and {\rm II}
\jour Trans. Amer. Math. Soc.
\vol 214
\yr 1975
\pages  1--74
\endref

\ref \key R
\by R. Ree
\paper Sur une famille de groupes de permutations 
doublement transitifs
\jour Canad. J. Math.
\vol 16
\yr 1964
\pages 797--820
\endref

\ref \key S1
\by J-P. Serre
\paper G\'eom\'etrie alg\'ebrique et g\'eom\'etrie 
analytique
\jour Ann. Inst. Fourier
\vol 6
\yr 1956
\pages 1--42
\endref

\ref \key S2
\bysame
\book Algebraic Groups and Class Fields
\publ Springer-Verlag, New York
\yr 1988
\endref

\ref \key St
\by R. Steinberg
\paper Variations on a theme of Chevalley
\jour Pacific J. Math.
\vol 9
\yr 1959
\pages 875--891
\endref

\ref \key Su1
\by M. Suzuki 
\paper On a class of doubly transitive groups
\jour Ann. Math.
\vol 75
\yr 1962
\pages 104--145
\endref

\ref \key Su2
\bysame
\book Group Theory {\rm I, II}
\publ Springer-Verlag, New York
\yr 1986
\endref

\ref \key V
\by B. L. van der Waerden
\book Modern Algebra \rm I, II
\publ Frederick Ungar Publishing Co., New York
\yr 1949
\endref

\ref \key W
\by J. H. M. Wedderburn
\paper A theorem on finite algebras
\jour Trans. Amer. Math. Soc.
\vol 6 \rm(1905)
\pages 349--352
\endref

\ref \key Wi
\by H. Wielandt
\book Finite Permutation Groups
\publ Academic Press, New York
\yr 1964
\endref

\ref \key Z1
\by H. Zassenhaus
\paper Kennzeichnung endlicher linearer Gruppen als 
Permutationsgruppen
\jour Abh. Math. Sem. Univ. Hamburg
\vol 11
\yr 1936
\pages 17--44
\endref

\ref \key Z2
\bysame
\paper \"Uber endliche Fastk\"orper
\jour Abh. Math. Sem. Univ. Hamburg
\vol 11
\year 1936
\pages 187--220
\endref
\endRefs
\enddocument